\documentclass[msom,nonblindrev]{informs3}

\OneAndAHalfSpacedXI 

\usepackage{comment}
\usepackage{tikz}
\usepackage{algorithm}
\usepackage{algorithmic}
\usepackage{subfiles}
\usepackage{subfigure}
\usepackage{enumitem}
\usepackage{hyperref}
\usepackage{multirow}
\usepackage{bbm}

\usepackage{natbib}
 \bibpunct[, ]{(}{)}{,}{a}{}{,}%
 \def\bibfont{\small}%
\TheoremsNumberedThrough     

\EquationsNumberedThrough    

\MANUSCRIPTNO{0000} 

\begin{document}



\RUNAUTHOR{Sun, Dai, and Shi}

\RUNTITLE{Inpatient Overflow with PPO}

\TITLE{Inpatient Overflow Management\\ with Proximal Policy Optimization}

\ARTICLEAUTHORS{%
\AUTHOR{Jingjing Sun\thanks{Code is available at \url{https://github.com/OverflowPPO/InpatientOverflow.git}.}}
\AFF{School of Data Science, The Chinese University of Hong Kong, Shenzhen, China, \EMAIL{jingjingsun@link.cuhk.edu.cn}, \URL{}}

\AUTHOR{J. G. Dai}
\AFF{School of Operations Research and Information Engineering, Cornell University, Ithaca, NY 14853, \EMAIL{jd694@cornell.edu}, \URL{}}

\AUTHOR{Pengyi Shi}
\AFF{Mitch Daniels School of Business, Purdue University, West Lafayette, IN 47907, \EMAIL{shi178@purdue.edu}, \URL{}}

} 

\ABSTRACT{%
\textbf{Problem Definition}:
Managing inpatient flow in large hospital systems is challenging due to the complexity of assigning randomly arriving patients -- either waiting for primary units or being overflowed to alternative units. Current practices rely on ad-hoc rules, while prior analytical approaches struggle with the intractably large state and action spaces inherent in patient-unit matching. Scalable decision support is needed to optimize overflow management while accounting for time-periodic fluctuations in patient flow.
\textbf{Methodology/Results}:
We develop a scalable decision-making framework using Proximal Policy Optimization (PPO) to optimize overflow decisions in a time-periodic, long-run average cost setting. To address the combinatorial complexity, we introduce atomic actions, which decompose multi-patient routing into sequential assignments. We further enhance computational efficiency through a partially-shared policy network designed to balance parameter sharing with time-specific policy adaptations, and a queueing-informed value function approximation to improve policy evaluation. Our method significantly reduces the need for extensive simulation data, a common limitation in reinforcement learning applications. Case studies on hospital systems with up to twenty patient classes and twenty wards demonstrate that our approach matches or outperforms existing benchmarks, including approximate dynamic programming, which is computationally infeasible beyond five wards.
\textbf{Managerial Implications}:
Our framework offers a scalable, efficient, and explainable solution for managing patient flow in complex hospital systems. More broadly, our results highlight that domain-aware adaptation is more critical to improving algorithm performance than fine-tuning neural network parameters when applying general-purpose algorithms to specific applications. 
}%



\maketitle

%


\section{Introduction}

Managing inpatient flow is a significant challenge in hospital operations, especially when patient demand fluctuates unpredictably. Hospital capacity often cannot be scaled up quickly enough when demand increases, resulting in intermittent bed shortages and overcrowding. Overcrowded inpatient units can lead to a range of problems, including denied admissions, emergency department (ED) boarding, ambulance diversions, surgery cancellations, and adverse patient outcomes caused by long wait times~\citep{mckenna2019emergency,jones2022association}. The issue of intermittent bed shortages is worsened by the common practice of grouping specific clinical services within designated inpatient units, each with a fixed number of beds~\citep{song2020capacity}. When the designated inpatient units of certain clinical services reach capacity, hospitals often have to find alternative ways for their incoming patients, which frequently involves assigning them to non-dedicated ``overflow'' units.

This overflow practice, also known as off-service placement~\citep{song2020capacity} or patient outlying~\citep{sprivulis2006association}, helps temporarily alleviate congestion in hospitals by assigning patients to non-primary units when their designated wards are full. This practice is commonly used worldwide, including in the United States~\citep{zhang2020simulation, lim2024spillover}, England~\citep{stylianou2017medical}, France~\citep{stowell2013hospital}, Spain~\citep{alameda2009clinical}, Australia~\citep{sprivulis2006association}, and Singapore~\citep{shi2016models}, with 10\% to 20\% of patients placed in overflow units according to these studies.
However, this practice introduces additional complexity into patient placement and may compromise the quality and continuity of care~\citep{stowell2013hospital, stylianou2017medical, song2020capacity}. Moreover, excessive reliance on overflow can reduce overall system efficiency by increasing coordination efforts and potentially triggering a cycle of subsequent overflow events, further exacerbating congestion~\citep{dong2019off}.

Effectively managing overflow requires nuanced decision-making to balance the benefits of reducing inpatient overcrowding with the costs of unnecessary or excessive overflow assignments. It is crucial to account for the long-term impact, given that the typical patient length-of-stay is 3 to 4 days. This means that an overflow decision not only affects the immediate system state but can also influence system congestion over several days. This complexity is further compounded by \emph{time-varying} dynamics, such as predictable hourly patterns of patient arrivals and discharges. For example, should patients wait a few more hours if more discharges are expected soon? Overflow assignments in practice, though widely used, often rely on anecdotal experiences or ad-hoc rules, which are suboptimal and lack the rigorous analytical support needed for effective decision-making. \cite{dai2019inpatient} is one of the earliest works that formulate the overflow decision problem as a periodic, long-run average-cost Markov Decision Process (MDP), aiding the decisions on whether to overflow a waiting patient and, if so, to which ward. Several other recent works have also adopted analytical methods to aid overflow decision making~\citep{schmidt2013decision,zhang2020simulation,bertsimas2024hospital,dean2024learning}.

These analytical approaches, however, share a common challenge: the intractably large state and action spaces inherent in the ``matching problems.'' In the context of patient overflow, this arises from multiple patients (customers) from different specialties waiting for admission, while a combinatorial number of bed assignment (server) options -- including the choice to continue waiting -- are available. Similar challenges are observed in other domains like ride-hailing, inventory management, and volunteer assignment, where the combinatorial nature of potential actions greatly complicates decision-making. \cite{dai2019inpatient} tackle the large state space using value function approximation but encounter difficulties with the large action space, which limited their method's scalability  -- their case study only dealt with systems involving up to five specialties and wards. Other approaches, such as Linear Programming or fluid-based optimization, manage large action spaces more efficiently but often oversimplify future stochastic variations~\citep{bertsimas2024hospital}; see a comprehensive literature review in Section~\ref{sec:lit-review}. 
These limitations highlight the broader difficulty in developing practical and scalable decision-support algorithms for complex customer-server matching problems, not just to healthcare but a wide range of different domains.

\subsection{Paper Overview and Contribution}

This paper addresses the existing gap by developing an algorithmic framework based on the Proximal Policy Optimization (PPO) method. Our key contributions are summarized as follows.

\begin{enumerate}[leftmargin=*]
\item \textbf{Scalable Algorithim}:
We develop the first known scalable algorithm for overflow decision-making in a time-periodic, long-run average setting. Our algorithm successfully handles systems with up to 20 medical specialties and 20 hospital units, greatly exceeding the scale of prior methods. Beyond healthcare, our framework is applicable to a broad class of matching problems in operations management with similar queueing structures and contributes to efficiently managing the combinatorial complexity in customer-server matching across domains.

The success of this framework hinges on two key elements: (a) leveraging the ``atomizable'' action space and (b) incorporating domain knowledge of queueing systems. 

\item  \textbf{Action Atomization}: 
As detailed in Section~\ref{Sec: Atomic actions}, atomic action decomposes class-level decisions into sequential individual-level assignments, which is critical for simplifying the original combinatorial action space. Despite its seemingly myopic nature, we show that this approach, when combined with PPO, achieves performance comparable to the ADP benchmark in a five-pool system, the largest setting considered in~\cite{dai2019inpatient}; i.e., this approach \emph{works}. 
The atomic formulation is inspired by \cite{feng2021scalable}, who coined the term \emph{atomic action decomposition} in the context of finite-horizon ride-hailing problems. 
Extending this idea to the more challenging infinite-horizon, long-run average cost setting is nontrivial, and to the best of our knowledge, this paper is the first to show that combining atomic actions with PPO is effective in this regime. We discuss related developments in this growing research area in Section~\ref{sec:lit-review}.

\item \textbf{Dealing with Time-periodicity}: 
To account for real-world daily and weekly periodic patterns in patient arrivals and discharges, we integrate time-periodicity into our PPO framework. First, we establish a theoretical guarantee for policy improvement of PPO under the periodic setting (Section~\ref{sec:ppo-method}). Second, we design a partially-shared policy network that balances parameter sharing and temporal specificity (Section~\ref{Subsec: Policy representation}). 
Lastly, we propose a new value function approximation tailored to the queueing dynamics under the atomic actions to improve computational efficiency (Section~\ref{Subsec: policy evaluation}).

\item \textbf{Actionable Insights}:
Case study results in Section~\ref{sec:case-study}, including ablation studies, confirm that our integrated framework enhances sample efficiency while maintaining superior performance over other benchmarks. These results highlight that domain-aware algorithm design is more impactful than tuning generic hyperparameters. Moreover, our policy visualization and interpretability analysis reveal that the learned policies are both intuitive (aligning with managerial expectations) and insightful (uncovering nuanced network effects that are difficult to capture with heuristic rules). These insights provide actionable guidelines for managing overflow while also highlight the necessity of sophisticated, data-driven tools in complex hospital networks.
  
\end{enumerate}


\section{Literature review}
\label{sec:lit-review}

\subsection{Bed Partition and Patient Routing}
 
Adding new beds to alleviate inpatient congestion is often not feasible due to constraints related to space, finances, and regulations. Therefore, we focus on literature that explores strategies for utilizing existing bed resources more efficiently~\citep{Zamani2024}. \cite{best2015managing} examine bed management at a strategic level, where the decision involves partitioning inpatient beds into separate wings dedicated to different specialties (overflow is not allowed). As a follow-up, \cite{izady2021clustered} consider a similar setting but introduce the option of ``pooled wings,'' which can accommodate patients from multiple specialties when dedicated wings reach capacity.

\cite{song2020capacity} highlight that wing designs can exacerbate inpatient crowding, particularly when combined with unpredictable patient volumes and limited bed capacity. This often results in frequent off-service (overflow) placements, especially in large teaching or public hospitals that handle high volumes of unpredictable emergency department demand~\citep{manning2023systematic}. However, it is crucial to manage these overflow assignments carefully to minimize adverse patient outcomes and ensure quality care~\citep{lim2024spillover}. 
Regarding patient assignment under fixed bed partitioning, the most relevant study is \cite{dai2019inpatient}, which serves as the foundation for this paper. They formulate the problem as an MDP and use approximate dynamic programming (ADP) to address the large state space. Their approach focuses on a deterministic policy space, which can handle systems with up to five server-pools but does not address the large action space. \cite{chen2021primal} develop a primal-dual method and use a similar MDP for inpatient assignments as one of their case studies. Their method tackles the large state space by decomposing the problem into multiple subproblems with an inverted-V structure (i.e., single-class multi-pool). However, the action space in each subproblem remains combinatorial. 
In contrast, the PPO algorithm we propose in this paper addresses both the large state and action spaces and can effectively solve problems in large systems (e.g., 20 server pools). 
\citet{dong2025multiclass} use ADP to optimize scheduling policies in multiclass, multiserver queues with wait-dependent service slowdowns, focusing on leveraging the problem structure to estimate value function differences directly rather than addressing action space complexity.
Other recent studies applying ADP methods to patient-provider matching in healthcare operations include~\cite{cire2022dynamic, gao2024shortening, khorasanian2024dynamic}, among others.

In contrast to studying the routing problem from a sequential (dynamic) decision-making perspective, \cite{thomas2013automated} formulate the patient assignment problem as a mixed-integer programming (MIP) model. Their objective is to minimize immediate costs without considering the long-term impact of current actions, which results in a myopic policy. Similarly, \cite{bertsimas2024hospital} propose a robust multistage optimization framework that minimizes the current cost along with an approximation of future costs based on a deterministic model. While their approach provides flexibility for individualized policies, it follows a different methodological framework. In particular, future costs are handled differently -- they are approximated as a weighted sum over the remainder of the day and the following week, based on the \emph{worst-case} scenarios within a predefined uncertainty set. This can lead to more conservative, risk-averse policies. In contrast, our approach directly optimizes policy decisions under uncertainty via PPO, focusing on the stochastic nature of the problem. Additionally, this paper emphasizes the importance of tailoring a general-purpose PPO algorithm to the structure of the inpatient overflow problem (e.g., specialized policy NN designs and basis functions) and demonstrates that these modifications can have a greater impact on performance than hyperparameter tuning.

\subsection{PPO Applications}
The PPO algorithm was first introduced in the seminal work by \cite{schulman2017proximal}, which provided theoretical justification for using PPO to solve MDP problems with bounded cost functions and discounted cost objectives. \cite{dai2022queueing} extend these theoretical guarantees to MDP problems with unbounded cost functions and long-run average cost objectives. Since then, numerous studies have successfully applied PPO to domains such as inventory management \citep{perez2021algorithmic} and transshipment \citep{vanvuchelen2024cluster}. 
However, our inpatient overflow MDP problem presents unique challenges that make the direct application of PPO algorithms difficult, specifically due to the large action space and the time-varying, periodic long-run average setting.

In conventional PPO algorithms, the policy space is directly parameterized by a Neural Network (NN), with the output size corresponding to the number of possible actions, typically a small number (e.g., \cite{schulman2017proximal, dai2022queueing}). However, in our problem, the action space is significantly larger, making it computationally prohibitive to train the policy network using traditional parameterization methods.
\citet{feng2021scalable} demonstrates the effectiveness of combining action decomposition with PPO for ride-hailing, but only in a finite-horizon setting.
While PPO has been successfully applied to infinite-horizon MDPs \citep{dai2022queueing}, to our knowledge, there are no established examples of its application to problems combining both large action spaces and infinite-horizon objectives.  Our work addresses both challenges simultaneously, handling large action spaces in a time-varying, infinite-horizon environment. The study of atomic-action-based PPO is actively growing. While this paper was under review, concurrent work by \cite{Huo2025Thesis} and \cite{dai2025atomic} have demonstrated the effectiveness of this framework in data switch scheduling, and in battery charging and electric vehicle scheduling for robo-taxi systems, respectively; the latter introduced the term \emph{Atomic-PPO}. Motivated by the empirical success, \cite{dai2025optimal} further provides a theoretical justification for atomic action decomposition.


\section{Model Formulation}
\label{sec:model-formulation}

In this section, we present the queueing model formulation abstracted from the overflow problem. While we describe a general queueing model, we ensure it retains the key features relevant to our application. Specifically, the model is based on a parallel server system with multiple customer classes and server pools. Each customer class represents a medical specialty (or a cluster of related specialties), and each server pool corresponds to an inpatient unit (ward), with servers representing individual inpatient beds. 
We detail the queueing discipline along with the arrival and departure processes in Section~\ref{Subsec: queue}. We formulate the overflow decision problem as an MDP in Section~\ref{Sec:MDP model}. 
Appendix~\ref{app:notation} summarizes all the notations.

\subsection{Parallel-server System}
\label{Subsec: queue}

Consider a queueing system with $J$ customer classes and $J$ parallel-server pools. Each pool has $N_j$ identical servers and is dedicated to serving one of the $J$ customer classes -- the primary pool for that class.  Without loss of generality, we assume that server-pool $j$ is the primary pool for class $j$ customers, $j\in \mathcal{J}=\{1,2,\dots,J\}$. Upon the arrival of a class $j$ customer, if there are idle servers in primary pool $j$, the customer is admitted immediately, known as a primary assignment. Otherwise, the customer may wait in buffer $j$ (which has no capacity limit) for service in the primary pool or, at a decision epoch, be assigned to a non-primary pool according to the overflow policy specified in Section~\ref{Sec:MDP model}. This latter scenario is called an overflow assignment, and the non-primary pool that can serve class $j$ customers is referred to as an overflow pool for class $j$.

\noindent\textbf{Arrival process.} The underlying system evolves on a continuous-time scale. We assume that the arrival process for each customer class follows a periodic Poisson process, with the associated arrival rate function $\lambda_j(t)$ being periodic with a period of $T$ for each $j\in\mathcal{J}$. For ease of exposition, we set $T=1$ day to capture the time-varying pattern within a day. The corresponding daily arrival rate for class $j$ is given by $\Lambda_j = \int_0^1 \lambda_j(t) dt$.

\noindent\textbf{Departure process. }
We adopt a \emph{two-time-scale} departure scheme as in \cite{dai2019inpatient}. For customers admitted to pool $j$, there is a probability $\mu_j$ that they will depart from the system on the current day, and a probability $1-\mu_j$ that they will stay (memoryless). If a customer departs, the time of departure within the current day is modeled by the random variable $h_\textrm{dis}$, which follows a general distribution with support on $[0,1)$ and cumulative distribution function (CDF) $F_j$. This scheme was first introduced in \cite{shi2016models} and has since been widely used in modeling inpatient flow, e.g., \cite{dai2019inpatient, bertsimas2024hospital, kim2015icu}.

\subsection{MDP Model}
\label{Sec:MDP model}

We use an MDP framework to decide on whether to assign a waiting customer to one of the overflow pools. The objective is to find the optimal overflow policy that minimizes the infinite-horizon, long-run average cost. We consider a discrete-time setup, where the decision maker (e.g., hospital manager) makes overflow decisions at $m$ predefined time epochs each day (e.g., hourly), denoted as $t_0, t_1, \dots, t_{m-1}$ for day $t$.

\paragraph{State.}
We follow \cite{dai2019inpatient} to define the state space. Denote the system state observed at decision epoch $t_k$, before taking any action, as $S(t_k)$, where 
$$ S(t_k)=\big( X_1(t_k),\dots,X_J(t_k),Y_1(t_k),\dots,Y_J(t_k), h(t_k) \big). $$ 
The state variable is a $(2J+1)$-dimensional vector, with each component as follows: 
 \begin{itemize}
     \item $X_j(t_k)$, $j=1, \dots, J$, denotes the \textit{customer count} for pool $j$. This count includes both the number of customers waiting in buffer $j$ and those in service in pool $j$; we show below on how to recover each part. The former (waiting count) only includes class $j$ customers, while the latter (in-service count) may include customers from other classes due to overflow assignments. 
     
     \item $Y_j(t_k)$, $j=1, \dots, J$, denotes the \textit{to-depart count} in pool $j$. This quantity represents the number of customers anticipated to depart from pool $j$ between the present time and the end of the current day. Tracking this quantity is essential for the two-time-scale inpatient departure scheme, where the ``planned discharge'' information is used to inform overflow decisions. 
     
     \item $h(t_k)=k\in\{0,1,2,\dots,m-1\}$ denotes the \textit{time-of-day} indicator. 
 \end{itemize} 
 We denote the state space as $\mathcal{S}$, and the sub-space that contains all states with time-of-day being $h$ as $\mathcal{S}^h$. Therefore, $(\mathcal{S}^0,\dots,\mathcal{S}^{m-1})$ is a partition for $\mathcal{S}$ without overlap between sub-spaces.

\paragraph{Action.}
We assume that all primary assignments are made immediately when customers arrive or when primary servers become available. Thus, at each epoch $t_k$, the decision maker only needs to decide whether to assign a waiting customer to an overflow pool or let them continue waiting. We use $f(t_k) = \{f_{i,j}(t_k), ~ i,j=1,\dots,J\}$ to represent the overflow/wait decision at time $t_k$. Each $f_{i,j}$ corresponds to the number of assignments from class $i$ to pool $j$. If $i \neq j$, it denotes the number of overflow assignments; if  $i = j$, it denotes the number of customers who remain waiting.

For a given state $S(t_k) = s$, action $f$ is called \emph{feasible} if it satisfies the following conditions: 
\begin{equation}
    \begin{aligned}
    \sum_{ \ell=1,\ell\ne j}^J f_{\ell,j}\le  (N_j-x_j)\vee 0 , 
\quad 
     \sum_{ \ell=1}^J f_{i,\ell}=  (x_i-N_i)\vee 0  , \quad i,j=1,\dots,J. 
\end{aligned}
\label{Equ: feasibility of action}
\end{equation}
To explain, the total number of overflow assignments from all non-primary classes to pool $j$ must not exceed the number of idle servers in pool $j$, which can be computed as $(N_j - x_j) \vee 0$, with $x \vee 0 = \max(x, 0)$ for any real number $x$. Similarly, the total number of assignments for class $i$, including overflow assignments to a pool $\ell \neq i$ and those kept waiting, must equal the number of waiting customers. This can be calculated as $(x_i - N_i) \vee 0$, which we also refer to as the queue length. We use $q_i = (x_i - N_i) \vee 0$ to denote this queue length for class $i$ (buffer $i$). 
Overflow decisions are often coordinated across patients and units rather than made purely myopically, as observed at our partner hospital. Prior
work has also considered delaying placement (e.g., via holding) to improve assignment quality~\citep{zhang2020simulation}. More broadly, delaying decisions can improve decision quality in operational contexts by allowing information to accumulate~\citep{xie2025benefits}. We use the decision-epoch structure as a modeling abstraction of this coordination process. Although assignment to not-yet-idle beds is not allowed, future availability
is anticipated through pending discharges, which allows forward-looking decisions that reflect cross-unit coordination. 

\paragraph{State transition.}
We denote the state at epoch $t_k$ before taking an action as $s$, i.e., the pre-action state.  
After an action is taken, the resulting state is denoted as $s^+$, i.e., the post-action state. We further denote by $s'$ the state after the exogenous events (random arrivals and departures) occurring between the current epoch $t_k$ and the next epoch $t_{k+1}$. The transition from $s$ to $s^+$ is the same across all epochs (adding overflow-in patients and subtracting overflow-out patients), whereas the transition from $s^+$ to $s'$ differs between the midnight epoch ($k=0$) and non-midnight epochs ($k\neq 0$). 
For non-midnight epochs, given an action $f = {f_{i,j}}$, the transition from $s$ to $s'$ follows: 
 \begin{equation}
     \begin{aligned}
     x_j'&=x_j^+ +a_j-d_j
     =\left(x_j+\sum_{i=1,i\ne j}^J f_{i,j}-\sum_{k=1,k\ne j}^J f_{j,k}\right) +a_j-d_j;\quad j=1,\dots,J, 
     \\
     y_j'&=y_j-d_j;\quad j=1,\dots,J, 
     \\
     h'&=(h +1)\mod m, 
 \end{aligned}
 \label{Equ: transition-current to next}
 \end{equation}
where $a_j$ and $d_j$ denote, respectively, the number of class-$j$ arrivals and departures from pool $j$ between the two epochs. The distributions of arrivals and departures are specified in Appendix~\ref{app: transition dynamics and probability}. For the midnight epoch, the main difference in the transition is in the  to-depart count, with
$y'_j \sim \text{Bin}(\min\{x_j,N_j\},\mu_j)$ for $ j=1,\dots,J$, meaning that each customer currently in service independently departs on the next day with probability $\mu_j$. 
We denote the transition probability from state $s$ to state $s'$, through $s^+$ under action $f$, is denoted as $p(s'|s,f)$. The detailed calculation of these transition probabilities is provided in Appendix~\ref{app: transition dynamics and probability}.  

\paragraph{Cost.}
As discussed, there are negative consequences associated with both excessive waiting (system congestion) and overflow assignments. To capture this trade-off, we consider two types of costs: 
(i) \emph{holding cost}: each waiting customer of class $j$ incurs a cost $C_j$ at each epoch; (ii) \emph{overflow cost}: each overflow assignment from class $i$ to pool $j$ incurs a one-time cost $B_{i,j}$. Given pre-action state $s$ and action $f$, the one-epoch cost is the sum of these two components, expressed as: 
 \begin{equation*}
g(s, f) = \sum_{j=1}^J C_j q_j^+ +\sum_{i=1}^J\sum_{j=1,j\ne i}^J B_{i,j}f_{i,j},
 \label{Equ: cost}
 \end{equation*}
where $q_j^+ =(x_j^+ - N_j)\vee 0$ denotes the post-action queue length for class $j$.

\paragraph{Policy and objective.}
We consider the long-run average cost objective since overflow decisions have long-term impacts on the system. We focus on all feasible stationary Markovian policies, and we denote the policy space as $\Pi$. For a policy $\pi \in \Pi$, we use $\pi(f\mid s)$ to denote the probability of selecting a feasible action $f$ in a given state $s$ (where deterministic policies are included via using the Dirac measure). The resulting stochastic process $\{S(t_k)\}$ forms a discrete-time Markov Chain (DTMC) with a one-epoch transition probability from state $s$ to $s'$ as $p_{\pi}(s'|s)=\mathbb{E}_{f\sim \pi(\cdot|s)} p(s'|s,f)$. The expected one-epoch cost is  $g_{\pi}(s)=\mathbb{E}_{f\sim \pi(\cdot|s)}g(s,f)$. Correspondingly, the objective is to identify the optimal policy 
\begin{equation}
\pi^{*}=\arg\min_{\pi\in\Pi}\lim_{T\rightarrow\infty}\frac{1}{T}\mathbb{E}_{\pi}\left[\sum_{t=1}^T\sum_{k=0}^{m-1}g_{\pi}\big( S(t_k) \big) \right].
\label{eq:obj}
\end{equation}
Under the periodic arrival assumptions, the system reaches a periodic steady state~\citep{liu2011large}. 
Without loss of generality, we assume that there exists a stationary Markovian policy that is optimal; see \cite{puterman2014markov} for conditions under which this holds.
We conclude this section with the following remarks on possible extensions. 

\noindent\textbf{Extensions.}  
We can use the concise state representation with customer count $X_j$ because we assume (i) primary assignments are handled automatically on-the-fly, allowing us to calculate both the waiting count and the in-service count from $X_j$, and (ii) the departures depend only on the server pool, so we can disregard customer class information once a customer is admitted into service. For (i), our analytical framework can be extended to a bed assignment problem that incorporates primary assignments (not assumed to be handled automatically). This would require augmenting the state space by separating $X_j$ into the waiting count $Q_j$ and the in-service count $Z_j$, where $Q_j$ tracks all customers who arrived between decision epochs, including those who could have been assigned to primary pools upon arrival but are now waiting. Assumption (ii) is made due to the relatively consistent average length of stay across medical specialties (approximately 4 days), while it could be relaxed in a similar way by augmenting the state space.

While our baseline setting does not consider transferring overflow customers back to their primary pool, this can be accommodated by adding actions that represent the number of in-service transfers from pool $j$ back to their primary pools. In addition, although we assume a linear cost function, our framework can handle more complex cost structures, such as polynomial costs or costs related to service levels (e.g., the proportion of patients waiting more than 4 hours). See the appendix of \cite{dai2019inpatient} for an example to account for service-level metrics. We note that these adjustments do not alter the core elements of our problem, and our algorithmic framework is both flexible and general enough to accommodate them. For ease of exposition, we will focus on the current baseline setting, where primary assignments are handled automatically, departures are pool-dependent, no transfer back of overflow customers, and the cost is linear; extensions are left as future research.


\section{Atomic Actions and Proximal Policy Optimization}

One of the core challenges in solving the MDP model developed in Section~\ref{Sec:MDP model} is the combinatorial explosion of the action space. 
Specifically, when there are enough idle servers in the system, the action space has a size of $\prod_{j=1}^{J} \left(\begin{array}{c}
     q_j+w_j  \\
     w_j 
\end{array}
\right)$, where $q_j$ is the number of waiting customers in class $j$ and $w_j$ is the number of overflow pools for class $j$ that have idle servers. To illustrate, in a twenty-pool system with each server-pool having a capacity of 60 and a state $(x_1, x_2, x_3, \dots, x_{20}) = (65, 63, 62, 50, 50, 50, \dots, 50)$, the number of potential ways to assign 10 waiting customers can exceed $10^9$. Even if the value function is computed (either exactly or via approximation), brute-force searching through all possible actions to find the one with the lowest cost-to-go is infeasible. To overcome this curse of dimensionality, we use atomic actions in combination with PPO that form the key foundation for our algorithm's scalability. In Section~\ref{Sec: Atomic actions}, we introduce the atomic action and policy setup. In Section~\ref{sec:ppo-method}, we present the PPO framework.

\subsection{Action Space Reduction: Atomic Actions}
\label{Sec: Atomic actions}

Atomic actions break the original action into a sequence of individual ones. That is, instead of solving the assignments for all waiting customers at once, the atomic action approach decomposes the problem into smaller, individual actions, sequentially making overflow/wait decisions for each customer. This sequential decision-making process is referred to as the atomic decision process.

To elaborate, for a given state $s = (x_1, \dots, x_J, y_1, \dots, y_J, h)$ observed at epoch $t_k$, the total number of waiting customers across all classes is $q = \sum_{j=1}^J q_j$, where $q_j = (x_j - N_j) \vee 0$. Following a pre-determined rule (e.g.,
first-come-first-serve), we determine an order for these $q$ waiting customers and conduct atomic actions $\mathbf{a} = \{a^1, a^2, \dots, a^q\}$ for each of them sequentially. Here, $a^n$ represents the assignment for the $n$th customer in the chosen order. The action space for each $a^n$ is simply $\mathcal{J} = \{1, \dots, J\}$.  
If the class of the $n$th customer is $j$, then $a^n = j$ indicates the customer will continue waiting until the next decision epoch; otherwise, the customer is assigned to an overflow pool (see details on feasible overflow assignment in Section~\ref{subsec:policy-para}).  

\begin{figure}[htp]
    \centering
    \includegraphics[]{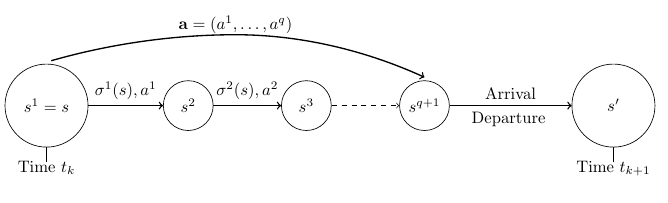}
    \caption{Illustration of the atomic decision process.
    }
    \label{fig:sub MDP}
\end{figure}

Figure~\ref{fig:sub MDP} illustrates the atomic decision process for deciding and executing the atomic actions. We denote the sequence of customer classes of the chosen order as $\boldsymbol{\sigma}(s)=(\sigma^1(s),\sigma^2(s),\cdots, \sigma^q(s))$, with $\sigma^n(s)$ representing the class of the $n$th waiting customer within the chosen order. Starting from the first customer, who is from class $\sigma^1(s)$, we execute the atomic action $a^1$ for this customer, then we update the state according to $a^1$ and sequentially move to the second customer, so on and forth till we execute the atomic action $a^q$ for the last waiting customer, who is from class $\sigma^q(s)$. 
The atomic action $a^n$ for the $n$th customer is chosen based on the \textit{atomic state}, $s^n$, which is observed after assigning the $(n-1)$th customer. This state has the same structure as the one introduced in Section~\ref{Sec:MDP model},  
i.e., $s^n = (x^n_1, \dots, x^n_J, y_1, \dots, y_J, h)$.  
After executing $a^n$, the atomic decision process updates the state to $s^{n+1}$ as 
\begin{equation}
    x^{n+1}_j=x^n_j+\sum_{k\ne j}\mathbbm{1}(a^n=j,\sigma^n(s)=k)-\sum_{k\ne j}\mathbbm{1}(a^n=k, \sigma^n(s)=j), \forall j\in \mathcal{J}, 
    \label{Equ: atomic dynamic}
\end{equation}  
where $\mathbbm{1}(\cdot)$ is the indicator function. This updated atomic state will be used for deciding $a^{n+1}$.

We note that even though the action is still state-dependent, the largest action space size is $|\mathcal{J}| = J$, which is drastically smaller compared to the combinatorial action space in the original MDP problem. This reduction forms the most critical foundation for developing an algorithm that scales effectively with the number of customer classes and server pools.

For clarity, and to distinguish between the MDP state and action defined in Section~\ref{Sec:MDP model} and the ones used for the atomic decision process here, we will refer to those in Section~\ref{Sec:MDP model} as the \emph{system-level} state and action throughout the rest of the paper. Meanwhile, we will refer to the ones in the atomic decision process as the \textit{atomic} state and action. It is important to note that the system-level and atomic states and actions are closely related. Specifically, the initial atomic state $s^1$ corresponds to the system-level (pre-action) state $s$, and after all atomic actions are executed, the final atomic state $s^{q+1}$ corresponds to the system-level post-action state $s^+$. For any given pair $(\boldsymbol{\sigma}(s), \textbf{a})$, we can recover the system-level action $f_{i,j}$ via $\sum_{n=1}^q \mathbbm{1}(a^n=j, \sigma^n(s)=i)$.

Theoretical support for the proposed atomic action approach comes from a contemporary paper \cite{dai2025optimal}. Their Theorems 1 and 2 show that, in certain stochastic processing networks, (i) decomposing joint actions into atomic steps incurs no loss in optimality, and (ii) even enforcing step-independent policies (each atomic action follows the same rule) still achieves the optimal long-run average reward. This step-independent structure mirrors our setup in Section 4.1, where each atomic action is sampled from the same policy $\pi(\cdot \mid s^\ell)$ regardless of step index $\ell$.


\subsubsection{Randomized Atomic Action}
\label{Subsec: randomization}

In this paper, we use randomized policies to prescribe the atomic actions. We will first outline the specifics of this randomization and then discuss its advantages over conventional deterministic policies at the end of Section~\ref{sec:ppo-method}.

For a given atomic state $s^n$, we define a probability distribution over all possible atomic actions for each customer class. We denote this distribution for a class $i$ customer, given $s^n$, as $\kappa(\cdot|s^n,i)$, and the set of distributions as $\kappa=\{\kappa(\cdot|s^n,i), i\in\mathcal{J}\}$. 
Consequently, starting from the initial state $s^1 = s$ (the system-level pre-action state) and following a given order $\boldsymbol{\sigma}(s)=\{\sigma^1(s),\dots,\sigma^q(s)\}$ for $q$ waiting customers, the probability of selecting a sequence of atomic actions $\textbf{a}=\{a^1,\dots,a^{q}\}$ is given by $\prod_{n=1}^q\kappa(a^n|s^n,\sigma^n(s))$, 
where $(s^1,\dots,s^{q})$ denotes the corresponding atomic state sequence generated within the atomic decision process, as specified in Equation~\eqref{Equ: atomic dynamic}.

According to the correspondence between atomic and system-level actions, the probability of taking a system-level action $f$ under the given distributions $\kappa$ can be expressed as: 
\begin{equation}
\begin{aligned}
\pi(f\mid s)
    =&\sum_{\textbf{a}} \prod_{n=1}^q\kappa(a^n|s^n,\sigma^n(s))\cdot\mathbbm{1}\Big( \sum_{n=1}^q \mathbbm{1}(a^n=j, \sigma^n(s)=i) = f_{i,j}, \forall i,j \Big) .
\end{aligned}
\label{equ: system-level probability of atomic action policy}
\end{equation}
To simplify notation and align with the broader PPO algorithm literature, we use the system-level policy $\pi=\{\pi(f\mid s)\}$ in the rest of the paper when we refer to a policy. However, it is important to clarify that this system-level policy is driven by the underlying atomic-level probability distributions. In other words, policy $\pi$ is the outcome of the sequence of atomic actions governed by the distributions $\kappa$. Our goal is to optimize $\kappa$ in such a way that the induced system-level policy $\pi$ achieves the performance objective in~\eqref{eq:obj}.
We discuss how the randomized policy can be practically applied in decision-making in Section~\ref{subsec: managerial insights}.

\subsubsection{Policy Parameterization}
\label{subsec:policy-para}

As discussed, our ultimate goal is to optimize the atomic action probabilities $\kappa$. To achieve this, we parameterize these probabilities using a set of parameters $\theta \in \Theta$. For a given $\theta$, we denote the resulting probabilities for atomic actions across customer classes as $\kappa_\theta = \{\kappa_\theta(\cdot | s, i), ~ i \in \mathcal{J}\}$. The induced system-level policy, denoted by $\pi_\theta$, is then derived from these parameterized atomic probabilities via~\eqref{equ: system-level probability of atomic action policy}. We refer to $\theta$ as the \textit{policy parameter}. 

While there are various ways to parameterize $\kappa_\theta$, in this paper, we use Neural Networks (NN), where $\theta$ represents the NN parameters. Before selecting the atomic action for the $n$th customer, given the atomic state $s^n$ and customer class $\sigma^n(s) = i$, the NN outputs the probability of choosing $a^n = j$ as: 
\begin{equation}
    \kappa_{\theta}\big( a^n=j|s^n,\sigma^n(s)=i \big) = \frac{\text{exp}( g_{i,j}(s^n|\theta) )}{\sum_{k=1}^J \text{exp}(g_{i,k}(s^n|\theta) )}, 
    \label{Equ: randomized atomic action}
\end{equation}
where $g_{i,j} (s|\theta)$ is one of the $J\times J$ output neurons of a dense NN with input $s$ and parameters $\theta$.

Equation~\eqref{Equ: randomized atomic action} could assign a positive probability to non-feasible atomic actions, such as when the selected server pool is full or cannot serve class $j$. To ensure the feasibility of the prescribed atomic action, we set $\kappa_\theta(a^n | s^n, \sigma^n(s)) = 0$ for $a^n \not\in \mathcal{A}(s^n, \sigma^n(s))$ (i.e., force the probability of non-feasible actions to 0), and then rescale the remaining probabilities to sum to 1. This adjustment is equivalent to resampling: if a sampled atomic action based on $\kappa_\theta$ is not feasible, we continue resampling until a feasible action is selected.

Now, the task of determining an effective overflow policy is reduced to learning the policy parameter $\theta$. We describe how the PPO method can accomplish this goal in the following section. Note that if the NN has $D$ hidden layers and $O(J^2)$ neurons per layer, the size of the parameter $\theta$ is $O(DJ^4)$. This is a polynomial growth rate rather than an exponential one, which ensures that the PPO method remains scalable when the number of customer classes and server pools increases.


\subsection{Finding Optimal Policy: Proximal Policy Optimization}
\label{sec:ppo-method}

We employ the state-of-the-art policy gradient method, PPO, to learn the parameter $\theta$ with custom modifications tailored to our specific queueing network context. We provide an overview of the PPO framework in this section; we defer the details of our tailored modifications to Section~\ref{sec:ppo-tailor}.

The relative value function under the system-level policy $\pi$ satisfies the Poisson equation  
\begin{equation}
v_\pi(s) = g_{\pi}(s) - \gamma_\pi+\mathbb{E}_{s'\sim p_{\pi}(s'|s)}[v_\pi(s')], 
    \label{equ: def of relative value function}
\end{equation} 
where $\gamma_\pi$ denotes the long-run average cost under policy $\pi$, $g_{\pi}(s)=\mathbb{E}_{f\sim \pi(\cdot|s)}g(s,f)$ is the expected one-epoch cost, and $p_{\pi}(s'|s)=\mathbb{E}_{f\sim \pi(\cdot|s)} p(s'|s,f)$ is the transition probability. 
Correspondingly, the advantage function of policy $\pi$ is defined as 
\begin{equation}
A_\pi(s,f)=g(s,f)-\gamma_\pi+\mathbb{E}_{s'\sim p (\cdot|s,f)}[v_\pi(s')]-v_\pi(s), 
    \label{equ: def adv}
\end{equation}
i.e., it measures the relative benefit (cost reduction) of taking a specific action $f$ compared to the expected performance of following the current policy $\pi$ at a given state $s$. In the cost minimization setting, a negative advantage function value means that the action leads to a lower expected cost (which is desirable). Intuitively, policy gradient methods, including PPO, update the policy parameters in the direction that increases the probability of actions with negative advantage values and decreases the probability of those with positive advantage values.

Specifically, assume that the current policy $\pi_\eta \in \Pi$ is parameterized by $\eta \in \Theta$, using the atomic action parameterization described in Section~\ref{subsec:policy-para}. For notational simplicity, let $v_\eta$ and $A_\eta$ denote the relative value function and advantage function under $\pi_\eta$, respectively. To improve upon $\pi_\eta$ to get a better policy $\pi_\theta$, \cite{schulman2017proximal} suggests \emph{minimizing} the following PPO objective:
\begin{equation}
L(\theta) := \mathop{\mathbb{E}}\limits_{\substack{ s\sim\mu_\eta\\ f\sim\pi_\eta(\cdot|s)}} \max \left\{ r_{\theta,\eta}(f\mid s) A_\eta(s, f),
\quad \text{clip}\left( r_{\theta,\eta}(f\mid s) , 1 -\epsilon , 1 + \epsilon \right) A_\eta(s,f) \right\}.
\label{Equ: PPO Obj}
\end{equation}
Here, the distribution $\mu_\eta$ denotes the stationary distribution under policy $\pi_\eta$; $r_{\theta,\eta}$ denotes the probability ratio between $\pi_\theta$ and $\pi_\eta$:  
\begin{align}
    r_{\theta,\eta}(f\mid s):=\frac{\pi_\theta(f\mid s)}{\pi_\eta(f\mid s)},
    \label{Equ: policy distance}
\end{align} 
where $\pi_\theta(f\mid s)$ and $\pi_\eta(f\mid s)$ are defined as in~\eqref{equ: system-level probability of atomic action policy}, with $\kappa$ replaced by $\kappa_\theta$ and $\kappa_\eta$, respectively;  
and the clipped function  
\begin{equation*}
    \text{clip}(x , 1-\epsilon , 1+\epsilon):=\left\{
    \begin{array}{cc}
         1-\epsilon,& \text{if } x<1-\epsilon \\
         1+\epsilon,&  \text{if } x>1+\epsilon \\
         x,&\text{otherwise}
    \end{array}
    \right. . 
\end{equation*}
The difference $\epsilon$ is a hyper-parameter that can be tuned for better algorithm performance.

For the theoretical guarantee, we establish \textbf{Proposition~1}, which extends the performance improvement guarantee of PPO from the time-stationary, long-run average cost setting in \cite{dai2022queueing} to our more complex periodic environment. 
See the complete statement of Proposition~1 (alongside needed technical conditions) in Appendix~\ref{App: PPO Motivation}.
The key challenge in the proof lies in the periodicity of state transitions, which makes directly applying previous results inapplicable. To address this, we introduce a new framework using daily Markov chains (MCs), where one period (e.g., one day) is treated as a single transition step. This transformation allows us to redefine the objective function in a way that ensures aperiodicity to use previous theoretical frameworks.  
Proposition~1 establishes that the change in long-run average cost can be decomposed into two terms: one that dominates the improvement and another that diminishes at a faster rate. Since PPO minimizes the dominant term while keeping the policy update conservative to control the secondary term, it leads to a provable policy improvement guarantee.

In implementation, we estimate the objective function~\eqref{Equ: PPO Obj} using simulation data samples, i.e., 
\begin{equation}
   \begin{aligned}
        \hat L(\theta, \mathcal{D}_\eta^T, \hat{\textbf{A}}_\eta(\mathcal{D}_\eta^T)) = 
        \frac{1}{Tm} \sum_{t=1}^T\sum_{k=0}^{m-1} \max &\left\{ 
      \hat{r}_{\theta,\eta}(f(t_k)|s(t_k)) \hat{A}_\eta(s(t_k),f(t_k)),\right.\\
    &\left.\text{clip}\left(\hat{r}_{\theta,\eta}(f(t_k)|s(t_k)) , 1-\epsilon , 1+\epsilon\right) \hat{A}_\eta(s(t_k), f(t_k))
    \right\}, 
   \end{aligned}
   \label{equ: estimated PPO obj}
\end{equation}
where the dataset $\mathcal{D}_\eta^T$ 
contains $T$ days of simulation data generated under policy $\pi_\eta$ (including the atomic states, atomic actions and order of patients within the simulation of sequential making process for each day), and the vector $\hat{\textbf{A}}_\eta(\mathcal{D}_\eta^T) = \left\{\left( \hat A_\eta(s(t_k),f(t_k))\right)_{k=0}^{m-1}\right\}_{t=1}^T$
consists of the estimated advantage functions for each decision epoch within the $T$ days using the collected data $\mathcal{D}_\eta^T$. 
The detailed steps of the PPO algorithm are given in Algorithm~\ref{Alg1: PPO}.

\begin{algorithm}[htp] 
\caption{Proximal policy optimization}
\label{Alg1: PPO}
\textbf{Output:} Policy $\pi_{\theta^*}$\\
\textbf{Input:} Policy $\pi_\eta$
\begin{algorithmic}[1]
\FOR{policy iteration $r=0,1,2,...,R-1$ }

    \STATE Run policy $\pi_\eta$ for $T$ days and collect simulation dataset $\mathcal{D}_\eta^T$. 
    
    \STATE Estimate the advantage functions $\hat{\textbf{A}}_\eta(\mathcal{D}_\eta^T)$. 
    
    \STATE Minimize the estimated objective function \eqref{equ: estimated PPO obj} w.r.t $\theta\in \Theta$ and obtain parameters $\theta^*$.
    
    \STATE Update $\eta\leftarrow \theta^*$.
\ENDFOR
\end{algorithmic}
\end{algorithm}

Computing $\pi(f \mid s)$ in \eqref{equ: system-level probability of atomic action policy} is computational expensive due to the required summation. In implementation, we approximate it using the product $\prod_{n=1}^{q(t_k)}\kappa \big( a^n(t_k)\mid s^n(t_k),\sigma^n(t_k) \big)$, evaluated along the realized processing order in simulation; correspondingly, we compute an approximate policy ratio $\hat r_{\theta,\eta}$. In Section~\ref{subsec:alg-implement}, we further introduce a batching scheme that eliminates order dependence and yields a simpler, more direct mapping between the system-level policy $\pi$ and the atomic-level distributions $\kappa$.

We conclude this section with two important remarks. First, our atomic policy parameterization differs from the conventional ones used in PPO framework, which directly parameterize a randomized policy over the entire system-level action space and would otherwise require an intractably large policy network. Second, performance improvement in PPO relies on conservative policy updates (Proposition~1 in Appendix~\ref{App: PPO Motivation}), which is infeasible under the deterministic policy setup. In particular, while atomic actions alone reduce the action space, combining them with traditional value-based ADP leads to severe policy ``chattering.'' This phenomenon arises because, when value functions are approximated rather than evaluated exactly, greedy policy updates may overreact to estimation errors, causing oscillations between improving and worsening policies~\citep{bertsekas2011dynamic}. PPO’s proximal updates mitigate this issue. In other words, the \emph{combination} of atomic actions and PPO is critical; see empirical results in Section \ref{sec:case-study}.


\section{Tailoring Design for Periodicity in MDP}
\label{sec:ppo-tailor}

A salient feature of our problem setting is the intrinsic time-dependency. The optimal action at different decision epochs within a day can vary significantly, even when the customer counts are similar. While PPO is well-established for stationary environments, applying it to time-nonstationary systems presents additional challenges. Specifically, our problem considers the long-run average cost in a periodic, time-varying setting. This makes ad hoc adjustments such as rolling-horizon methods ineffective. We address the challenges by integrating time-periodicity into the PPO framework and leveraging the queueing system structure. 
In Section~\ref{Subsec: Policy representation}, we introduce a tailored NN design to represent policies in a periodic setting, capturing both differences and similarities between policies at different epochs. This design is critical for improving sample efficiency (i.e., requiring less simulation data and enabling more efficient computation).
In Sections~\ref{Subsec: policy evaluation} and~\ref{subsec:alg-implement}, we introduce two more elements to further improve algorithm performance and reduce computational time: (i) value function approximations that leverage the queueing dynamics under the randomized atomic policy, and (ii) a batching setup for atomic actions.

\subsection{Policy Network Design}
\label{Subsec: Policy representation} 

\begin{figure}[htbp]
\centering 
\includegraphics[scale=0.6]{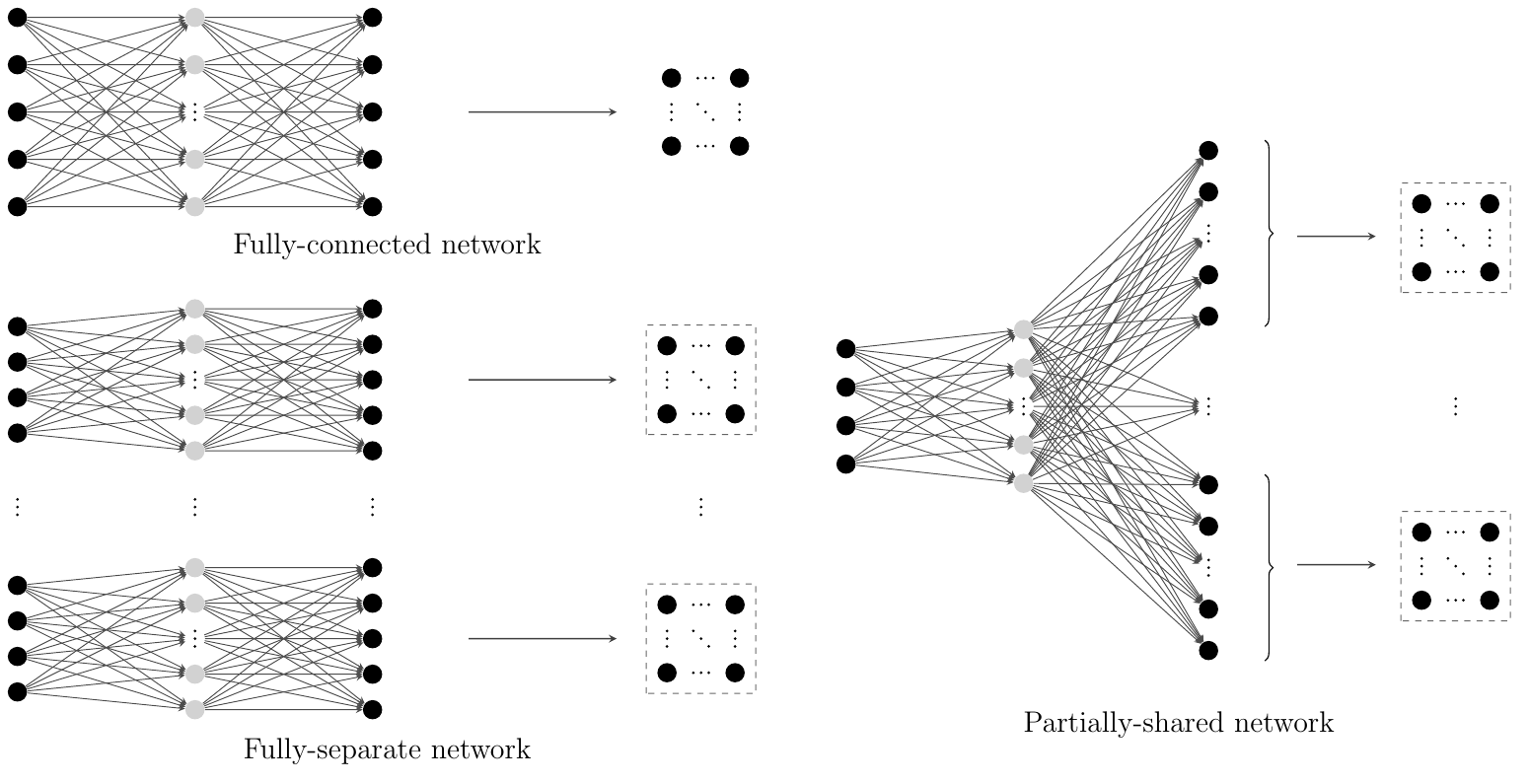}
\vspace{-0.2in}
\caption{Three policy network structures. The fully-connected network takes epoch index as part of the input and uses a single network to learn actions across all epochs. The fully-separate network trains $m$ separate networks for each of the $m$ epochs. The partially-shared network excludes the epoch index from input but uses different blocks of output neurons for actions of different epochs. 
}
\label{fig:Policy network structures}
\end{figure}

As discussed in Section~\ref{Subsec: randomization}, we use NNs to parameterize the atomic action probabilities $\kappa_\theta= \{\kappa_\theta(j|s,i), i,j\in\mathcal{J}\}$, which induce the system-level policy $\pi_\theta$. We refer to these NNs as the ``policy network'' and $\theta$ as the NN parameters. 
To accommodate the periodic setting, we consider three different designs for the policy network, as shown in Figure~\ref{fig:Policy network structures}. The first two are intuitive adaptations of policy networks from conventional time-stationary settings, while the third is our own innovation. We describe each design below. 

\vspace{-0.1in}
\begin{itemize}[leftmargin=*]
    \item The \textbf{fully-connected} network structure includes the epoch index as part of its input, and generates atomic action probabilities of different epochs with the same set of NN parameters. Specifically, given state $s \in \mathcal{S}$, we use one-hot encoding to include the epoch index $h$ in the input, encoding $s$ as $(\tilde{s}, e^h)$, where $\tilde s=(x_1,\dots,x_J,y_1,\dots,y_J)$ is the state excluding epoch index, and $e^h\in\mathbb{R}^m$ is the unit vector with the $h$th element being 1 and all others as 0. The action probabilities are then given by:
    \begin{equation*}
    \kappa_{\theta}(j|s=(\tilde s,h),i) = \frac{\text{exp}( g_{i,j}(\tilde s,e_h|\theta) )}{\sum_{k=1}^J \text{exp}(g_{i,k}(\tilde s,e_h|\theta))}, \quad i,j\in\mathcal{J}, s\in\mathcal{S},
\end{equation*}
where $g_{i,j}$ is one of the $J\times J$ output neurons of the fully-connected network. This design uses the same set of parameters to generate outputs for different epochs, which captures temporal differences by incorporating time as part of the input. 
    \item The \textbf{fully-separate} network structure uses $m$ independent networks for each of the $m$ different epochs, so the epoch index does not need to be part of the input. Specifically, the NN parameters $\theta = (\theta^0, \dots, \theta^{m-1})$, where $\theta^h$ is the parameters for the network corresponding to epoch $h$. The atomic action probabilities are given by: 
\begin{equation*}
    \kappa_{\theta}(j|s=(\tilde s, h),i) = \frac{\text{exp}( g_{i,j}^h(\tilde s|\theta^h) )}{\sum_{k=1}^J \text{exp}(g_{i,k}^h(\tilde s|\theta^h))}, \quad i,j\in\mathcal{J}, s\in\mathcal{S},
\end{equation*}
where $g_{i,j}^h$ is one of the $J \times J$ output neurons of the network for epoch $h$. This structure allows each epoch to have its own independent set of network parameters to capture the differences between them, but it does not exploit any potential similarities across epochs.
 
    \item The \textbf{partially-shared} network structure is a new structure we propose in this paper. This structure takes a middle ground between the fully-connected and fully-separate structures: it excludes the epoch index from the input while structuring the NN parameters such that the first several layers are shared across epochs and the later layers are specific to each epoch. Specifically, we define the NN parameters as $\theta = (\tilde\theta, \theta^1, \dots, \theta^m)$, where $\tilde\theta$ represents the shared parameters (typically for the first few layers), and $(\theta^1, \dots, \theta^m)$ correspond to the epoch-specific parameters (typically for the final layer). The atomic action probabilities are given by: 
     \begin{equation*}
    \kappa_{\theta}(j|s=(\tilde s, h),i) = \frac{\text{exp}( g_{i,j}^h(\tilde s|\tilde \theta,\theta^h))}{\sum_{k=1}^J \text{exp}(g_{i,k}^h(\tilde s|\tilde \theta,\theta^h))}, \quad i,j\in\mathcal{J}, s\in\mathcal{S}.
\end{equation*}
In other words, this structure does not differentiate states across epochs in the input but selects the corresponding block of output neurons for each epoch. This allows us to leverage data from all epochs to train the network while still prescribing different action probabilities for different times of the day.
\end{itemize}
To ensure feasibility, given state $s = (\tilde{s}, h)$, for a non-feasible pool $k \not\in \mathcal{A}(s,j)$ of class $j$, we can set $g_{j,k}(s|\theta) = 0$ for the fully-connected structure, $g_{j,k}^h(\tilde{s}|\theta^h) = 0$ for fully-separate, and $g_{j,k}^h(\tilde{s}|\tilde\theta, \theta^h) = 0$ for partially-shared, followed by re-normalizing the remaining probabilities. 

\smallskip  
\noindent\textbf{Combining advantages. }
The fully-connected and fully-separate designs each have their pros and cons, while the partially-shared policy network combines their strengths. The fully-connected network design uses the same set of parameters across all time epochs, naturally capturing similarities between epochs and enabling training with data from all time epochs. However, since the policy differences between epochs can be highly complex, a more sophisticated network (e.g., with additional hidden layers) may be needed to capture these nuances and hence more training data. The fully-separate network design accounts for differences between epochs with independent networks, but it would also require more time and data to train, as each network is only trained using data from its corresponding epoch.
The partially-shared network design strikes a balance: shared parameters in the early layers capture commonalities between epochs for efficient training, while the separated last layer adapts to epoch-specific policy outputs.  
In the case study of a ten-pool system, the partially-shared design requires only a few thousand simulation data samples for training, yet outperforms the other two designs that need \emph{100 times more} data; see more in Section~\ref{subsec:important-two-design}.


\subsection{Policy Evaluation: Linear Value Function Approximation}
\label{Subsec: policy evaluation}

In the standard PPO implementation, a separate neural network, known as the \textit{value network}, is used alongside the policy network to approximate value functions. The value network is trained with targets estimated via Monte Carlo evaluation, which can have high variance in long-run average settings. 
In this paper, we adopt a linear value function approximation and leverage the queueing system structure to design proper basis functions. To train the value approximation, we use the Least-Squares Temporal Difference (LSTD) method proposed by \cite{Bradtke1996}, which provides more stable estimates in tabular and linear function approximation settings, as shown in \cite{Boyan1999} and \cite{Lagoudakis2003}.

The key to accurate linear value function approximations lies in the choice of basis functions. Our proposed basis functions are: 
\begin{equation}
    V_d(s), ~x_j, ~x_j^2, ~y_j, ~y_j^2, ~x_jy_j, \quad j\in\mathcal{J}. 
    \label{equ: basis function choices}
\end{equation}
Here, the polynomial terms are the same as those in~\cite{dai2019inpatient}. We introduce a new design, $V_d(s) = \sum_{j=1}^J \hat v_{\pi,j}(s_j)$, which approximates the relative value function under the randomized atomic policy. Each $\hat{v}_{\pi,j}(s_j)$ relies only on pool-$j$’s state information, $s_j = (x_j, y_j, h)$, and estimates the relative value function based on the holding and overflow cost components for class-$j$ customers. 

Specifically, under the randomized atomic policy $\pi$, we decompose the original $J$-pool system into $J$ single-pool systems, and the relative value function $\hat v_{\pi,j}$ for the decomposed pool-$j$ system is solved via the post-action Poisson equation  
\begin{equation} 
\hat v_{\pi,j}(s_j) = C_j q_j -\hat\gamma_{\pi,j}+\mathbb{E} \left[\sum_{k=1}^J B_{j,k}\bar f_{j,k}+\hat v_{\pi,j}(s_j')\right]
    \label{Equ: decomposed Poisson equation-main}
\end{equation}
for a given post-action state $s_j$ (the superscript $+$ dropped for simplicity). Here, $C_j q_j$ is the post-action holding cost with $q_j=\max(s_j-N_j, 0)$, $\hat{\gamma}_{\pi,j}$ is the average cost of the decomposed pool-$j$ system, and the cost-to-go term consists of: (a) the expected overflow cost $\mathbb{E}[\sum_{k=1}^J B_{j,k}\bar f_{j,k}]$, where $\bar{f}_{j,k}$ approximates the overflow assignments from class $j$ to other pools, and (b) the expected value function at the next post-action state $\mathbb{E}\left[\hat{v}_{\pi,j}(s’_j)\right]$.
To evaluate (a) and (b), we need to properly adjust the queueing dynamics for each decomposed pool-$j$ system by incorporating overflow assignments from other classes while excluding overflow assignments from class $j$ to other pools (treated as ``diversions''). To achieve this, we approximate the probability of assigning a class-$j$ customer to pool $k$ under policy $\pi$ using $\bar\kappa_j(k|s_j,j)$, which depends only on pool-$j$’s state information. Given the arrivals are Poisson and leveraging the Poisson thinning property, the outflows $\bar{f}_{j,k}$ and inflows $\bar{f}_{i,j}$ can be approximated as Poisson distributions with rates $\lambda_j \bar{\kappa}_j(k|s_j,j)$ and $\lambda_i \bar{\kappa}_i(j|s_j,i)$, respectively. 
Details on computing $\bar{\kappa}_j$ and deriving $\hat{v}_{\pi,j}$ are provided in Appendix~\ref{App: pool-wise decomposition}.

\subsection{Batching Setup and Other Implementation Details}
\label{subsec:alg-implement}

In the general atomic decision process, each of the atomic actions sequentially depends on the atomic states as illustrated in Figure~\ref{fig:sub MDP}. 
We find that decision ordering has limited impact on overall performance; see Appendix~\ref{sub-app: order}. Motivated by this observation, we propose a batching setup to simplify the decision process, which uses the same initial state $s=s^1$ for all atomic actions in a batch. Specifically, given $s$ and a chosen customer class order $\boldsymbol{\sigma}(s)$, the probability of conducting atomic action sequence $\textbf{a}$ becomes $\prod_{n=1}^{q}\kappa(a^n|s,\sigma^n(s))$. 
In other words, during this batch, since the input state $s$ no longer changes, we prescribe atomic actions to the customers from the same class $i$ following the same distribution $\kappa(\cdot|s,i)$. 
Under the assumption that any system-level action $f$ satisfying $\sum_j f_{i,j} = q_i, \forall i$ is feasible, we have the following lemma. 
\begin{lemma}
Under the batching setup, for a given system-level state $s$ and any pre-determined rule for choosing customer order $\boldsymbol{\sigma}(s)$, the probability of taking action $f$  follows a multinomial distribution that is independent of the order $\boldsymbol{\sigma}(s)$, i.e.,  
\begin{equation}
   \begin{aligned}
\pi(f\mid s)=\prod_{i=1}^J\frac{q_i!}{\prod_{j=1}^J f_{i,j}!}\prod_{j=1}^J \kappa(j|s,i)^{f_{i,j}}.
   \end{aligned}
    \label{equ: system-level action probability}
\end{equation}  
\label{lemma:batching}
\end{lemma}
\begin{proof}{Proof. }
From~\eqref{equ: system-level probability of atomic action policy}, we have that 
\begin{equation*}
    \begin{aligned}       \pi(f\mid s)=&\sum_{\textbf{a}}\prod_{n=1}^{q}\kappa(a^n|s,\sigma^n(s))\cdot\mathbbm{1}\left(f_{i,j}=\sum_{n=1}^q \mathbbm{1}(a^n=j, \sigma^n(s)=i), i,j\in\mathcal{J}\right) \\       =&\left(\prod_{i=1}^J\prod_{j=1}^J \kappa(j|s,i)^{f_{i,j}}\right)\cdot \sum_{\textbf{a}}\mathbbm{1}\left(f_{i,j}=\sum_{n=1}^q \mathbbm{1}(a^n=j, \sigma^n(s)=i), i,j\in\mathcal{J}\right)\\
       =&\prod_{i=1}^J\frac{q_i!}{\prod_{j=1}^J f_{i,j}!}\prod_{j=1}^J \kappa(j|s,i)^{f_{i,j}} . 
    \end{aligned}
\end{equation*} 
Here, the second equality uses the fact that the atomic actions for class $i$ follow the same distribution $\kappa_\theta(\cdot|s,i)$, so that we can group the atomic actions for each of the $q$ waiting customers simply by their class $i$ and the assigned pool $j$, which connects to the $f_{i,j}$.
The last equality comes from the combinatorial partition of assigning the $q_i$ waiting customers from class $i$ to each pool $j$. \hfill $\square$
\end{proof}

The first advantage of this batching approach is that the atomic policy for the same class of customers remains \emph{consistent} within the same decision epoch, which ensures a degree of fairness regardless of the pre-selected processing order $\boldsymbol{\sigma}(s)$. The second advantage is the improved time efficiency in both the simulation and policy update processes. Unlike the original sequential atomic decision setup, we do not need to update each atomic state or compute the state-dependent atomic action probabilities for each individual customer. 
The third advantage lies in the great simplification of evaluating~\eqref{equ: estimated PPO obj}. Specifically, recall that we need to compute the probability ratio 
$$r_{\theta,\eta}(f\mid s)=\frac{\pi_\theta(f\mid s)}{\pi_\eta(f\mid s)}
$$
at each decision epoch, where $\pi_\theta(f\mid s)$ and $\pi_\eta(f\mid s)$ are defined as in~\eqref{equ: system-level probability of atomic action policy}, with $\kappa$ replaced by $\kappa_\theta$ and $\kappa_\eta$, respectively.
Directly calculating this ratio is computationally intensive because it involves summing over all possible atomic action sequences $\boldsymbol{a}$ corresponding to the system-level action $f$ for a given customer order $\boldsymbol{\sigma}(s)$ (i.e., the summation over $\boldsymbol{a}$ part in~\eqref{equ: system-level probability of atomic action policy}). In contrast, under the batching setup, using~\eqref{equ: system-level action probability} from Lemma~\ref{lemma:batching}, this ratio is greatly simplified to:
\begin{equation}
    r_{\theta,\eta}(f\mid s)=\prod_{i=1}^J\prod_{j=1}^J \left(\frac{\kappa_\theta(j|s,i)}{\kappa_\eta(j|s,i)}\right)^{f_{i,j}},
    \label{equ:ratio with batching}
\end{equation}
since the combinatorial factor $\frac{q_i!}{\prod_{j=1}^J f_{i,j}!}$ is a constant for given $f_{i,j}$'s and cancels out in the ratio calculation. This also establishes a direct \emph{equivalence} between the ratio of atomic action probabilities ($\kappa_\theta$, $\kappa_\eta$) and the ratio of the induced system-level policies ($\pi_\theta$, $\pi_\eta$).

As in the original setup, it is also possible that a sampled action under the batching setup could become infeasible if the selected pool no longer has an available server or cannot accommodate the customer class. To handle this, we apply the same resampling technique to ensure a feasible action is selected. 
Other implementation details are as follows. We use Tensorflow v1 to build the policy network training pipeline and use the Ray package for parallel simulation data generation. The PPO algorithm iterates through three major modules: simulation data generation, policy evaluation, and policy improvement. For data generation, multiple actors (CPUs) run parallel simulations of equal length, and then the simulated data is concatenated into one long trajectory, which constitutes the entire training dataset. For policy evaluation, we compute the linear value function approximations and the estimated advantage functions. For policy improvement, we update the policy network parameters $\theta$ using the Adam optimization method, running for $E$ training epochs (one training epoch corresponds to the number of complete passes through the entire training dataset). This process repeats until the average costs from two iterations differ by less than a small $\delta$. Table~\ref{tab:PPO hyper-parameters choice} summarizes important PPO hyper-parameters used in a ten-pool system in the case study in Section~\ref{sec:case-study}. More details on tuning these hyper-parameters are provided in Appendix~\ref{Subsec: additional NN structure results}.

\begin{table}[htbp]
    \centering
\scalebox{0.85}{
    \begin{tabular}{|c|c|}
    \hline
         Parameters & Baseline Choice \\ \hline
         Neural network structure & Partially-shared\\ 
         Neural network complexity &One hidden layer with $34$ neurons\\
         Basis function  &$(V_d,X,X^2,Y,Y^2,XY)$\\
         Initial policy  & Complete-overflow\\
         Simulation days per actor&10,000\\
         Number of actors  &10\\
         Number of training epochs $(E)$  &15\\
         Clipping parameter $(\epsilon)$  &0.5\\
         Gap tolerance $(\delta)$ &0.1\\
         \hline
    \end{tabular}
}
\caption{Baseline choice of PPO hyper-parameters in the ten-pool system.}
\label{tab:PPO hyper-parameters choice}
\end{table}


\section{Case Study}
\label{sec:case-study}

In this section, we evaluate the effectiveness of our tailored PPO algorithm through a case study on inpatient overflow assignment in large hospital systems, including settings with up to twenty patient classes and twenty wards (pools). In Section~\ref{sec:case-model-setting}, we introduce the case study settings. In Section~\ref{sec:five-pool-case}, we show that our PPO algorithm achieves similar performance to the ADP method from~\cite{dai2019inpatient}, considered the ``state-of-the-art,'' while significantly reducing computational time in the five-pool baseline, the largest system handled by~\cite{dai2019inpatient}. In Section~\ref{subsec:large-system-case}, we show the scalability in ten- and twenty-pool systems. In Section~\ref{subsec:important-two-design}, we examine the algorithm's performance across different hyper-parameter settings. We also conduct a detailed ablation study to highlight the critical role of our atomic policy design and the value of each tailored component. Finally, in Section~\ref{subsec: managerial insights} we showcase that the policy from the algorithm output is interpretable, and we summarize key takeaways for managing inpatient overflow.

\subsection{Model Settings}
\label{sec:case-model-setting}

The model settings for our case study are based on the five-pool system from~\cite{dai2019inpatient}, which is calibrated using aggregate patient-flow statistics from a large teaching hospital in Asia. The calibrated system has 315 servers (inpatient beds) and a daily arrival rate of 70 customers (patients). Time-varying arrival and discharge patterns are detailed in Appendix~\ref{app:case-study}. Each day is divided into eight decision epochs. Each class (medical specialty) has a dedicated primary ward, one preferred overflow ward, and two secondary overflow wards. The corresponding overflow costs are $B = (B_1, B_2)$, with $B_1 < B_2$ to reflect ward preferences. 
Following~\cite{dai2019inpatient}, we set $C = 6$ and test overflow costs ranging from 5 to 65, with $B = (30, 35)$ as the baseline.\footnote{Using estimates from prior studies linking waiting and overflow to inpatient LOS extensions (Singer et al. 2011, Howlett et al. 2026, Song et al. 2020), one overflow event is approximately comparable to 6.4--10.2 hours of waiting. Under our normalization of a holding cost of 6 per patient per 3-hour epoch (i.e., 2 per patient-hour), this corresponds to an overflow cost of 12.8--20.4 as an order-of-magnitude benchmark. Additional operational burdens, such as coordination effort and spillovers, may justify larger values. We therefore vary overflow costs over a broad range rather than relying on a single calibrated value.} 
As these costs are difficult to estimate accurately, we plot Pareto curves by varying them to show the trade-off between congestion and overflow rates, allowing managers to choose policies aligned with their target performance goals; see Appendix~\ref{sub-app:additional_metrics}.

\begin{figure}[htbp]
    \centering
\includegraphics{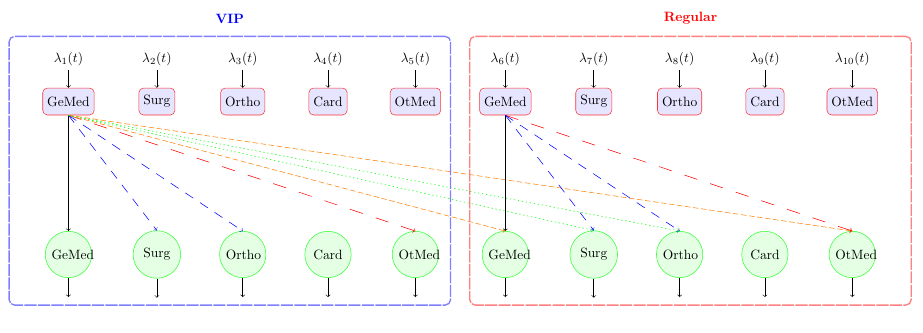}
    \caption{Illustration of the ten-class ten-pool system based on practice in our partner hospital.  The solid black arrow depicts the primary assignment route, while red and blue dashed arrows show preferred and secondary overflow assignments within a department. Additionally, orange and green dotted arrows represent cross-department overflow assignments (only from VIP to regular).
    }
    \label{fig: ten-pool routing}
\end{figure}

We extend this model to ten-pool and twenty-pool systems. The ten-pool system consists of two departments: a VIP department with smaller capacity and a regular department with larger capacity. Both departments follow the five-pool structure, and the total capacity is double that of the baseline five-pool system. Arrival rates are adjusted proportionally, as shown in Table~\ref{tab: parameter settings}. Overflow routes within each department remain unchanged (one preferred and two secondary overflow wards), but we also allow VIP patients to overflow to the regular department (not vice versa). Figure~\ref{fig: ten-pool routing} provides a schematic representation of the feasible routes in the ten-pool system, using the GeMed (General Medicine) specialty as an example. The diagram illustrates both within-department routes (black solid lines and red/blue dashed lines) and cross-department routes (orange and green dotted lines), highlighting the overflow pathways available within and between departments. Overflow costs are set at $B=(B_1, B_2, B_3, B_4)=(25, 30, 35, 40)$, with lower costs for within-department overflows $(B_1, B_2)$ and higher costs for cross-department overflows $(B_3, B_4)$. The unit holding cost is set at $C=6$ for regular patients and $C=7$ for VIP patients. For the twenty-pool system, we assume it consists of two hospitals, each modeled as a ten-pool system with different capacities, as shown in Table~\ref{tab: parameter settings}. Patients can overflow within their hospital or to the other hospital. Overflow costs are set at $B=(B_1,\dots,B_8)=(25,30,35,49,40,45,50,55)$, with lower costs for within-hospital overflows $(B_1,\dots,B_4)$ and higher costs for cross-hospital overflows $(B_5, \dots, B_8)$. The unit holding costs remain the same.

\vspace{-0.15in}
\begin{table}[htp]
\centering
\scalebox{0.75}{ 
\begin{tabular}{|ccc|c|c|c|}
\hline
\multicolumn{3}{|c|}{}                                                                                                & Number of inpatient beds & Daily arrival rate & Nominal utilization \\ \hline
\multicolumn{3}{|c|}{Five-pool}                                                                                       & 315                     & 70                 & 0.889      \\ \hline
\multicolumn{2}{|c|}{\multirow{3}{*}{Ten-pool}}                                                       & VIP dept.     & 210                     & 50                 & 0.962      \\ 
\multicolumn{2}{|c|}{}                                                                                & Regular dept. & 420                     & 90                 & 0.857      \\ 
\multicolumn{2}{|c|}{}                                                                                & Total         & 630                     & 140                & 0.889      \\ \hline
\multicolumn{1}{|c|}{\multirow{7}{*}{Twenty-pool}} & \multicolumn{1}{c|}{\multirow{3}{*}{Hospital 1}} & VIP dept.     & 175                     & 50                 & 1.143      \\ 
\multicolumn{1}{|c|}{}                             & \multicolumn{1}{c|}{}                            & Regular dept. & 385                     & 90                 & 0.835      \\ 
\multicolumn{1}{|c|}{}                             & \multicolumn{1}{c|}{}                            & Total         & 560                     & 140                & 1          \\ \cline{2-6} 
\multicolumn{1}{|c|}{}                             & \multicolumn{1}{c|}{\multirow{3}{*}{Hospital 2}} & VIP dept.     & 245                     & 50                 & 0.816      \\ 
\multicolumn{1}{|c|}{}                             & \multicolumn{1}{c|}{}                            & Regular dept. & 455                     & 90                 & 0.791      \\ 
\multicolumn{1}{|c|}{}                             & \multicolumn{1}{c|}{}                            & Total         & 700                     & 140                & 0.800      \\ \cline{2-6} 
\multicolumn{1}{|c|}{}                             & \multicolumn{1}{c|}{Total}                       & -             & 1260                    & 280                & 0.889      \\ \hline
\end{tabular}
}
\caption{Parameter settings for five-pool, ten-pool and twenty-pool models;
the latter two are formed by clustering and proportionally scaling the five-pool model in arrivals and capacities.}
\label{tab: parameter settings}
\end{table}

\subsection{Five-pool System: Comparing Performance with ADP}
\label{sec:five-pool-case}

Table~\ref{tab:five pool PPO} compares the performance of our PPO algorithm with the ADP method in the five-pool model. We consider two setups: a balanced setting where each pool has equal capacity, and an unbalanced setting where some wards have less capacity and higher congestion. 
These experiments show that the PPO method achieves performance comparable to ADP (within 4\% difference) -- the best-performing method for the five-pool model to date -- while significantly reducing computational time. Specifically, PPO completes one iteration in about two hours, compared to over ten hours for ADP, since its atomic policy design eliminates the need for exhaustive search over all feasible actions. Furthermore, in the simpler two-pool systems where the exact optimal policy is known, PPO achieves near-optimal results (within 3\% difference).

\begin{table}[htp]
\centering
\scalebox{0.85}{ 
    \begin{tabular}{|c|c|c|c|c|c|c|}
    \hline
                & \multicolumn{3}{c|}{Eight-epoch; Unbalanced} 
                & \multicolumn{3}{c|}{Eight-epoch; Balanced}\\ \cline{2-7}
                
                & B=(15,25) & B=(30,35) & B=(40,45) &B=(15,25) & B=(30,35) & B=(40,45)\\ \hline

        ADP &152.68$\pm$0.76&196.31$\pm$0.92& 250.39$\pm$0.45&132.46$\pm$0.32&193.03$\pm$0.62&248.48$\pm$0.71\\\hline
        
        PPO   
        &153.17$\pm$0.34 
        &202.37$\pm$0.65
        & 260.87$\pm$0.40
        
        &134.59$\pm$0.82 
        &196.78$\pm$1.06
        &258.32$\pm$0.98\\\hline
        
    \end{tabular}
}
\caption{\textbf{Five-pool system: long-run average costs.} We set $N_1=60,N_2=64,N_3=67,N_4=N_5=62$ in unbalanced setting and $N_1=\cdots=N_5=63$ in balanced setting. Original atomic decision (non-batched) is used here. 
} 
\label{tab:five pool PPO}
\end{table}

\subsection{Large Systems: Algorithm Scalability}
\label{subsec:large-system-case}
 
\begin{figure}[htp]
\subfigure[Ten-pool]{        
\includegraphics[width=0.37\linewidth]{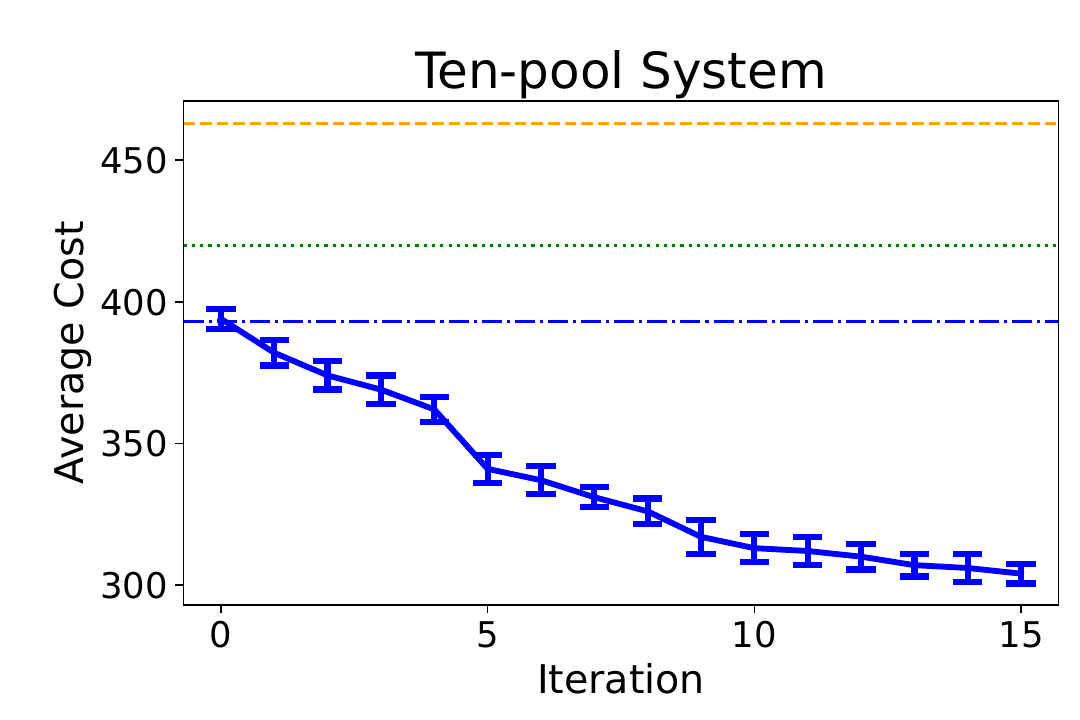}
\label{fig:10pool PPO}
}
\subfigure[Twenty-pool]{  
\includegraphics[width=0.38\linewidth]{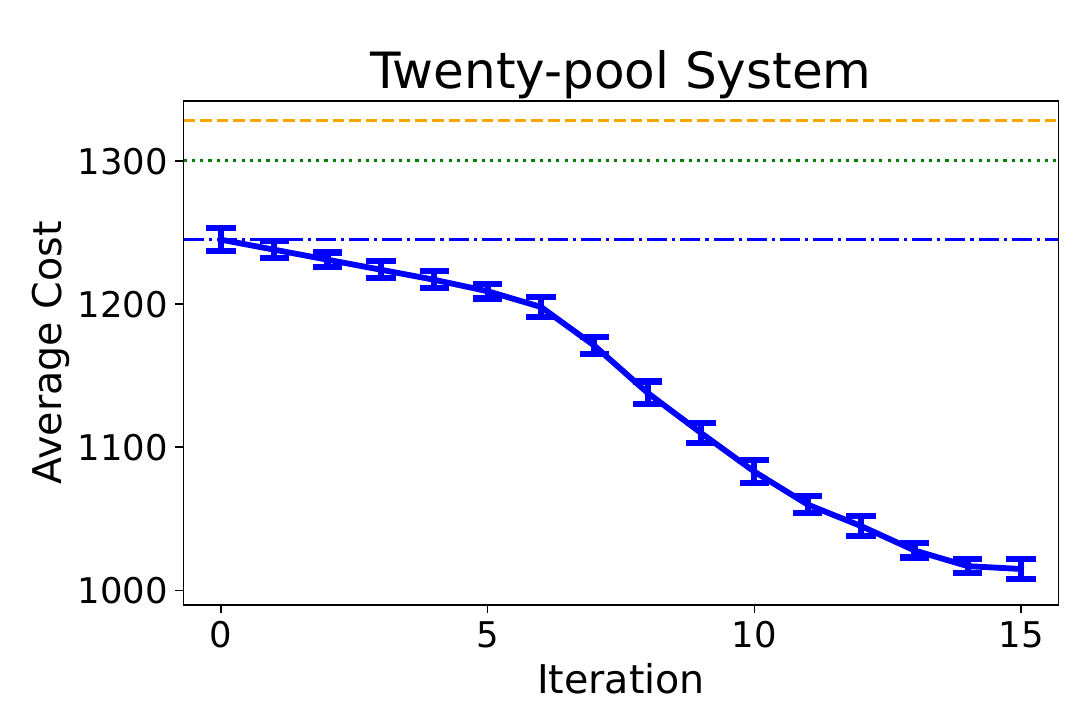}
\label{fig:20pool PPO}
}
\raisebox{1cm}{
\subfigure{
\hspace*{-0.5cm}
\includegraphics[width=0.2\linewidth]{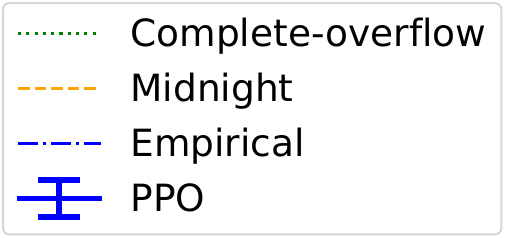} 
}}
\caption{\textbf{Performance comparison in the ten-pool and twenty-pool systems.} 
} 
\label{fig:scalability}
\end{figure}

In larger systems with more than five pools, the ADP method becomes infeasible due to the need to exhaustively search over the combinatorially large action space. Therefore, we compare the performance of the PPO algorithm against three benchmark policies introduced in \cite{dai2019inpatient}: (a) the complete-overflow policy, where overflow assignments are allowed at each decision epoch; (b) the midnight policy, allowing overflow assignments only at midnight each day; and (c) the empirical policy, which reflects current practice by allowing overflow assignments only during nighttime (7 pm to 7 am the next day). In all three policies, when overflow assignments are allowed, waiting patients are assigned one by one: first to their preferred overflow ward if feasible, and then to a randomly chosen secondary ward, until no waiting patients remain or no idle beds are available.
Figure~\ref{fig:scalability} shows the performance comparison in ten and twenty-pool systems, where PPO significantly outperforms all benchmark policies. In the ten-pool system, PPO achieves an average cost of 300, approximately 23\% lower than the best-performing benchmark (the empirical policy with an average cost of 390). For the twenty-pool system, the improvement is similar: PPO reduces the average cost to 1000, about 25\% lower than the empirical policy (which has an average cost of 1250).
The computational time is about 3.5 hours per iteration for the ten-pool system and about 8 hours per iteration for the twenty-pool systems. Note these are training times; once the policy networks are trained, the execution of the policy is almost instantaneous in real-time deployment. 
These results highlight the scalability and effectiveness of our PPO algorithm in large-scale systems. 

\vspace{-0.15in}

\subsection{Importance of Tailored Designs}
\label{subsec:important-two-design}

In this section, we evaluate the performance of different hyper-parameter combinations in the ten-pool system, using the values from Table~\ref{tab:PPO hyper-parameters choice} as the baseline. While hyper-parameters for NNs play a role in algorithm performance as traditionally understood, we find that they are \emph{less} critical than the design of the policy network, basis functions, and batching setup -- the tailored elements introduced in this paper. These results highlight the importance of adapting a general-purpose algorithm to the specific problem domain for optimal performance and computational efficiency.  Table~\ref{tab:ablation summary} summarizes our ablation study and the relative impact of each tailored component on performance and computation time.

\begin{table}[htp]
\centering
\scalebox{0.95}{
\begin{tabular}{|l|cc|}
\hline
\textbf{Configuration} 
& \textbf{Performance Change} 
& \textbf{ Time Change} \\ \hline
Full algorithm 
& --- & --- \\\hline
\textit{Remove PPO (atomic + greedy search)} 
& $\downarrow\downarrow$  & $\uparrow\uparrow$ \\
\textit{Remove tailored policy network (alternative NN designs)} 
& $\downarrow\downarrow$ & $\sim$ \\
\textit{Remove tailored basis (basis in \cite{dai2019inpatient})} 
& $\downarrow$ & $\sim$ \\
\textit{Remove LSTD (NN critic)} 
& $\downarrow$ & $\uparrow$ \\
\textit{Remove batching (sequential updates, various orders)} 
& $\sim$ & $\uparrow\uparrow$ \\
\textit{Fewer epochs} 
& $\downarrow$ & $\downarrow$ \\
\textit{Fewer simulation days} 
& $\sim$ & $\downarrow$ \\ \hline
\end{tabular}}
\caption{Ablation study summary. For each variant, the term in parentheses indicates the replacement mechanism. ``$\uparrow\uparrow$'' or ``$\downarrow\downarrow$'' indicate large increase or decrease in performance (lower average cost) and computation time ($>5\%$ cost or $>25\%$ time difference), compared with the full algorithm (integrating all elements: Atomic + PPO + Tailored Policy Network + Tailored LSTD + Batching); ``$\uparrow$'' or ``$\downarrow$'' moderate increase or decrease ($2\sim5\%$ cost or $10\%\sim 25\%$ time difference); and ``$\sim$'' negligible change ($<2\%$ cost or $<10\%$ time difference).}
\label{tab:ablation summary}
\end{table}

\noindent\textbf{Is PPO necessary? } 
Row 2 of Table~\ref{tab:ablation summary} examines the effect of removing PPO and replacing it with a conventional value-based method (ADP) paired with atomic actions. In this setup, policy updates are performed via greedy search over all possible atomic actions. While atomic actions reduce the action space and make updates tractable, combining them with greedy search introduces significant instability due to the well-known ``chattering'' phenomenon~\citep{bertsekas2011dynamic}. The average cost {varies} between 259.24 and 282.77, more than 27\% increase from our PPO algorithm (see Appendix~\ref{subsubsec: atomic_adp} for details). The chattering arises from the compounding effect of (i) the greedy policy updates and (ii) the imperfect approximation of the value function, even with carefully chosen basis functions. These challenges are especially pronounced in long-run average cost settings, where accurately estimating the average cost is difficult, which in turn makes reliable value function estimation more challenging (Bertsekas et al. 2011, Vol. 2, Chapter 6.6).
Moreover, the atomic setup, which evaluates actions at an individual (rather than system) level, further exacerbates the chattering since each action’s marginal impact is highly sensitive to estimation noise. In contrast, pairing atomic actions with PPO mitigates these issues via incremental policy updates (as formalized in Proposition~1) and  updating actions simultaneously via the system-level policy ratio~\eqref{Equ: policy distance}. 
Together, these results highlight that the combination of atomic actions with PPO is essential in our infinite-horizon, long-run average context.

\begin{figure}[htp]
\centering  
\subfigure[Impact of simulation days]{   

\includegraphics[width=0.35\linewidth]{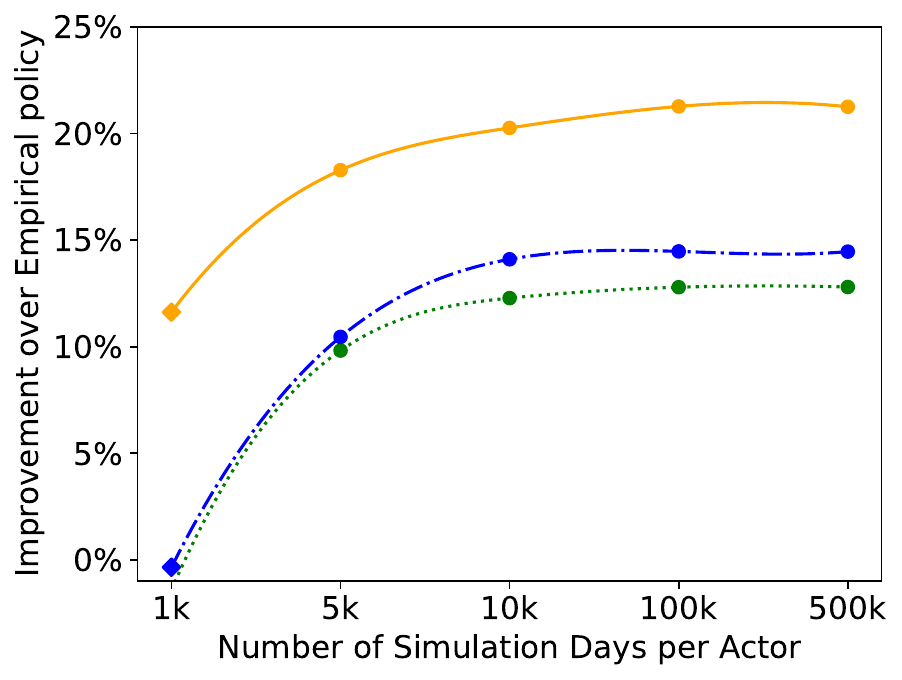} 
}
\subfigure[Impact of training epochs\protect\footnotemark]{ 
\includegraphics[width=0.35\linewidth]{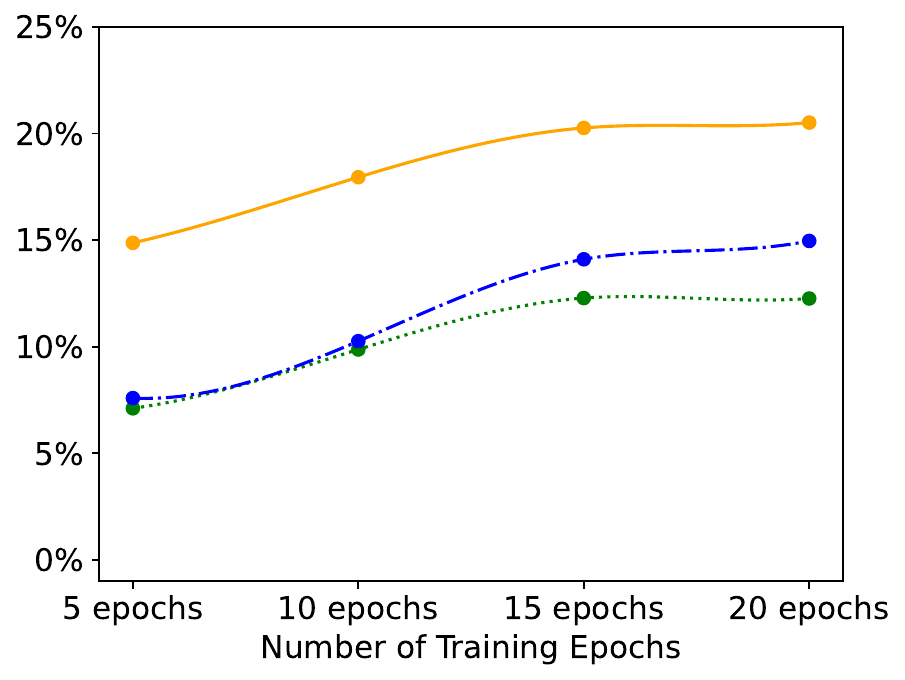}
}
\raisebox{1cm}{
\subfigure{
\includegraphics[width=0.2\linewidth]{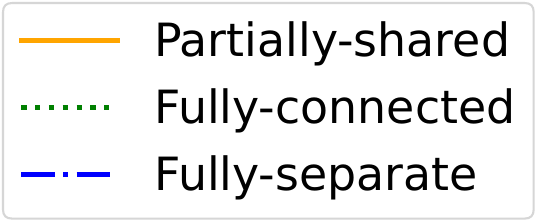} 
}}
\caption{Impact of training days and epochs. Performance is measured as the \% improvement in average cost between the PPO policy and the best-performing benchmark (the empirical policy).}
\label{fig:trade-off1}    
\end{figure}

\begin{figure}[htp]
\centering  

\subfigure[Impact of NN depth]{  

\includegraphics[width=0.35\linewidth]{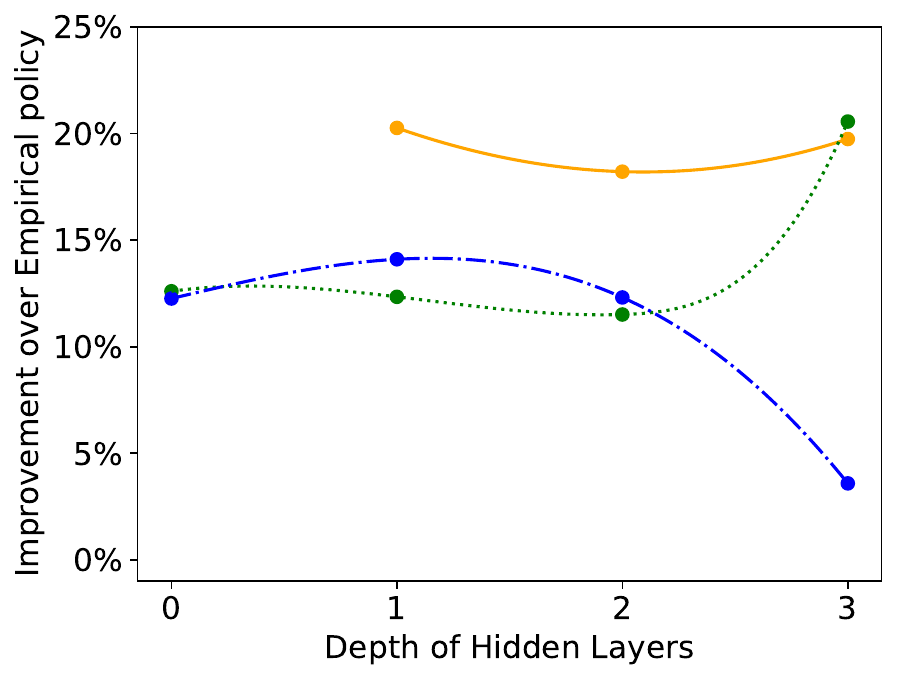}  
}
\subfigure[Impact of NN width]{ 
\includegraphics[width=0.35\linewidth]{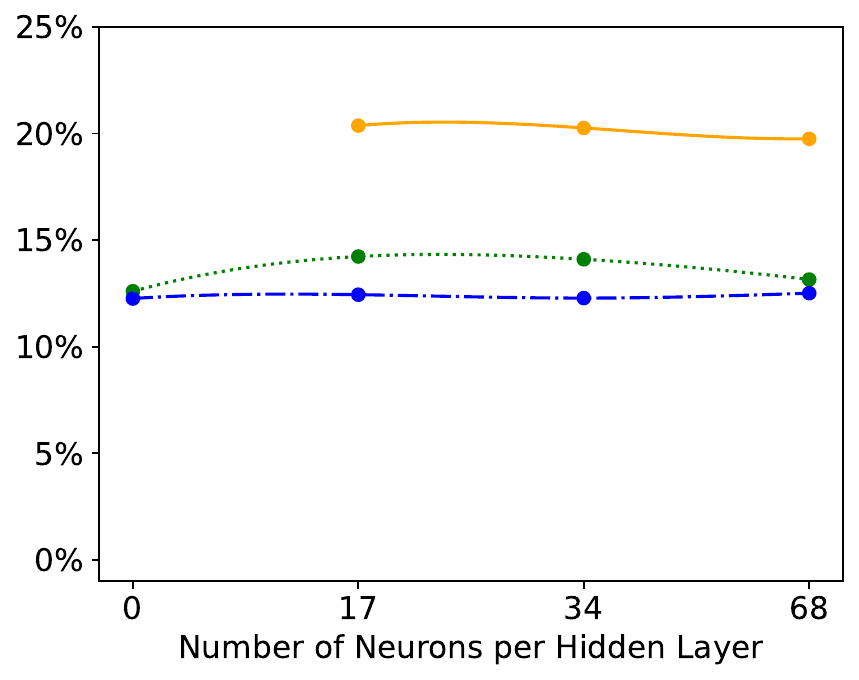}
}
\raisebox{1cm}{
\subfigure{
\centering
\includegraphics[width=0.2\linewidth]{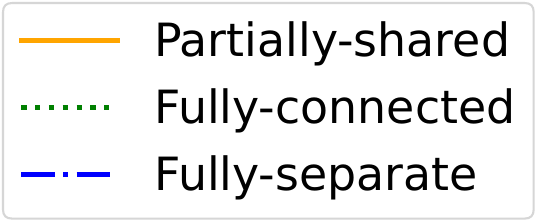} 
}}
\caption{Impact of NNs' depth and width. Performance is measured as the \% improvement in average cost between the PPO policy and the best-performing benchmark (the empirical policy). A value of ``0'' in the depth or width indicates no hidden layer, which is possible for fully-connected or fully-separate networks but not for the partially-shared network due to its design structure.}
\label{fig:trade-off2}    
\end{figure}

\noindent\textbf{Policy network design. }
We compare the performance of three policy NN designs in the ten-pool system across various neural network widths, depths, and training settings; Row 3 of Table~\ref{tab:ablation summary} summarizes the impact; Figures~\ref{fig:trade-off1}-\ref{fig:trade-off2} show the detailed results. The figures plot the percentage improvement in average cost achieved by the PPO policy (after 10 iterations) compared to the best-performing benchmark, the empirical policy, thus, a higher percentage indicates greater improvement (more desirable). We make the following observations.

First, the {partially-shared network}, represented by the solid (orange) line, consistently outperforms the other two designs across all settings. This holds even when the partially-shared network is trained with (i) fewer simulation days (1,000-5,000 days) compared to $\geq 50,000$ days for the other designs, (ii) fewer training epochs (5 epochs) compared to $\geq 10$ epochs for the other designs, and (iii) simpler NN structures (1 hidden layer with 17 neurons) versus more complex networks. These findings emphasize that the choice of policy network representation has a first-order effect—it significantly impacts algorithm performance improvement. Moreover, the results suggest that our partially-shared design effectively reduces average costs while also cutting down computational time, as it is more sample-efficient and requires a less complex NN to train.

Second, Figure~\ref{fig:trade-off1} shows that the fully-separate design is the most sensitive to increases in simulation days and training epochs among the three designs. It is because this design learns a separate policy for each epoch and requires significantly more data and longer training to achieve good performance. In contrast, Figure~\ref{fig:trade-off2} indicates that the fully-connected design benefits from increased network depth, since deeper networks better capture the complex differences across epochs while using a shared policy representation. These numerical findings align with the advantages and limitations of each design discussed in Section~\ref{Subsec: Policy representation}.

\footnotetext{When training data is very small (1,000 simulation days per actor), to achieve good performance within fewer iterations, we use the no-overflow policy as the initial policy (instead of complete-sharing in the baseline) and set the training epochs to 30 (instead of 15 in the baseline). }

\vspace{-0.4in}

\noindent\textbf{Other tailored designs. } 
Rows 4-5 of Table~\ref{tab:ablation summary} report the impact of replacing our tailored LSTD with (a) using basis functions in \cite{dai2019inpatient} and (b) using a NN-based value network, respectively. Table~\ref{tab:PPO compare VdVs} further compares the average costs of using our proposed basis functions against those in \cite{dai2019inpatient}, under different amounts of simulation data. The results show that our tailored basis functions are much less sensitive to the number of simulation days: while performance is similar with a large dataset (10,000 days), the advantage of our proposed functions becomes more evident with fewer simulation days (1,000 days), which is more practical for real-world implementation. 
Using the NN-based value network, on the other hand, shows worse final performance and requires extra time for training. This is likely due to higher variance and less stable value estimates for training the NN. See Appendix~\ref{subsubapp: NN critic} for detailed comparison.

Row 6 of Table~\ref{tab:ablation summary} compares the batching approximation and the original sequential atomic updates. The batching setup substantially reduces computational time ($\sim50\%$) while maintaining similar performance. 
These findings show the value of our tailored designs in enhancing the time efficiency of the PPO algorithm. 
This not only speeds up the algorithm significantly but also makes it more practical for real-world applications where large-scale simulation data may be limited.

\begin{table}[htp]
    \centering
\scalebox{0.95}{
    \begin{tabular}{|c|c|c|c|}
    \hline
    Basis function & 1k days per actor&10k days per actor\\
\hline
 Proposed      
 & 311.75$\pm$3.74& 302.01$\pm$1.28 \\     
 \cite{dai2019inpatient} 
& 325.71$\pm$4.07 & 307.71$\pm$1.31 \\ 
\hline  
\end{tabular}
}
\vspace{0.1in}
\caption{Comparing the average costs from PPO algorithm using different basis functions. 
}
\label{tab:PPO compare VdVs}
\end{table}

We present in Appendix~\ref{Subsec: additional NN structure results} a detailed performance comparison under various combinations of hyper-parameters (remaining rows of Table~\ref{tab:ablation summary}). The key takeaway remains consistent with our earlier discussion: our tailored PPO design captures the most salient features of the problem and is of first-order importance in reducing long-run average costs and improving computational efficiency. While tuning other hyper-parameters, such as the number of simulation days, training epochs, and policy update iterations, can further enhance performance, their impact is relatively marginal.

\subsection{Practical Implementation and Policy Interpretability}
\label{subsec: managerial insights}

In this section, we discuss how to implement our approach in practice, examine the policy interpretability since this is key for future implementation, and summarize actionable insights.

\noindent\textbf{Practical implementation.}  Since our algorithm is model-based, the hospital first needs to calibrate the queueing system to reflect specific patient flows, arrival rates, and discharge patterns. Such calibration has been successfully applied in various hospital settings, e.g, \cite{bertsimas2024hospital}. Once calibrated, managers can vary parameters such as holding and overflow costs, which act as levers to balance competing objectives: minimizing patient waiting or reducing overflow occurrences. By adjusting these values, managers can explore trade-offs and generate Pareto-optimal solutions to meet their desired operational goals. For instance, increasing overflow costs could reduce the amount of overflow, or vice versa; see, e.g., Appendix~\ref{sub-app:additional_metrics}.

Note that the PPO algorithm employs a randomized policy that provides action probabilities for patient routing decisions. Decision-makers can use these probabilities, $\pi(f\mid s)$, to guide overflow strategies, such as selecting the most likely action or using de-randomization techniques (e.g., adding a final layer to the policy network for deterministic choices). The randomized policy also offers flexibility for managers to adjust routing decisions in real time. For example, while the algorithm suggests the broad direction such as which wards or departments require the most overflow, managers can refine these suggestions based on real-world factors like patient conditions or staff availability. Ultimately, decisions regarding which individual patients to overflow will depend on the judgment of hospital managers, hospitalists, and charge nurses; our algorithm supports, rather than replaces, the overflow decision-making.

\noindent\textbf{Policy Interpretability.}
Using visualization and explainable AI techniques, we analyze the policy generated by our PPO algorithm and extract actionable insights for hospital managers. 
First, the learned policy aligns with managerial intuition in several key areas: it overflows more patients when source wards are congested and/or destination wards are underutilized; it is less aggressive with overflow during the day, anticipating discharges, and more proactive in the evening when discharges are rare; and it dynamically balances medical closeness between the primary
and overflow wards against system-wide congestion levels when selecting overflow destinations. Beyond these intuitive behaviors, the analysis reveals more subtle, non-obvious \emph{network-level} effects. 
For example, the policy reroutes patients based on the occupancy and discharge dynamics of their source wards to mitigate future congestion. It also accounts for indirect interactions among patient classes that share common overflow destinations, so congestion in one part of the network
can affect overflow decisions elsewhere.
These results demonstrate that the PPO policy not only aligns with operational intuition but also leverages system-level complexity to achieve performance gains that would be difficult to replicate via manual decision rules. Full visualizations and detailed analyses are provided in Appendix~\ref{app:shap}.

\noindent\textbf{Practical and Algorithmic Insights.}
A key strength of our approach lies in its domain adaptation—
leveraging problem-specific characteristics to develop more effective and efficient algorithms tailored in the hospital overflow context. By integrating a tailored policy network, customized basis functions, and a batching setup, our PPO algorithm requires significantly less simulation data than typical reinforcement learning or ADP-based approaches. This not only accelerates training but also makes the method practical in environments with limited data or computational resources. The case study highlights that domain-aware design choices are much more critical to algorithm performance than ad-hoc neural network tuning. Our findings offer a systematic path for adapting general-purpose RL frameworks to complex operational systems.

Beyond algorithmic contributions, our policy visualization analysis also provides actionable heuristics for hospital managers, for example, overflow intensity should vary by time of day, and decisions should depend on the relative congestion of both source and target units rather than follow rigid rules. While such principles can support better decision-making in the absence of sophisticated analytics, our study also highlights the limitations. Human-devised heuristics often fall short in capturing the anticipatory dynamics and interconnected network effects of large hospital systems. Scalable, data-driven tools like ours can systematically learn and adapt to these complexities, going beyond what ad hoc rules can achieve.

\section{Conclusion}

Overflow management becomes increasingly complex as the number of hospital wards and patient classes increases, imposing challenges that traditional methods like ADP struggle to handle at scale. The approach developed in this paper addresses the scalability issue by leveraging atomic actions and a partially-shared policy network. Additionally, our tailored designs for the hospital overflow context allow the training of our algorithm efficient and reduce the need for extensive simulation data, which overcomes a common limitation of reinforcement learning algorithms. The case study results demonstrate that our PPO algorithm matches or outperforms existing benchmarks in systems with up to twenty hospital wards.

{The atomic action + PPO framework is particularly effective in settings with the following structural features: (i) combinatorial action spaces involving multi-customer assignment decisions to multiple service pools; (ii) intertemporal congestion effects, where current decisions impact future system state and cost through queue dynamics; (iii) periodic, long-run optimization with arrival and departure patterns follow strong time-of-day or day-of-week cycles, which requires policies to adapt dynamically over time while still optimizing long-run performance. While we apply the method to inpatient overflow management, similar structures arise in other healthcare and service systems such as dynamic nurse assignment~\citep{Sun2025fluidppo}, call centers or customer support systems, emergency responder dispatching, ride-hailing~\citep{feng2021scalable}, and long-term care admissions~\citep{dong2025multiclass}. In these contexts, our approach offers a scalable alternative to full enumeration or problem-specific heuristics.}

There are several promising directions for future research building on the framework developed in this paper. One potential area is exploring the integration of patient-level contextual data, such as individual health conditions or treatment plans, to further personalize overflow decisions for better patient outcomes. Another avenue of interest is incorporating uncertainty in patient arrival rates, e.g., through online learning, to adapt to sudden surges in demand during public health crises. Moreover, our algorithm testing has thus far been on simulated systems, an essential first step in algorithm development. We are actively seeking hospital partners to deploy the algorithm in real-world settings and gather pilot data. Incorporating a human-in-the-loop component will be crucial for balancing algorithmic recommendations with human decision-making to make the decision-support adaptable to the complexities of hospital operations.


\ACKNOWLEDGMENT{%
}

%
%
%
%


\bibliographystyle{informs2014} 
\bibliography{refs} 


\newpage
\begin{APPENDICES}
\section{Notations}
\label{app:notation}
\begingroup

\renewcommand{\arraystretch}{1.25} 

\begin{table}[htp]
    \centering
\scalebox{0.85}{
    \begin{tabular}{|c|p{15cm}|}
    \hline

    $\mathcal{J}=\{1,2,\dots,J\}$ & set of pools with $J$ being the number of pools 
    \\\hline

    $N_j$ & number of servers in pool $j$, $j\in\mathcal{J}$\\\hline 

    $\lambda_j(t)$ & arrival rate function of class $j$ on a continuous-time scale\\\hline

    $\Lambda_j=\int_0^1\lambda_j(t)$ & daily arrival rate\\\hline

    $\mu_j$& probability of each customer in pool $j$ to be discharged on the current day\\\hline

    $h_{dis}$& random variable for the time of leaving on the current day in case of departure, with CDF $F_j$ for pool $j$\\\hline

    $t_0, t_1,\dots,t_{m-1}$& decision epochs on day $t$, with $m$ being the number of epochs per day \\\hline

    \begin{tabular}{c}
        $S(t_k)=(X_1(t_k),\dots,X_J(t_k)$  \\
         $,Y_1(t_k),\dots,Y_J(t_k), h(t_k))\in \mathbb{N}^{2J+1}$
    \end{tabular}& system state at decision epoch $t_k$ ($k$-th decision time-point on day $t$)\\\hline

    $X_j(t_k)$& customer count in class $j$ at decision epoch $t_k$\\\hline

    $Y_j(t_k)$&  to-depart count in pool $j$ at decision epoch $t_k$\\\hline

    $h(t_k)$&  time-of-day indicator for decision epoch $t_k$\\\hline

    $\mathcal{S}$& state space\\\hline

    $\mathcal{S}^h$& sub-space which contains all states with time-of-day being $h$\\\hline

    $f(t_k)=\{f_{i,j}(t_k), ~ i,j=1,\dots,J\}$& overflow/wait decision at decision epoch $t_k$\\\hline

    $Q_j=(X_j-N_j)^+$&  waiting count in class $j$\\\hline

    $Z_j=\max\{N_j,X_j\}$&  in-service count in pool $j$\\\hline

    $s$ &pre-action state\\\hline

    $s^+$ &post-action state\\\hline

    $p(s'|s,f)$ & one-epoch transition probability from $s$ to $s'$ under action $f$\\\hline

    $C_j$& unit holding cost of class $j$ waiting customers\\\hline

    $B_{i,j}$& unit overflow cost of assignments from class $i$ to pool $j$\\\hline

    $g(s,f)$& one-epoch cost at state $s$ under action $f$\\\hline
    
    $\pi(f\mid s)$& probability of choosing action $f$ at state $s$ according to policy $\pi\in \Pi$\\\hline

    $p_{\pi}(s'|s)=\mathbb{E}_{f\sim \pi(\cdot|s)} p(s'|s,f)$& expected transition probability from $s$ to $s'$ under policy $\pi\in\Pi$\\\hline

    $g_{\pi}(s)=\mathbb{E}_{f\sim \pi(\cdot|s)}g(s,f)$& expected one-epoch cost under policy $\pi\in\Pi$\\\hline

    $\mathbf{a}=\{a^1,a^2,\dots,a^{q}\}$& atomic action sequence at state $s$, with $a^n$ denoting the atomic action for the $n$th waiting customer \\\hline

     $\boldsymbol{\sigma}=\{\sigma^1(s),\dots,\sigma^{q}(s)\}$& an order of waiting customer's classes, with $\sigma^n(s)$ denoting the  class of the $n$th waiting customer\\\hline
    
     $s^n$&  ``atomic state'' within the atomic decision process after making the assignment for the $n-1$th customer  \\\hline

     $\mathcal{A}(s^n,\sigma^n(s))$& feasible atomic action space for $a^n$ with atomic state $s^n$ and customer class $\sigma^n(s)$\\\hline

 \end{tabular}
}
\caption{Notations in the main paper (in the order of appearance). }
    \label{tab:Notation complete order}
\end{table}

\begin{table}[htp]
    \centering
\scalebox{0.85}{
    \begin{tabular}{|c|p{15cm}|}
    \hline
    
      $\pi_\theta$& randomized atomic policy parameterized by $\theta$\\\hline

     $ \kappa_{\theta}(a^n|s^n,\sigma^n(s))$& probability for choosing each $a^n$ given the current state $s^n$ and the customer class $\sigma^n(s)$ according to policy $\pi_\theta$\\\hline

    $\pi_\theta(f\mid s)$& probability of choosing a system-level action $f$ at state $s$ under policy $\pi_\theta$\\\hline 
    
    $g_{j,k} (s|\theta)$& output of a fully-connected network with parameters $\theta$ and input $s$\\\hline

    $\gamma_\pi$& average cost under policy $\pi$\\\hline

    $v_\pi$& relative value function under policy $\pi$\\\hline

    $A_\pi,A_\eta$& advantage function under policy $\pi$ and policy $\pi_\eta$\\\hline

    $\mu_\theta$& stationary distribution under policy $\pi_\theta$\\\hline

    $r_{\theta,\eta}$&probability ratio between policies $\pi_\theta$ and $\pi_\eta$\\\hline

    $L(\theta)$& objective of PPO\\\hline

    $\text{clip}(x , 1-\epsilon , 1+\epsilon)$&clip function\\\hline

$\hat L(\theta, \mathcal{D}_\eta^T, \hat{\textbf{A}}_\eta(\mathcal{D}_\eta^T))$& estimated objective function using the simulation data collected under
policy $\pi_\eta$, with $\mathcal{D}_\eta^n$ denoting trajectory of $T$ simulated days generated under policy $\pi_\eta$, and  $\hat{\textbf{A}}_\eta(\mathcal{D}_\eta^T)$ denoting estimated advantage functions at each decision epoch of trajectory $\mathcal{D}_\eta^T$  \\\hline

     \end{tabular}
}

    \caption{Notations for atomic policy and PPO (in the order of appearance). 
    }
    \label{tab:Notation complete order-continued}
\end{table}

\endgroup

\section{Transition Dynamics under Randomized Policy}

\label{app: transition dynamics and probability}

The transition dynamics are different for non-midnight epochs and midnight epoch. 

\textit{Non-midnight epochs ($h\ne 0$).} From the assumptions on the arrivals and departures in Section~\ref{Subsec: queue}, the number of arrivals, denoted as $a_j$, is a realization from the random variable $A_j^h$, which follows $\text{Poi}\big( \int_{h}^{h'}\lambda_j(s)ds \big)$; and the number of departure, denoted as $d_j$, is a realization of the random variable $D_j^h$, which follows $\text{Bin}(y_j,p_j^h)$. Here, $p_j^h$ represents the probability of a to-depart patient who is still in hospital at epoch $h$ and will be discharged between $h$ to $h'$: 
 $$p_j^h=\frac{F_j(h')-F_j(h)}{1-F_j(h)},
 $$
where $F_j(h)$ is the CDF for the discharge time.

From the transition dynamics in Equation~\eqref{Equ: transition-current to next}, we can specify the transition probability given action $f$ at a non-midnight epoch $h$ as
 $$p(s'|s,f)=\prod_{j=1}^J\mathbb{P}(D_j^h=y_j-y_j')\mathbb{P}\left(A_j^h=x_j'-\left(x_j+\sum_{i=1,i\ne j}^J f_{i,j}-\sum_{k=1,k\ne j}^J f_{j,k}\right)+y_j-y_j'\right),\quad h=1,2,\dots,m-1.
 $$
 
\textit{Midnight epoch ($h=0$).}
Compared with non-midnight epochs, the main differences of the dynamics in midnight epoch are: (i) there is no departure between midnight epoch and the next epoch as in \cite{dai2019inpatient}; (ii) random draw of new to-depart patients who will leave in the coming day, i.e., $y_j'$ is a realization of $B_j\sim Bin(\min\{x_j,N_j\},\mu_j)$. The transition probability given action $f$ at midnight epoch follows
$$p(s'|s,f)=\prod_{j=1}^J\mathbb{P}\left(A_j^0=x_j'-x_j-\sum_{i=1,i\ne j}^J f_{i,j}+\sum_{k=1,k\ne j}^J f_{j,k}\right)\mathbb{P}(B_j=y_j'),\quad h=0.
 $$
 \label{subsec: transition dynamics under randomized policy}
For a randomized policy $\pi$ with the corresponding action probability $\pi(f\mid s)$ defined in \eqref{equ: system-level probability of atomic action policy}, the state transition probability follows  
$$p_\pi(s'|s)=\mathbb{E}_{f\sim \pi(\cdot|s)}p(s'|s,f)
$$

\medskip

\section{PPO in Periodic Setting}
\label{App: PPO Motivation}

\cite{dai2022queueing} provides theoretical guarantee of performance improvement using the PPO method in the time-stationary, long-run average cost setting. However, extending it to the setting in our paper is not straightforward due to a key challenge: periodicity. In Section~\ref{subsec: challenges}, we elaborate on this challenge. To address the periodicity, in Section~\ref{subsec: daily transition MC} we introduce the Markov chains (MC) with one period (e.g., one day) as one step under the periodic setting and redefine the objective functions accordingly. We call them as the \emph{daily MC}s to contrast with the original MC where one epoch corresponds to an hour/few hours within a day (one period contains $m$ epochs). In Section~\ref{subsec: average cost gap}, we prove that, under certain assumptions about the current policy and the updated policy,  (i) the difference between the long-run average cost of the two policies can be decomposed into two terms; (ii) the decay rate of the second term is faster than that of the first term. 
Based on these results, we show in 
Section~\ref{subsec: relate upper bound to PPO} the improvement guarantee of PPO for the daily MCs.

\subsection{The Key Challenge}
\label{subsec: challenges}

There are two differences between the setting in \cite{dai2022queueing} and our setting: (i) their per-epoch cost is action-independent (i.e., the action does not change the immediate holding cost), and (ii) the resulting DTMC under their considered policy class is irreducible and aperiodic. For (i), the overflow actions affect the per-epoch cost in our setting, which requires some different algebra in derivation. The critical difference lies in (ii). That is, our setting has time-varying periodicity in the state transitions. As a result, if we directly incorporate the epoch $h$ (hour-of-day indicator) into the state $s$, the resulting DTMC becomes periodic. Specifically, a state $s^h \in \mathcal{S}^h$ can only transition to a state of the next epoch $s^{h+1}\in\mathcal{S}^{h+1}$ (recall that $\mathcal{S}^h$ denotes the subset of state space containing states with the time-of-day indicator $h$). Therefore, the results in \cite{dai2022queueing} cannot be directly applied to our periodic environment.

\subsection{Daily MC and Revised Objective}
\label{subsec: daily transition MC}

Given any policy $\pi_\eta$ parameterized by $\eta$, we denote the transition matrix and cost vector of the original MDP, respectively, as $\textbf{P}_\eta=\{p_\eta(s|s'),s,s'\in\mathcal{S}\}$ and  $\textbf{g}_\eta=\{g_\eta(s),s\in\mathcal{S}\}$. 
The matrix $\textbf{P}_\eta$ has the following structure:
\begin{equation}
        \mathbf{P}_\eta=\left(\begin{array}{ccccc}
         0&\mathbf{P}^{0,1}_\eta&0&\cdots&0  \\
         0&0&\mathbf{P}^{1,2}_\eta&\cdots&0  \\
         \vdots&\vdots&\vdots&\ddots&\vdots\\
         0&0&0&\cdots&\mathbf{P}^{m-2,m-1}_\eta\\
         \mathbf{P}^{m-1,0}_\eta&0&0&\cdots&0
\end{array}\right),
\label{equ: block structure of periodic P}
\end{equation}
where $\mathbf{P}^{h,h'}_\eta
$ denotes the sub-matrix for transitions from states in $\mathcal{S}^{h}$ to states in $\mathcal{S}^{h'}$ under policy $\pi_\eta$. The expected cost vector can be written as 
$$\mathbf{g}_\eta=\left(
    \begin{array}{c}
        \mathbf{g}^0_\eta\\
        \vdots\\
        \mathbf{g}^{m-1}_\eta
    \end{array}\right),
$$
where $\mathbf{g}^h_\eta
$ denotes the one-epoch cost vector for states in $\mathcal{S}^h$.

To address the challenge of periodicity in the original MDP setting, we introduce the \emph{daily MC}s, which are aperiodic. Under a given policy $\pi_\eta$, for each time-of-day indicator $h=0,1,...,m-1$, if we observe states only at epoch $h$ of each day, the resulting stochastic process $\{S(t_h), t=0,1,...\}$ is still an MC. The corresponding \emph{daily} transition matrix and expected cost vector are defined as
\begin{equation}
    \begin{aligned}
        \tilde {\mathbf{P}}^h_\eta&=\mathbf{P}^{h,h+1}_\eta\mathbf{P}^{h+1,h+2}_\eta\cdots\mathbf{P}^{h-1,h}_\eta,\\
      \tilde {\mathbf{g}}^h_\eta&=\mathbf{g}^h_\eta+\mathbf{P}^{h,h+1}_\eta\mathbf{g}^{h+1}_\eta+\cdots+ \mathbf{P}^{h,h+1}_\eta\mathbf{P}^{h+1,h+2}_\eta\cdots\mathbf{P}^{h-2,h-1}_\eta\mathbf{g}^{h-1}_\eta.
    \end{aligned}
    \label{equ: def for one-day Markov chain under pi_eta}
\end{equation}
Without loss of generality, we restrict policies such that the induced daily MCs are irreducible and aperiodic. This allows us to use the proof framework in \cite{dai2022queueing} to extend the results in time-stationary settings to our daily MC. 

To start, the relative value function $v^h_\eta$ under $\pi_\eta$ is given by 
\begin{equation}
    v^h_\eta(s):=\mathbb{E}\Big[ \sum_{t=1}^\infty \Big( \tilde{g}_\eta^h(S(t_h))-(\mu_\eta^h)^T\tilde{\mathbf{g}}_\eta^h \Big) ~|~ S(1_h)=s \Big], ~\forall s\in\mathcal{S}^h. 
    \label{equ: relative vf of h}
\end{equation}  
We redefine the PPO objective as: 
\begin{equation}
L^h(\theta) := \mathop{\mathbb{E}}\limits_{\substack{ s\sim\mu_\eta^h\\ f\sim\pi_\eta(\cdot|s)}} \max \left\{ r_{\theta,\eta}(f\mid s) A^h_\eta(s, f),
\quad \text{clip}\left( r_{\theta,\eta}(f\mid s) , 1 -\epsilon , 1 + \epsilon \right) A^h_\eta(s,f) \right\},
    \label{equ: PPO obj for one epoch h}
\end{equation}
where $A^h_\eta$ is the advantage function of the daily MC, whose definition will be specified later in this section; the vector $\mu_\eta^h$ is the stationary distribution of this daily MC; and $r_{\theta,\eta}$ is the policy ratio. We update the policy from $\pi_\eta$ by minimizing this redefined objective $L^h$. Moreover, we update the policy only at states in $\mathcal{S}^h$. We denote the updated policy as $\pi_{\theta,h}$, which follows 
 \begin{equation}
        \pi_{\theta,h}(f\mid s)=\left\{
        \begin{array}{cc}
          \pi_\theta(f\mid s),   & \forall s\in\mathcal{S}^h, f\in\mathcal{A}(s), \\
           \pi_\eta(f\mid s),  & \forall s\in\mathcal{S}\setminus\mathcal{S}^h, f\in\mathcal{A}(s)
        \end{array}
        \right..
    \end{equation}
We denote the one-day transition matrix and one-day cost vector of the induced daily MC, under this policy update scheme, as $\tilde{\mathbf{P}}_{\theta}^{h}$ and $\tilde{\mathbf{g}}_{\theta}^{h}$, respectively. They follow 
\begin{equation}
    \begin{aligned}
&\tilde{\textbf{P}}_{\theta}^h=\textbf{P}_\theta^{h,h+1}\textbf{P}_\eta^{h+1,h+2}\cdots\textbf{P}_\eta^{h-1,h}, 
    \end{aligned}
    \label{equ: P for daily MC for h}
\end{equation}
and 
\begin{equation}
    \begin{aligned}
\tilde{\textbf{g}}_\theta^h   =     & \textbf{g}_\theta^h+\textbf{P}_\theta^{h,h+1}\textbf{g}_\eta^{h+1}
        +
        (\textbf{P}_\theta^{h,h+1}\textbf{P}_\eta^{h+1,h+2})\textbf{g}_\eta^{h+2}
        +\cdots +(\textbf{P}_\theta^{h,h+1}\textbf{P}_\eta^{h+1,h+2}\cdots\textbf{P}_\eta^{h-2,h-1})\textbf{g}_\eta^{h-1}.  
    \end{aligned}
    \label{equ: g for daily MC for h}
\end{equation}
Note that the transition from epoch $h$ to $h+1$ is given by $\textbf{P}_\theta^{h,h+1}$ with $\theta$, while all other transitions remained to be $\textbf{P}^{\cdot, \cdot}_\eta$ with $\eta$. For notational simplicity, we use $\theta$ here in the subscript but emphasize that we only update the policy at epoch $h$. 

Let $\tilde g_\eta^h(s,f)$ represent the one-day expected cost starting from state $s\in\mathcal{S}^h$, given that action $f$ is taken at $s$ and subsequent actions are determined by the policy $\pi_\eta$. 
The advantage function for the daily MC (under our considered policy update) is given by: 
\begin{equation}
A_\eta^h(s,f)=\tilde{g}_\eta^h(s,f)-(\mu_\eta^h)^T\tilde{\mathbf{g}}_\eta^h+\mathbb{E}_{s'\sim \tilde{p}_\eta^h (\cdot|s,f)}[v_\eta^h(s')]-v_\eta^h(s),
    \label{equ: def adv one-day h}
\end{equation}
where $\tilde p_\eta^h(s'|s,f)$ denotes the one-day transition probability from state $s\in\mathcal{S}^h$ in the current day to state $s'\in\mathcal{S}^h$ in the next day, given that action $f$ is taken at $s$ and subsequent actions are determined by  $\pi_\eta$.

\subsection{Average cost gap}
\label{subsec: average cost gap}

We denote the \emph{unclipped} PPO objective function as $N_1^h(\theta,\eta)$, which follows
\begin{equation}
    \begin{aligned}
N_1^h(\theta,\eta):=&\mathop{\mathbb{E}}\limits_{\substack{ s\sim\mu_\eta^h
\\ f\sim\pi_\eta(\cdot|s)}}
    [r_{\theta,\eta}(f\mid s) A_\eta^h(s,f)]
    = 
\mathop{\mathbb{E}}\limits_{\substack{ s\sim\mu_\eta^h\\ f\sim\pi_\theta(\cdot|s)}}
[ A_\eta^h(s,f)]
\\ =&(\mu_\eta^{h})^T \Big( \tilde{\mathbf{g}}_{\theta}^{h}-(\mu_\eta^{h})^T\tilde{\mathbf{g}}_\eta^{h}\mathbf{e}+(\tilde{\mathbf{P}}_{\theta}^{h}-I)\mathbf{v}_\eta^{h} \Big) .  
    \end{aligned}
    \label{Equ: relate N1 to adv h}
\end{equation}
Our eventual goal is to examine whether minimizing $N_1^h$ while controlling $r_{\theta,\eta}$ to be near 1 can guarantee policy improvement. To achieve this goal, we study the average-cost gap between the current policy $\pi_\eta$ and the new policy $\pi_{\theta,h}$. We impose the following assumptions for the stability of the current policy $\pi_\eta$; these assumptions are from~\cite{dai2022queueing}. 
\begin{assumption}
    For any given $h=0,1,...,m-1$, 
    \begin{itemize}
        \item the daily MC with transition matrix $\tilde {\mathbf{P}}^h_\eta$ is irreducible and aperiodic. 
        \item there exists some vector $\mathcal{V} = \{\mathcal{V}(s), ~s\in\mathcal{S}^h\}\ge 1$, some constants $b\in(0,1)$, $d\ge 0$ and finite subset $C\subset \mathcal{S}^h$ satisfying   
      \begin{equation}
            \tilde{\mathbf{P}}^h_\eta \mathcal{V}\le b \mathcal{V}+d \mathbf{1}_{C},
        \end{equation}
        where each element of the vector $\mathbf{1}_{C}=\{\mathbbm{1}_{C}(s)\in\{0,1\},s\in\mathcal{S}^h\}$ is $\mathbbm{1}_{C}(s)=1$ if $s\in C$ and $0$ if $s\not\in C$.
        \item the one-day cost vector satisfies $|\tilde {\mathbf{g}}^h_\eta|\le\mathcal{V}$.
    \end{itemize}
    \label{ass: drift}
\end{assumption}

Under these conditions, one can prove that there exists a unique stationary distribution $\mu_\eta^h$ for the MC with transition matrix $\tilde{P}_\eta^h$; see for example, Theorem 11.3.4 and Theorem 14.3.7 in \cite{meyn2012markov}, which is also restated in Lemma 1 of \cite{dai2022queueing}.
Moreover, for any $h=0,1,...,m-1$, we assume that the new policy (using the considered updating mechanism) satisfies the following condition: 
\begin{equation}
    \|(\tilde{\mathbf{P}}_{\theta}^{h}-\tilde{\mathbf{P}}_{\eta}^{h})\sum_{n=0}^{\infty}(\tilde{\mathbf{P}}_{\eta}^{h}-\Pi_\eta^{h})^n\|_{\mathcal{V}}<1,
    \label{equ: ass for theta}
\end{equation}
where every row of $\Pi_\eta^{h}$ equals the stationary distribution of the MC with transition matrix $\tilde{\mathbf{P}}_\eta^{h}$, i.e.,  $\Pi_\eta^{h}(s,s')=\mu_\eta^{h}(s')$, and the $\mathcal{V}$-norm of a matrix $\Omega \in \mathcal{S}^{h}\times\mathcal{S}^{h}$ is defined as 
\begin{equation}
    \|\Omega\|_{\mathcal{V}}=\sup_{ s\in\mathcal{S}^h} \frac{1}{\mathcal{V}(s)} \sum_{ s'\in\mathcal{S}^{h}}|\Omega(s,  s')|\mathcal{V}( s' ). 
    \label{equ: V-norm matrix def-1}
\end{equation} 
Then one can show that the MC with transition matrix $\tilde{\mathbf{P}}_{\theta}^{h}$ also has a unique stationary distribution $\mu_{\theta}^{h}$. For the proof, see Lemma 4 of \cite{dai2022queueing}.

In the following proposition, we decompose the difference in the long-run average cost between the current policy $\pi_\eta$ and the updated policy $\pi_{\theta,h}$ into two terms, and analyze their decay rates when $\theta$ approaches $\eta$.  
\begin{proposition}
Assume that the current policy satisfies Assumption~\ref{ass: drift} and the new policy  satisfies~\eqref{equ: ass for theta}. The difference between the long-run average costs of the two policies equals 
$$(\mu_{\theta}^{h})^T\tilde{\textbf{g}}_{\theta}^{h}-(\mu_\eta^h)^T\tilde{\textbf{g}}_{\eta}^h=N_{1}^{h}(\theta,\eta)+N_{2}^{h}(\theta,\eta) , 
$$
where $N_1^h$ is defined in Equation~\eqref{Equ: relate N1 to adv h}, and $N_2^h$ is defined as
\begin{equation}
    N_{2}^{h}(\theta,\eta):=(\mu_{\theta}^{h}-\mu_\eta^{h})^T \Big( \tilde{\mathbf{g}}_{\theta}^{h}-(\mu_\eta^{h})^T\tilde{\mathbf{g}}_\eta^{h}\mathbf{e}+(\tilde{\mathbf{P}}_{\theta}^{h}-I)\mathbf{v}_\eta^{h} \Big) . 
\end{equation}
Moreover, we have $ N_1^h(\theta,\eta)=O(||\textbf{r}^h_{\theta,\eta} - 1||_\infty )$ and  $N_2^h(\theta,\eta)=O(||\textbf{r}^h_{\theta,\eta} - 1||_\infty ^2)$, where 
$$\big\|\textbf{r}^h_{\theta,\eta}-1\|_\infty:=\sup_{s\in\mathcal{S}^h}\sum_{f\in\mathcal{A}(s)}|
r_{\theta,\eta}(f\mid s)-1 \big\|.
$$ 
\label{prop: periodic average cost upper bound}
\end{proposition} 
The proof of this proposition is in Appendix~A of the Technical Companion \citep{tech2025comp}.

\subsection{Improvement Guarantee of PPO}
\label{subsec: relate upper bound to PPO}

According to Proposition~\ref{prop: periodic average cost upper bound}, if we update policy from $\pi_\eta$ to $\pi_{\theta,h}$, the change in the average cost equals  $N_1^h(\theta,\eta)+N_2^h(\theta,\eta)$. Then, the negativity of $N_1^h(\theta,\eta)+N_2^h(\theta,\eta)$ guarantees that the new policy yields an improved performance comparing with the current policy $\pi_\eta$ in terms of the average daily cost (where this average daily cost is the same regardless the epoch we are considering). Since  
$N_1^h(\theta,\eta)+ N_2^h(\theta,\eta)\le N_1^h(\eta,\eta)+ N_2^h(\eta,\eta)=0$, to achieve maximum improvement, we can pick $\theta=\theta^*$ as
\begin{equation}
    \theta^*=\arg\min _{\theta\in\Theta} N_1^h(\theta,\eta)+ N_2^h(\theta,\eta) . 
\label{Equ: original obj}
\end{equation}  
However, this optimization problem cannot be directly solved because $N_2^h$ depends on the stationary distribution of the new policy, which has no closed form and cannot be estimated when $\theta$ is to be determined. Therefore, the PPO framework proposes using an alternative way to approximately solve problem~\eqref{Equ: original obj}. That is, obtain $\theta$ via the the objective function $L^h(\theta)$, as defined in Equation~\eqref{equ: PPO obj for one epoch h}. 

Proposition~\ref{prop: periodic average cost upper bound} implies that, as the policy ratio $r^h_{\theta,\eta}$ is close to 1, $N_2^h(\theta,\eta)$ is of a smaller order compared to $N_1^h(\theta,\eta)$. 
Therefore, if we get $\theta$ via (i) minimizing $N_1^h(\theta,\eta)$, the unclipped objective function and (ii) keeping $r_{\theta,\eta}(f\mid s)$ to be close to 1, then policy improvement is guaranteed as $N_1^h(\theta,\eta)<0$ and $||\textbf{r}^h_{\theta,\eta} - 1||_\infty$ is sufficiently small so that
$N_1^h(\theta,\eta)+N_2^h(\theta,\eta)<0$. These two goals are achieved simultaneously through the conservative updates induced by minimizing the clipped objective function $L^h(\theta)$, as we explain below.

Note that in adapting to the daily MC, $L^h(\theta)$ is not equivalent to \eqref{Equ: PPO Obj} unless $m=1$. In our algorithm implementation, we nonetheless adopt \eqref{Equ: PPO Obj}, since our primary objective here is to illustrate the rational behind PPO, namely, its clipped objective effectively controls the two terms $N_1$ and $N_2$. 
For notational simplicity, we now explain how the max operator and clipping function in \eqref{Equ: PPO Obj} work together to enforce conservative policy updates (dropping the index $h$). 
We can rewrite~\eqref{Equ: PPO Obj} in the following explicit piecewise form:
\begin{align*}
    &\max \left\{ r_{\theta,\eta}(f\mid s) A_\eta(s, f),
\quad \text{clip}\left( r_{\theta,\eta}(f\mid s) , 1 -\epsilon , 1 + \epsilon \right) A_\eta(s,f) \right\}\\
=&\left\{ \begin{array}{cc}
  (1+\epsilon)A_\eta(s,f),   &\text{if } A_\eta(s,f)<0 \text{ and } r_{\theta,\eta}(f|s)>1+\epsilon,  \\
 (1-\epsilon)A_\eta(s,f),   &\text{if } A_\eta(s,f)>0 \text{ and } r_{\theta,\eta}(f|s)<1-\epsilon,  \\
 r_{\theta,\eta}(f|s)A_\eta(s,f), &\text{otherwise}.
\end{array}\right.
\end{align*}
In the first case, $A_\eta(s,f)<0$ indicates that action $f$ reduces cost (desirable), so the policy update encourages increasing $r_{\theta,\eta}(f|s)$; the clipping at $1+\epsilon$ prevents an excessively large increase. In the second case, $A_\eta(s,f)>0$ indicates that action $f$ increases cost (not desirable), and the policy update encourages reducing $r_{\theta,\eta}(f|s)$; when $r_{\theta,\eta}(f|s)<1-\epsilon$, the clipping term constrains this decrease. In the remaining case, the unclipped objective is used.

\section{Basis Function Design}
\label{App: pool-wise decomposition}
\begingroup
\renewcommand{\arraystretch}{1.25}
\begin{table}[htp]
    \centering
\scalebox{0.85}{
    \begin{tabular}{|c|p{15cm}|}
    \hline
     
      $s_j=(x_j,y_j,h)\in \mathcal{S}_j$   & pre-action state of pool $j$ with its state space as $\mathcal{S}_j$ \\\hline
      
      $s_j^+$&post-action state of pool $j$\\\hline
      
      $s_j'$&pre-action state of pool $j$ at the next epoch\\\hline
      
      $s_j^{'+}$&post-action state of pool $j$ at the next epoch\\\hline

       $s_{-j}\in \mathcal{S}_{-j}$& pre-action state of pools excluding  pool $j$ with its state space as  $ \mathcal{S}_{-j}$\\\hline
       
      $f'=(f_{i,j}',i,j=1,\dots,J\})$& system-level overflow action of the next epoch\\\hline

      $\bar f=(\bar f_{i,j},i,j=1,\dots,J\})$& approximation of $f'$\\\hline
      
      $\pi_j$ &
      inflow and outflow policy of decomposed pool-$j$ system under policy $\pi$
    \\\hline

    $\bar\kappa_j(i|s_j,j)$& approximation of overflow action probability at the next epoch under policy $\pi$ which only uses the current post-action state information for the pool-$j$ system\\\hline

       $p_j(s_j'|s_j^+)$  & 
            transition probability of pool $j$ from post-action state $s_j^+$ to next post-action state $s_j'$\\\hline
       
       $p(s'|s^+,0)$& 
      transition probability from post-action state $s^+$ to next post-action state $s'$, with $p(s'|s^+,0)=\prod_{j=1}^J p_j(s_j'|s_j^+)$ 
      \\\hline

       $\bar v_\pi$&
       post-action value function of the original $J$-class $J$-pool system under overflow policy $\pi$ 
  \\\hline
       
       $v_{\pi,j}$& pre-action value function of pool $j$ under policy $\pi$\\\hline
       
       $\bar v_{\pi,j}$& post-action value function of of pool $j$ under policy $\pi$\\\hline
       
       $\hat v_{\pi,j}$& 
      post-action value function of the approximated decomposed pool-$j$ system under the approximated pool-independent overflow probabilities $\bar\kappa_j$ \\ \hline

       $V_d(s) = \sum_{j=1}^J \hat v_{\pi,j}(s_j)$& new design that tailors to the randomized policy\\\hline
       
       $\hat \gamma_{\pi,j}$& 
       average cost of the approximated decomposed pool-$j$ system under the approximated pool-independent overflow probabilities $\bar\kappa_j$ \\\hline

    \end{tabular}
}

    \caption{Notations in Appendix~\ref{App: pool-wise decomposition}}
    \label{tab: Notations in Appendix decomp}
\end{table}
\endgroup

In Section~\ref{Subsec: policy evaluation} of the main paper, we propose a new basis function $V_d(s)=\sum_{j=1}^J \hat v_{\pi,j}(s_j)$ that decomposes by each pool $j$ with $s_j=(x_j,y_j,h)$. Here, each $\hat v_{\pi,j}(s_j)$ denotes the \emph{post-action} value function of a single-pool system which approximates the dynamics of pool $j$ under policy $\pi$, and it is computed using the post-action Poisson equation for the decomposed pool-$j$ system in Equation~\eqref{Equ: decomposed Poisson equation-main}. To highlight the difference between the pre- and post-action states in this section, we use $s^+_j=(x_j^+,y_j^+,h)$ to denote the post-action state for the pool-$j$ system and restate Equation~\eqref{Equ: decomposed Poisson equation-main} as follows: 
\begin{equation} 
\hat v_{\pi,j}(s_j^+) = C_j q_j^+ -\hat \gamma_{\pi,j}+\mathbb{E} \left[\sum_{i=1}^J B_{j,i}\bar f_{j,i}+\hat v_{\pi,j}(s_j^{'+})\right],  
\label{Equ: decomposed Poisson equation-restate}
\end{equation}
where $q_j^+ = (x^+_j - N_j)^+$ denotes the post-action queue length, and $s_j^{'+}$ denotes the next post-action state. In this section, we provide detailed explanations on (i) how the pool-wise decomposition is performed, and (ii) how to estimate the parameters in Equation~\eqref{Equ: decomposed Poisson equation-restate} to approximate the dynamics in the original (non-decomposed) system. Table~\ref{tab: Notations in Appendix decomp} summarizes the notations used in this section.

\subsection{Pool-wise Decomposition}
\label{sub-app: pool-wise decomp}
In order to obtain pool-wise decomposed Poisson equation~\eqref{Equ: decomposed Poisson equation-restate}, we first derive the post-action Poisson equation for the original $J$-pool system. We then derive the decomposed version. 

\noindent\textbf{Post-action Poisson equation. } 
Under a given randomized, stationary policy $\pi$, the pre-action Poisson equation follows 
\begin{equation}
    \begin{aligned}
    v_\pi(s)=&\mathbb{E}_{f\sim \pi(\cdot|s)} \left[ g(s, f)-\gamma_\pi+\sum_{s'\in \mathcal{S}}p(s'|s,f) v_\pi(s') \right] 
\\
    =&\mathbb{E}_{f\sim \pi(\cdot|s)}\left[ \sum_{i,j=1}^JB_{i,j} f_{i,j} +\sum_{j=1}^J C_jq_j^{+} -\gamma_\pi +\mathbb{E}_{s'\sim p(\cdot|s^+)}v_\pi(s') \right], 
\end{aligned}
\label{Equ: Poisson equation expectation}
\end{equation}
where $s^+ = \{s^+_j\}_{j\in\mathcal{J}}$ denotes the current post-action state in the original $J$-pool system, $q_j^+$ is the post-action queue length for pool $j$, and $p(s'|s^+)$ denotes the transition probability from $s^+$ to the next pre-action state $s'$. Note that here, the transitions only depend on the arrivals and departures that occurred between the current and the next epochs.

We denote the post-action value function as $\bar{v}_\pi$, where 
\begin{equation}
    \bar v_\pi(s^+)=\sum_{j=1}^J C_jq_j^++\mathbb{E}_{s'\sim p(\cdot|s^+)}v_\pi(s'), 
    \label{Equ: def. post-action value function by pre-action value function}
\end{equation} 
using which we can rewrite (\ref{Equ: Poisson equation expectation}) as  
\begin{align*} 
    v_\pi(s) 
    = & - \gamma_\pi+\mathbb{E}_{f\sim \pi(\cdot|s)}\left[\sum_{i,j=1}^J B_{i,j}f_{i,j}+\bar v_\pi(s^+) \right]. 
\end{align*}
Plugging the above back to Equation~\eqref{Equ: def. post-action value function by pre-action value function}, we get the following post-action Poisson Equation: 
\begin{equation}
    \begin{aligned}
    \bar v_\pi(s^+)
    = & \sum_{j=1}^J C_j q_j^+-\gamma_\pi +\mathbb{E}_{s'\sim p(\cdot|s^+)}\mathbb{E}_{f'\sim \pi(\cdot|s')}\left[\sum_{i,j=1}^J B_{i,j}f_{i,j}'+\bar v_\pi(s'^+)\right], 
\end{aligned}
\label{post-action Poisson equation}
\end{equation}
with the pair $(\bar{v}_\pi,\gamma_\pi)$ being its solution.

\noindent\textbf{Decomposition. }
In the post-action Poisson Equation~\eqref{post-action Poisson equation}, the first two terms $\sum_{j=1}^J (C_j q^{+}_j-\gamma_{\pi,j})$ can be naturally decomposed by each pool, where $\gamma_{\pi,j}$ corresponds to the average cost of pool $j$ under policy $\pi$. We then focus on considering the decomposition for the terms within the expectation, i.e., 
\begin{align*}
\mathbb{E}_{s'\sim p(\cdot|s^+)}\mathbb{E}_{f'\sim \pi(\cdot|s')}\left[\sum_{i,j=1}^J B_{i,j}f_{i,j}'+\bar v_\pi(s'^+)\right]
= 
\sum_{j=1}^J\mathbb{E}_{s'\sim p(\cdot|s^+)}\mathbb{E}_{f'\sim \pi(\cdot|s')}\left[\sum_{i=1}^J B_{j,i}f'_{j,i}+\bar v_{\pi,j}(s_j'^+)\right]. 
\end{align*} 
We can observe here that the key is to decompose the probabilities, i.e., the transition probabilities $p$ and overflow action probabilities $\pi$, by each pool. Note that the transition probability  $p(s'|s^+,0)$ from post-action state $s^+$ to state of next epoch $s'$ can be decomposed by each pool naturally, that is,   
\begin{equation}
    p(s'|s^+,0)
= \prod_{j=1}^J  p_j(s_j'|s_j^+), 
\label{equ: transition matrix pool-wise decomp}
\end{equation}
where $p_j(s_j'|s_j^+)$ denotes the transition, for each pool $j$, from the post-action state to the next epoch pre-action state. This decomposition holds because after the overflow assignments, the state transitions only depend on the arrivals to and departures from each pool $j$, which are independent among $j$'s.

Next, let $s_{-j}^+=\{(x_i^+,y_i^+,h)\}_{i\ne j}$ denote the current post-action state of the other pools (excluding pool $j$) and $\mathcal{S}_{-j}$ the corresponding state space. Then, the $j$-th component within the expectation term can be rewritten as 
\begin{align}
\mathbb{E}_{s'\sim p(\cdot|s^+)}\mathbb{E}_{f'\sim \pi(\cdot|s')}\left[\sum_{i=1}^J B_{j,i}f'_{j,i}+\bar v_{\pi,j}(s_j'^+)\right]
= & \sum_{s_j'\in\mathcal{S}_j} p_j(s_j'|s_j^+)\mathbb{E}_{f'\sim \pi_j(\cdot|s_j',s_{-j}^+)}\left[ \sum_{i=1}^J B_{j,i} f'_{j,i} + \bar v_{\pi,j}(s_j'^+)\right]
\label{eq:exp-decomp-exact}
\end{align}
for any $j=1, \dots, J$, with 
$$\pi_j(f'|s_j',s_{-j}^+)=\sum_{s_{-j}'\in \mathcal{S}_{-j}} \prod_{i\ne j} p_i(s'_i|s_i^+)\pi(f' | s') \quad \textrm{for } s' = (s'_j, s'_{-j}). 
$$
The important step for this decomposition lies in the $\pi_j(f'|s_j',s_{-j}^+)$ term, which 
evaluates the marginal probability of taking overflow action $f'$ under the policy $\pi$ by taking an expectation over all the possible transition from $s_{-j}^+$ to $s'_{-j}$. This expectation is needed because in the decomposed equation~\eqref{eq:exp-decomp-exact} we only evaluated the transition for pool $j$ via $p_j(s'_j|s_j^+)$ and did not evaluate the transitions for other pools.

To fully decompose, we still need to approximate this marginal probability $\pi_j(f'|s_j',s_{-j}^+)$ since it requires information on the states of other pools besides pool $j$ (to get the probability of prescribing action $f'=\{f'_{i,j}\}$ based on policy $\pi$). To do so, we go back to the atomic action setup and approximate the probability of prescribing $f'_{i,j}$ with $\bar f_{i,j}\sim Poi(\lambda_i \bar \kappa_{j}(j|s'_j,i))$, and the probability of prescribing $f'_{j,i}$ with $\bar f_{j,i} \sim Poi(\lambda_j \bar \kappa_{j}(i|s'_j,j) )$, using some approximated probability $\bar \kappa_{j}(\cdot|s'_j,\cdot)$ that does not require any state information from other pools beyond pool $j$. Note that the probability $\bar \kappa_{j}(j|s'_j,i)$ is the probability for the \emph{atomic} action for assigning one patient from class $i$ to $j$. Under the batched routing setting with the approximate probability $\bar \kappa_{j}(j|s'_j,i)$, the aggregate action $\bar{f}_{i,j}$ becomes equivalent to the Poisson random variables defined here using the Poisson thinning property given state $s'_j$. The rational for $\bar{f}_{j,i}$ is the same.

Given these distributions of $\bar{f}$'s, for each pool, we could then approximate the remaining state transition in~\eqref{eq:exp-decomp-exact} -- from $s'_j$ to $s^{'+}_j$ -- by properly adjusting the inflows and outflows via: (a) overflow assignments from other classes to pool $j$, $\bar{f}_{i,j}$; (b) the overflow (diversion) from class $j$ to other pools, $\bar{f}_{j,i}$, i.e., $x_j^{'+}=x_j'+\sum_{i=1}^J \bar{f}_{i,j} - \sum_{i=1}^J \bar{f}_{j,i}$.
This eventually leads to our (approximated) pool-wise decomposed post-action Poisson equation: 
\begin{equation}
     \begin{aligned}
     \hat v_{\pi,j}(s_j^+) 
&=C_j q_j^+-\gamma_{\pi,j}
+
\mathbb{E}_{s_j'\sim p_j(\cdot|s_j^+)}
\left[
\sum_{i=1}^J B_{j,i}\bar f_{j,i} + \mathbb{E}_{\bar{f}}\left[\hat v_{\pi,j}(s_j^{'+})\right]
\right] , 
    \end{aligned}
    \label{Equ: decomposed Poisson equation}
\end{equation} 
Since (\ref{Equ: decomposed Poisson equation}) is derived from decomposing~\eqref{post-action Poisson equation} by pools, then summarizing over all $j$'s, the summation $\left( \sum_{j=1}^J \hat{v}_{\pi,j},\sum_{j=1}^J \hat\gamma_{\pi,j}\right)$ can serve as an approximation for the solution to~\eqref{post-action Poisson equation}. Therefore, we could use $V_d=\sum_{j=1}^J \hat{v}_{\pi,j}$ as one of the basis functions to approximate relative value functions.

\section{Additional Numerical Results}
\label{app:case-study}

For the baseline five-pool model, we use the same setting of feasible routes as in \cite{dai2019inpatient}, and the detailed description of the feasible routes in ten- and twenty-pool models are specified in Section~\ref{sec:case-model-setting}. Figure~\ref{fig: arr dep patterns} shows the time-varying arrival and discharge patterns of the five-pool model. 
The ten- and twenty-pool models can be considered as combination of two and four baseline five-pool settings, respectively, with identical discharge patterns and proportionally adjusted arrival rates. Detailed specifications are also given in Section~\ref{sec:case-model-setting}.  
The number of servers per pool in the baseline five-pool model is $(N_1,\dots,N_5)=(60,64,67,62,62)$. For the ten-pool model, we set $(N_1,\dots,N_{10})=(39,43,46,41,41,81,85,88,83,83)$, and for the twenty-pool model, we set $(N_1,\dots,N_{20})=(32,36,39,34,34,74,78,81,76,76,46,50,53,48,48,88,92,95,90,90)$.

\begin{figure}[htp]
\centering  
\scalebox{1}{
\subfigure[Hourly arrival rate]{   
\centering    
\includegraphics[width=0.4\linewidth]{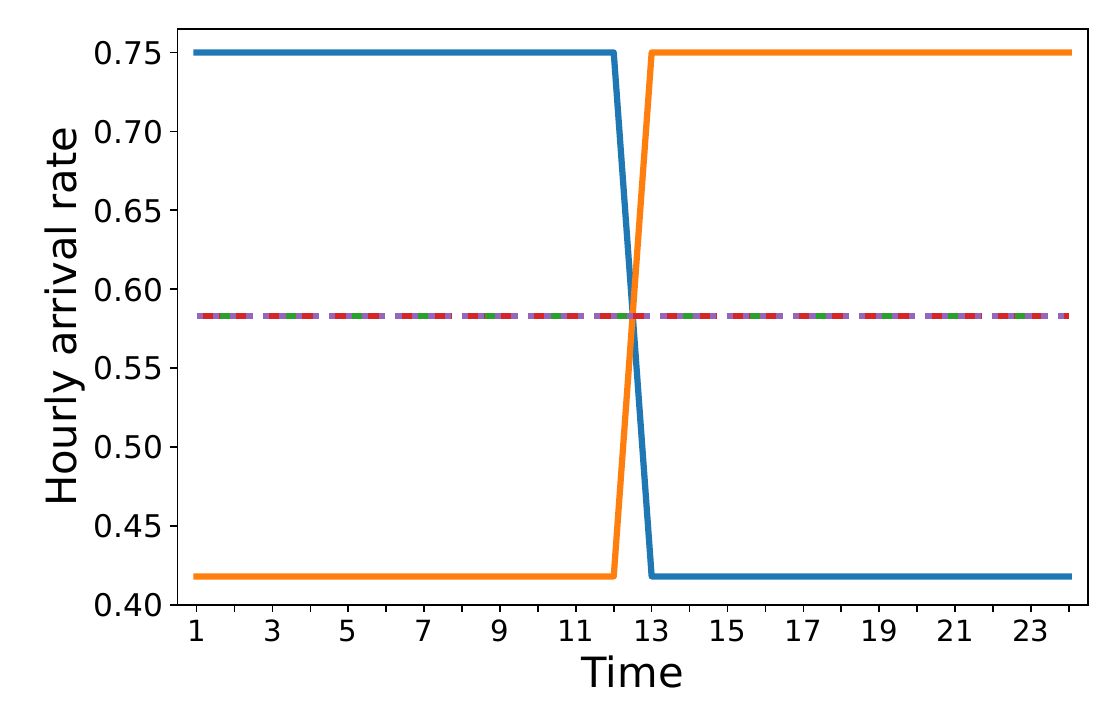}  
}
\subfigure[Hourly discharge distribution]{ 
\centering    
\includegraphics[width=0.4\linewidth]{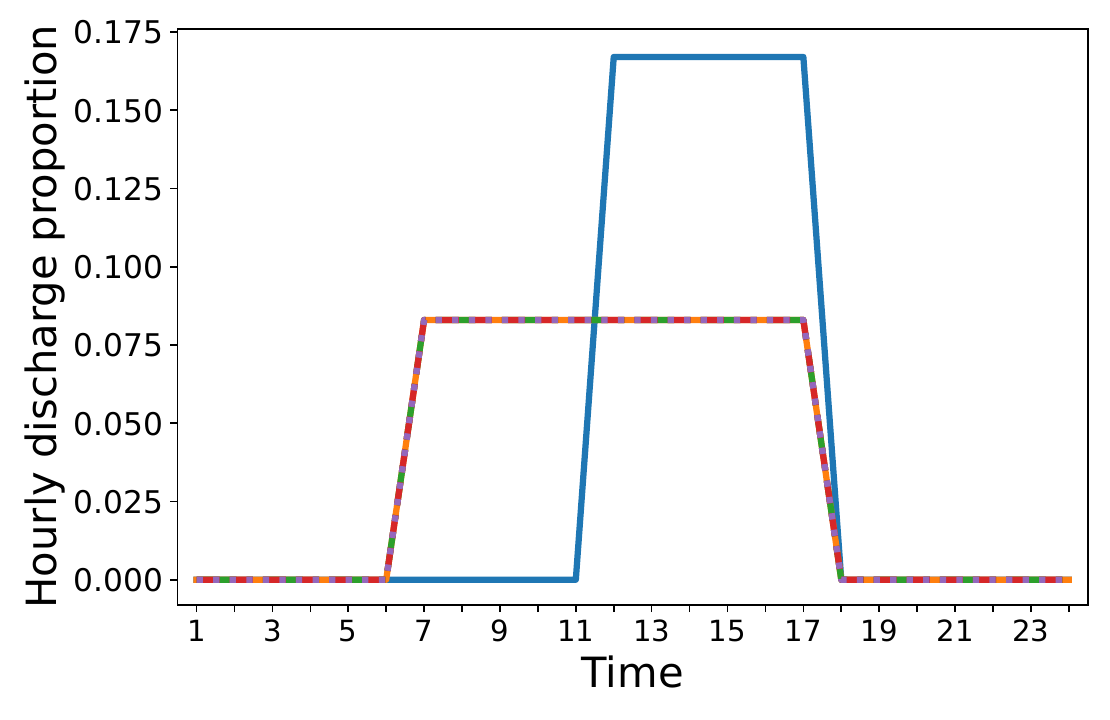}
}
\raisebox{1cm}{
\subfigure{
\centering
\includegraphics[width=0.1\linewidth]{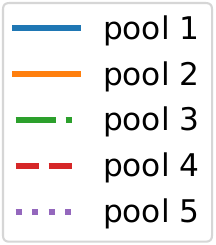} 
}}}
\caption{Time-varying arrival and discharge patterns in five-pool models.
}    
\label{fig: arr dep patterns}    
\end{figure}

\subsection{Impact of Hyper-parameters for Neural Network Training}
\label{Subsec: additional NN structure results}

Tables~\ref{tab: policy network comparison 1 in 10pool8epoch} and~\ref{tab: policy network comparison 2 in 10pool8epoch} report the numeric results for those plotted in Figures~\ref{fig:trade-off1} and \ref{fig:trade-off2}, conducted in the baseline ten-pool model. 
The ten-pool model is too large to be solved by ADP. Thus, we use the best of the three heuristic policies as the benchmark policy (with an average cost of 393.36). 

These results also offer guidance for tuning hyperparameters. Table~\ref{tab: policy network comparison 1 in 10pool8epoch} shows that comparable performance gains can be achieved by either increasing the number of simulation days or training epochs. However, increasing simulation days incurs significantly higher computational cost: not only due to longer data generation time but also because larger datasets slow down policy network training. Therefore, we recommend increasing the number of training epochs first. Table~\ref{tab: policy network comparison 2 in 10pool8epoch} further suggests that increasing network width generally yields better performance than increasing depth; in some cases, deeper networks may even degrade performance. Thus, we advise prioritizing wider architectures over deeper ones.

\begin{table}[htp]
\centering
 \scalebox{0.65}{\begin{tabular}{|cc|ccc|ccc|ccc|}
\hline
\multicolumn{2}{|c|}{(34 neurons per hidden layer)}                                                                         & \multicolumn{3}{c|}{0 depth}                                         & \multicolumn{3}{c|}{1 depth}                                             & \multicolumn{3}{c|}{2 depth}                                         \\ \hline
\multicolumn{2}{|c|}{\begin{tabular}[c]{@{}c@{}}\# simulation days \\ per actor\end{tabular}}              & \multicolumn{1}{c|}{5k}     & \multicolumn{1}{c|}{10k}      & 100k   & \multicolumn{1}{c|}{5k}       & \multicolumn{1}{c|}{10k}      & 100k     & \multicolumn{1}{c|}{5k}     & \multicolumn{1}{c|}{10k}      & 100k   \\ \hline
\multicolumn{1}{|c|}{\multirow{2}{*}{\begin{tabular}[c]{@{}c@{}}Full \\ separate\end{tabular}}} & 10 epoch & \multicolumn{1}{c|}{375.93$\pm$2.76} & \multicolumn{1}{c|}{349.88$\pm$1.79}   & 344.72$\pm$0.82 & \multicolumn{1}{c|}{358.82$\pm$2.98}   & \multicolumn{1}{c|}{342.98$\pm$1.34}   & 337.02$\pm$0.65   & \multicolumn{1}{c|}{374.42$\pm$2.24} & \multicolumn{1}{c|}{349.98$\pm$1.23}   & 341.24$\pm$0.47 \\ \cline{2-11} 
\multicolumn{1}{|c|}{}                                                                          & 15 epoch & \multicolumn{1}{c|}{364.79$\pm$2.87} & \multicolumn{1}{c|}{345.12$\pm$1.45 } & 343.30$\pm$0.46 
& \multicolumn{1}{c|}{352.21$\pm$2.49 } & \multicolumn{1}{c|}{337.89$\pm$1.32 } & 336.48$\pm$0.62 
& \multicolumn{1}{c|}{369.90$\pm$2.50} & \multicolumn{1}{c|}{344.33$\pm$1.07 } & 342.27$\pm$0.76 \\ \hline
\multicolumn{1}{|c|}{\multirow{2}{*}{\begin{tabular}[c]{@{}c@{}}Full \\ connect\end{tabular}}}  & 10 epoch & \multicolumn{1}{c|}{359.45$\pm$3.02} & \multicolumn{1}{c|}{343.47$\pm$1.34}   & 343.89$\pm$0.45 
& \multicolumn{1}{c|}{357.72$\pm$2.45}   & \multicolumn{1}{c|}{347.63$\pm$1.90}   & 342.29 $\pm$0.42 
& \multicolumn{1}{c|}{368.76$\pm$2.79} & \multicolumn{1}{c|}{352.29$\pm$2.07}   & 346.09$\pm$0.29 \\ \cline{2-11} 
\multicolumn{1}{|c|}{}                                                                          & 15 epoch & \multicolumn{1}{c|}{357.79$\pm$2.61} & \multicolumn{1}{c|}{343.78$\pm$1.20 } & 342.29$\pm$0.41 
& \multicolumn{1}{c|}{354.73$\pm$3.00 } & \multicolumn{1}{c|}{345.05$\pm$1.98 } & 343.05$\pm$0.39  
& \multicolumn{1}{c|}{359.90$\pm$2.78} & \multicolumn{1}{c|}{348.42$\pm$1.22 } & 346.27$\pm$0.66 \\ \hline
\multicolumn{1}{|c|}{\multirow{2}{*}{\begin{tabular}[c]{@{}c@{}}Partial \\ share\end{tabular}}} & 10 epoch & \multicolumn{3}{c|}{\multirow{2}{*}{-}}                              
& \multicolumn{1}{c|}{329.97$\pm$2.77}   & \multicolumn{1}{c|}{320.09$\pm$1.67}   & 312.98$\pm$0.52   & \multicolumn{1}{c|}{334.98$\pm$2.28} & \multicolumn{1}{c|}{329.98$\pm$1.45}   & 311.09$\pm$0.22 \\ \cline{2-2} \cline{6-11} 
\multicolumn{1}{|c|}{}                                                                          & 15 epoch & \multicolumn{3}{c|}{}                                                
& \multicolumn{1}{c|}{321.45$\pm$1.90 } & \multicolumn{1}{c|}{313.35$\pm$1.54 } & 309.69$\pm$0.45 
& \multicolumn{1}{c|}{327.71$\pm$2.35} & \multicolumn{1}{c|}{321.09$\pm$1.55 } & 312.29$\pm$0.24 \\ \hline
\end{tabular}}
\caption{Policy network comparison in 10-pool 8-epoch }
    \label{tab: policy network comparison 1 in 10pool8epoch}
\end{table}

\begin{table}[htp]
\centering
 \scalebox{0.8}{
\begin{tabular}{|c|c|c|c|c|c|}
\hline
\begin{tabular}[c]{@{}c@{}}(10k data\\ 15 epoch)\end{tabular} &            & 0 depth                   & 1 depth  & 2 depth  & 3 depth \\ \hline
\multirow{3}{*}{fully-separate} 
& 17 neurons & \multirow{3}{*}{345.12$\pm$1.45 } & 337.92$\pm$1.34   & 349.23$\pm$1.02   & 371.29$\pm$2.11  \\ \cline{2-2} \cline{4-6} 
& 34 neurons &     & 337.89$\pm$1.32 & 344.33$\pm$1.07  & 379.59$\pm$1.45  \\ \cline{2-2} \cline{4-6} 
 & 68 neurons &   & 341.63$\pm$1.29  & 342.12$\pm$1.47  & 377.72$\pm$1.55  \\ \hline
\multirow{3}{*}{fully-connected}
& 17 neurons & \multirow{3}{*}{343.78$\pm$1.20 } & 344.42$\pm$2.09   & 347.72$\pm$1.96   & 312.21$\pm$1.08  \\ \cline{2-2} \cline{4-6} 
& 34 neurons &    & 345.05$\pm$1.98  & 348.42$\pm$1.22  & 311.48$\pm$1.35  \\ \cline{2-2} \cline{4-6} 
& 68 neurons &    & 344.15$\pm$1.20  & 345.23$\pm$1.07  & 308.23$\pm$2.42  \\ \hline
\multirow{3}{*}{partially-shared}                                
& 17 neurons & \multirow{3}{*}{-}        & 313.23$\pm$1.76   & 320.09$\pm$1.55   & 313.32$\pm$1.27  \\ \cline{2-2} \cline{4-6} 
& 34 neurons &                           & 313.35$\pm$1.54  & 321.09$\pm$1.55  & 315.78$\pm$1.20  \\ \cline{2-2} \cline{4-6} 
& 68 neurons &                           & 315.41$\pm$1.98  & 318.82$\pm$1.45   & 312.28$\pm$1.22  \\ \hline
\end{tabular}}
\caption{network complexity tuning in 10-pool 8-epoch }
    \label{tab: policy network comparison 2 in 10pool8epoch}
\end{table}

\subsection{Comparison with Additional Benchmarks}

\label{sub-app: benchmarks}

This section introduces several benchmark policies for comparative evaluation. Section~\ref{subsubapp: threshold} introduces a heuristic threshold-based policy. Section~\ref{subsubsec: atomic_adp},~\ref{subsubapp: NN critic} and~\ref{sub-app: order} treat several variants tested in the ablation study as the benchmark policies. The results reported in the subsequent sections are from the baseline five-pool setting.

\subsubsection{Threshold-Based Policy}
\label{subsubapp: threshold}

We introduce a heuristic threshold-based policy triggers overflow when the queue lengths exceeds a threshold. To simplify, we divide the eight decision epochs into two periods (daytime and nighttime) following the empirical policy structure. For each class and period, we search for optimal threshold values from the range 0–10, as this range covers the most frequently observed queue lengths in simulation. The resulting threshold policies yield average costs ranging from 232.87 to 356.29, with the best policy using thresholds $[0,0,2,0,0]$ at night and $[1,4,6,4,4]$ during the day. This achieves a 2.87\% improvement over the empirical policy’s performance (239.76). However, given the high computational cost of threshold tuning and the relatively modest gains, we use the empirical policy as the primary benchmark in the main paper.

\subsubsection{Atomic + ADP}
\label{subsubsec: atomic_adp}

To evaluate whether the atomic decomposition can be paired with ADP instead of PPO, we implement an Atomic+ADP benchmark.  
As shown in Figure~\ref{fig:adp_atomic}, this method performed poorly and exhibited significant ``chattering'' (large performance oscillations). It is due to the greedy search mechanism in ADP: when the value function is inaccurately estimated, the policy can deteriorate dramatically. 

\begin{figure}[ht]
\centering
\includegraphics[width=0.6\textwidth]{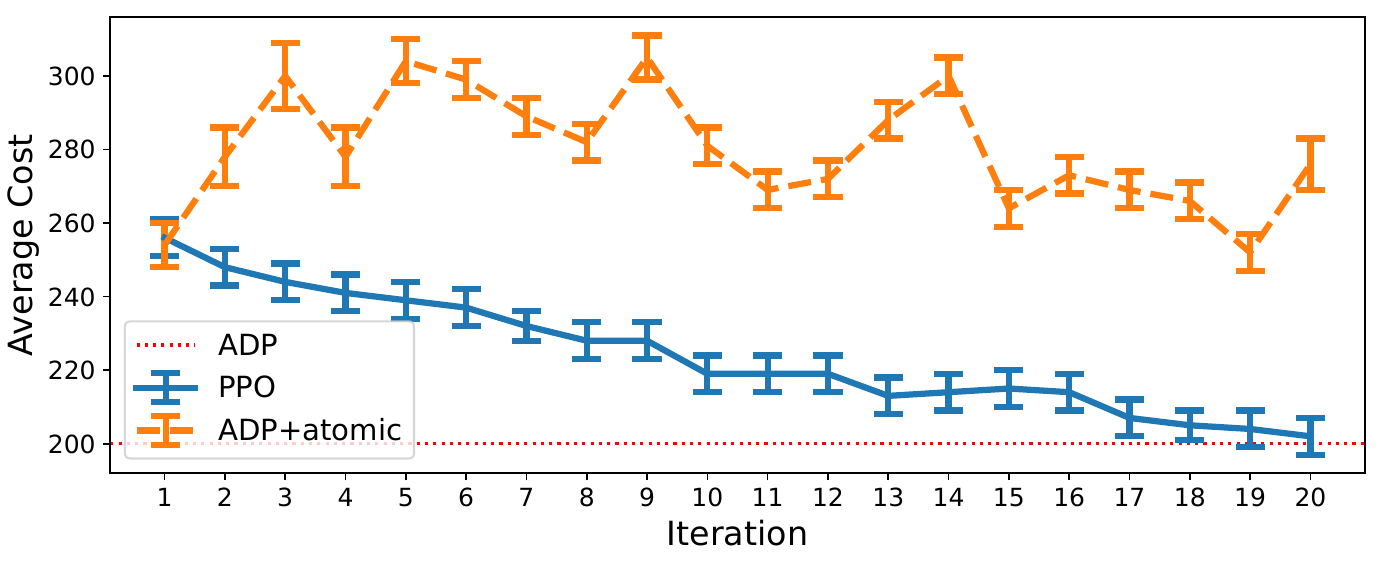}
\caption{Performance comparison of ADP, Atomic+ADP, and PPO (full algorithm) in the five-pool setting. Atomic+ADP exhibits poor performance and significant chattering.}
\label{fig:adp_atomic}
\end{figure}

To fairly assess the performance of the Atomic+ADP approach, we evaluate two configurations of action representation and a few different basis function choices. The results reported in Figure~\ref{fig:adp_atomic} correspond to the best-performing setup. Specifically, in the first configuration, we adopt the same atomic action structure as in our PPO algorithm: actions are executed in FIFO order, with each action corresponding to selecting an overflow pool $[1,2,\ldots,J]$ for the current patient. The value function approximation uses the basis set:
$$
\{X_j, Y_j, W_j, X^2_j, Y^2_j, W^2_j, V_s, e_{\text{epoch}}, e_{\text{class}}\},
$$
which augments the basis functions from~\cite{dai2019inpatient} by including linear and quadratic terms for waiting patients ($W_j$), one-hot encodings for the epoch ($e_{\text{epoch}} \in \{0,1\}^m$), and patient class ($e_{\text{class}} \in \{0,1\}^J$). This setup results in poor performance (best average cost: 282.77) and severe chattering, likely due to instability introduced by the high-dimensional one-hot features.
In the second configuration, inspired by~\cite{feng2021scalable}, actions are defined as ordered pairs $(i,j)$, where a class~$i$ patient is assigned to ward~$j$. This removes the need to order patients in advance. We test two basis sets: (i) the original from~\cite{dai2019inpatient}, ${X_j, Y_j, X^2_j, Y^2_j, V_s}$, which achieves an average cost of 268.81; and (ii) an augmented version that includes ${W_j, W^2_j}$, which further lowers the average cost (259.24) but still suffers from great instability. This latter one corresponds to the curve shown in Figure~\ref{fig:adp_atomic}.

We conjecture that designing effective basis functions is particularly difficult for the Atomic+ADP approach, as the value function must capture the impact of individual action ordering (which is different from the system-level actions considered in \cite{dai2019inpatient} and adds an inherently complex dependency). In contrast, PPO uses proximal, conservative updates that reduce sensitivity to approximation errors and eliminate chattering, leading to significantly more stable performance in our long-run average cost setting.

\subsubsection{NN-based Value Network}
\label{subsubapp: NN critic}

Following the standard setup of actor-critic RL algorithms, we have conducted experiments using an NN-based critic (value network) in place of LSTD. All experiments use 10k simulation days per iteration and we test both five-pool and ten-pool settings under sequential and batched atomic action updates. Across the board, the NN critic yields higher average costs and longer computation times, with performance gaps more pronounced in the ten-pool setting. For example, in the ten-pool case without batching, the NN critic results in an average cost of 324.23 (vs. 304.52 with LSTD), requiring 31\% more computation time. Using batching, the NN critic yields an average cost of 213 (vs. 207 for LSTD) in the five-pool setting with a 14\% increase in runtime; and 328.67 (vs. 309.02) in the ten-pool setting, with a 19\% increase in runtime.

The performance degradation under the NN critic likely stems from higher variance and less stable advantage estimates when simulation data are limited. In contrast, LSTD with tailored basis functions provides lower-variance, more accurate estimates that support more stable learning. Moreover, NN critics require extra time for training, while LSTD adds minimal overhead. 

\subsubsection{Action Ordering}
\label{sub-app: order} 

When performing sequential updates under the atomic action setup, one natural question is whether the order in which patients are processed affects performance. To investigate this, we have tested five ordering rules: Random, First-In–First-Out (FIFO), and three priority-based rules: (i) Class 1 to Class 5, (ii) Class 5 to Class 1, and (iii) Class 1,4,5,2,3 (from most- to least-loaded class). We compare performance under these different ordering rules with those from the batching approximation in two five-pool settings: the baseline setting and a more congested variant (with 3 fewer beds in ward 1 and 2 fewer in ward 2); the latter leads to longer queues and more batching. 

\begin{table}[htp]
    \centering
    \scalebox{0.8}{\begin{tabular}{|c|c|c|c|c|c|c|}
    \hline
        Ordering & Random & FIFO& Priority Rule (i) & Priority Rule (ii)& Priority Rule  (iii)&Batching \\\hline
        Baseline & 207.19$\pm$4.27 & 207.22$\pm$3.82& 203.24$\pm$4.51 &205.33$\pm$4.13 &203.11$\pm$4.00& 208.82$\pm$4.78\\\hline
        Crowded & 248.27$\pm$4.22 & 245.19$\pm$4.02& 236.38$\pm$4.11 & 244.89$\pm$4.43 &232.14$\pm$3.97& 247.33$\pm$4.09\\\hline
    \end{tabular}}
    \caption{Comparison of Different Ordering. }
    \label{tab: order comparison}
\end{table}

Results in Table~\ref{tab: order comparison} show minimal performance differences (about 2\%) across ordering rules in the baseline case. In the congested setting, differences became more pronounced (around 8–9\%). While batching performs comparably to Random and FIFO, it is generally worse than the best priority-based rule (Priority Rule iii). In summary, ordering does influence performance, but the effect is minimal when queue lengths are low (3–5 patients per class on average). When the system is more congested, the impact grows, but this could be mitigated by increasing decision frequency to reduce batch size. These findings support our use of the batching approximation as proposed in Section~5, which offers substantial computational savings (reducing per-iteration runtime by roughly 50\%) with only a minor trade-off in performance.

\subsection{Additional Performance Metrics}
\label{sub-app:additional_metrics}

Figure~\ref{fig:performance_metrics} plots the trade-off between system congestion and overflow rates for the five- and ten-pool systems. We fix the holding cost at $C=6$ and vary the overflow costs to construct the blue curve. Plots (a) and (c) show the daily overflow rate versus peak hourly queue length; plots (b) and (d) show the daily overflow rate versus the proportion of patients waiting longer than 4 hours, referred to as the 4-hour non-service level. Across all settings, the PPO policy achieves a Pareto improvement: it either reduces waiting at a similar overflow rate, lowers the overflow rate without increasing congestion, or improves both dimensions simultaneously. The improvement is more pronounced in the larger ten-pool system. For example, in Figures~\ref{fig:10pool_ql_of} and~\ref{fig:10pool_sl_of}, relative to the empirical policy, PPO achieves about 4 fewer overflow patients per day (approximately 16 fewer patients in non-primary wards, assuming an average LOS of 4 days) at the same peak queue length, or about 4 fewer patients waiting at the peak hour and a 15\% absolute reduction in the peak 4-hour non-service level at the same overflow rate.

\begin{figure}[htp]
\centering
\subfigure[Five-pool]{\includegraphics[width=0.35\linewidth]{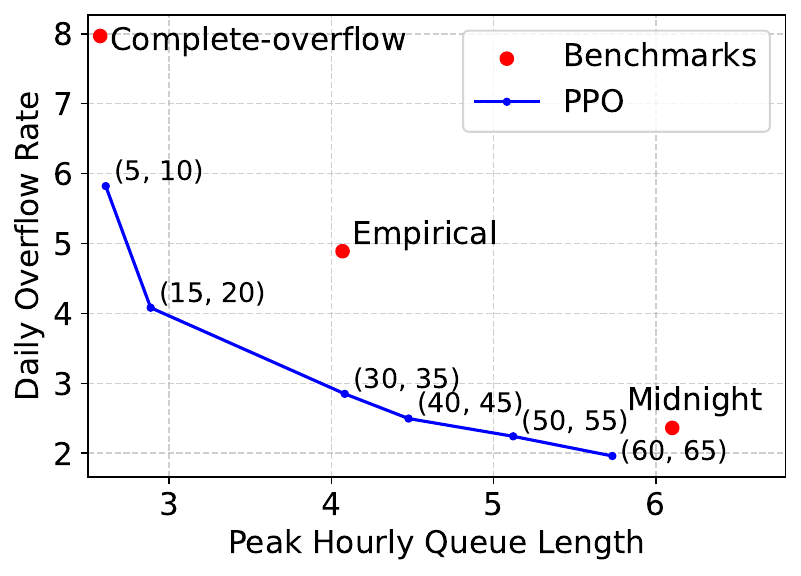}
\label{fig:5pool_ql_of}
}
\subfigure[Five-pool]{\includegraphics[width=0.35\linewidth]{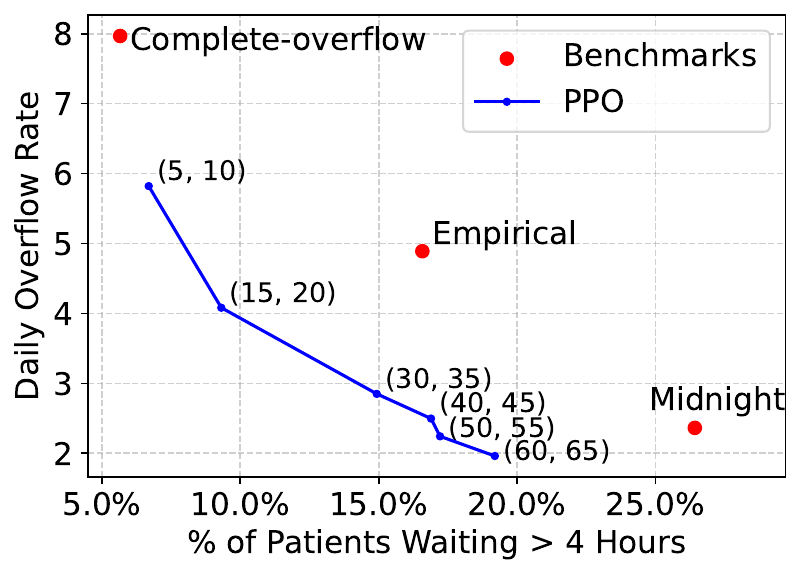}
\label{fig:5pool_sl_of}
}
\subfigure[Ten-pool]{\includegraphics[width=0.35\linewidth]{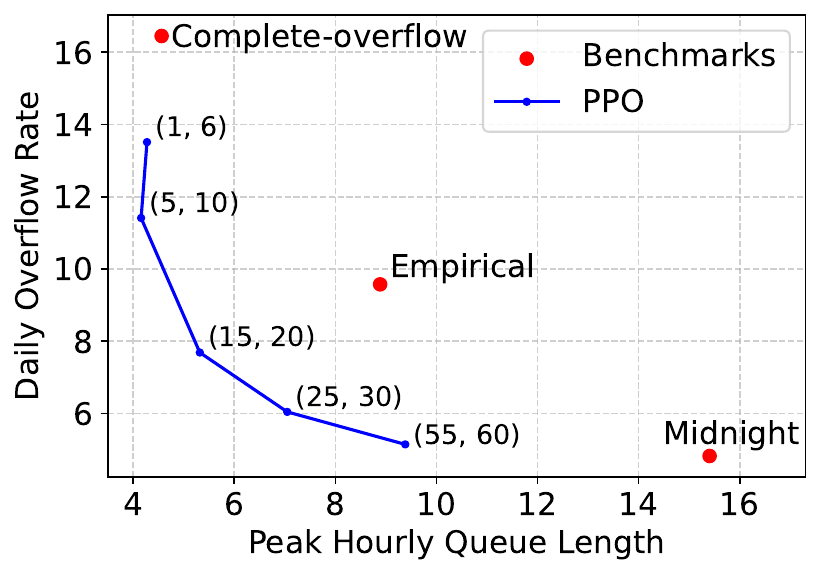}
\label{fig:10pool_ql_of}
}
\subfigure[Ten-pool]{\includegraphics[width=0.35\linewidth]{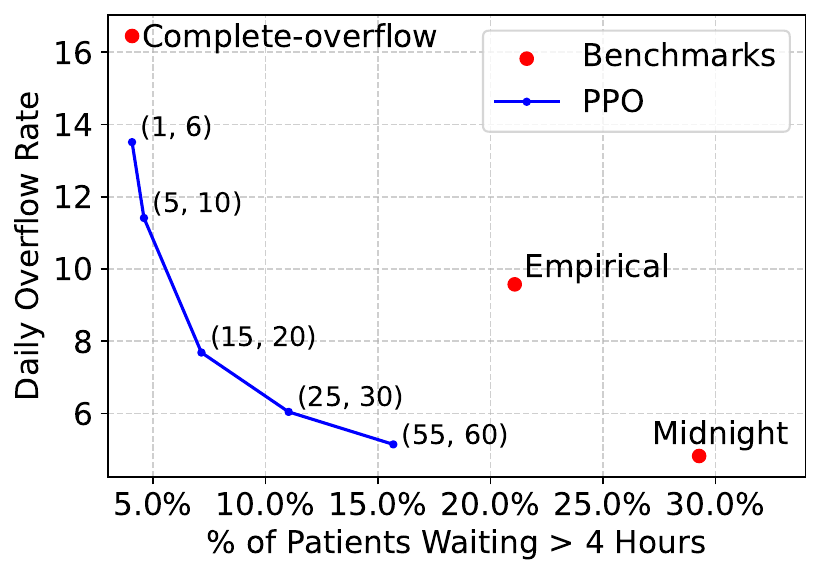}
\label{fig:10pool_sl_of}
}
\caption{Additional Performance Metrics. The numbers in parentheses along the blue line indicate the corresponding cost parameter pairs $(B_1,B_2)$. For the ten-pool settings, the additional cross-department overflow cost parameters are set as $(B_3, B_4)=(B_1+10,B_2+10)$.}
\label{fig:performance_metrics}
\end{figure}

\vspace{-0.2in}
Across the experiments, we fix the overflow cost gap $B_2-B_1$ at 5. We also conduct additional experiments with $C=6$ and $B_1=15$, while varying the gap $B_2-B_1$. As $B_2$ increases, we find that the PPO policy (i) uses fewer secondary overflow assignments and (ii) shifts toward more preferred overflow assignments, while allowing more patients to wait. Thus, increasing the overflow cost gap changes not only the overall use of overflow, but also the composition of overflow assignments.

\section{Policy Visualization and SHAP Analysis Details}
\label{app:shap}
This appendix provides a detailed analysis of the operational behavior of the learned PPO policy. Using visualization and explainable AI techniques, we illustrate its decision-making process and derive actionable managerial insights. The PPO policy examined here is trained under the baseline five-pool setting. Following \cite{dai2019inpatient}, the baseline five-pool setting assigns each pool a capacity of 
$(N_1,\dots,N_5) = (60, 64, 67, 62, 62)$. 
The preferred overflow wards for classes $(1,2,3,4,5)$ are $(5,3,2,2,1)$, respectively; the secondary choices are $\{(2,3), (4,5), (1,5), (1,3), (2,3)\}$, respectively.

Section~\ref{subapp:direct visualization} presents direct policy visualizations that reveal how different system factors influence the policy’s decisions. Section~\ref{subapp:shap} further employs SHAP (SHapley Additive exPlanations) analysis to generate additional insights into the learned policy, in particular, uncover underlying network effects. Section~\ref{sub-app: visual compare with adp} compares PPO policy and ADP policy by direct visualization.

\subsection{Direct Policy Visualization}
\label{subapp:direct visualization}
We visualize the overflow probabilities of class~1 patients to illustrate how key state variables influence overflow decisions. In each visualization, two variables are varied (shown on the axes), while all other components are fixed at representative baseline values: for fixed $X$, we set $X_j=N_j$ if $j\neq 1$ and $X_1=N_1+5$.

\subsubsection{Impact of patient count. }
Figure~\ref{fig:policy_viz} shows how the patient count ($X_j$) influence the policy.
The first three subfigures display how the probability that a class~1 patient is not overflowed ($P_{11}$)  changes with $X_1$ and $X_2$ at midnight, 12~pm, and 9~pm (Epoch $h=$ 0, 4, 7); the fourth subfigure shows how the overflow probability to ward~5 ($P_{15}$) changes with $X_5$ and $X_2$ at 9~pm ($h=7$).

Three key patterns emerge: (i) the policy prescribes a higher probability of overflow (i.e., a lower $P_{11}$) when $X_1$ is high, indicating congestion in the source ward, or when $X_2$ is low, suggesting greater available capacity in the overflow ward—this aligns with operational intuition; (ii) the policy is more aggressive in overflowing patients during midnight and 9~pm ($h=$0, 7), but more conservative during daytime hours such as 12~pm ($h=4$), reflecting anticipation of discharges that typically occur during the day and free up capacity; and (iii) the right-most plot demonstrates that the policy dynamically balances medical closeness between the primary
and overflow wards against system-level congestion, that is, increases the chance of assigning patients to secondary overflow wards (ward~2) when the preferred one (ward~5) is congested. These observations collectively indicate that the PPO policy captures realistic and adaptive operational behavior.

\vspace{-0.24in}
\begin{figure}[ht]
\centering
\includegraphics[width=0.8\textwidth]{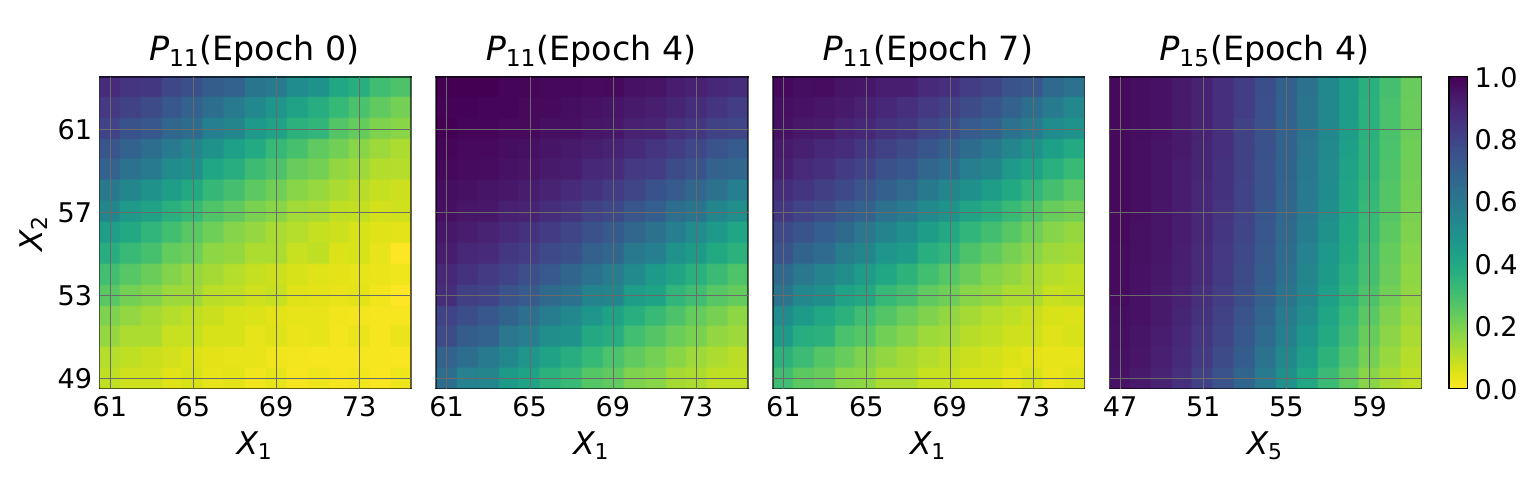}
\caption{{Policy visualization of class 1 patients on varying $X$. The vertical axis shows the value of $X_2$; the horizontal axis shows the value of $X_1$ in the left three plots and of $X_5$ in the right-most one. 
}}
\label{fig:policy_viz}
\end{figure}

\subsubsection{Impact of to-depart quantity. }
Figure~\ref{fig:policy_viz_y} shows that the number of patients to be discharged by the end of the day ($Y_j$) is also an important factor. The left figure illustrates how $P_{11}$ changes with $Y_1$ and $Y_2$ at 12~pm (Epoch $h=4$), while the right figure shows how $P_{15}$ changes with $Y_2$ and $Y_5$ at the same time. For fixed $Y$ (the ones we do not vary), we set $Y_j=0$.

The left figure shows that the policy prescribes a lower probability of overflow (i.e., a higher $P_{11}$) when $Y_1$ is high, indicating that more capacity is anticipated to free up soon; or when $Y_2$ is low, implying that more beds in the overflow ward will remain occupied. The right figure reveals that 
the policy increases the change of overflowing patients to one of the secondary wards (ward~2) when the preferred overflow ward (ward~5) has fewer patients to be discharged.

\begin{figure}[htp]
 \centering
\includegraphics[width=0.42\textwidth]{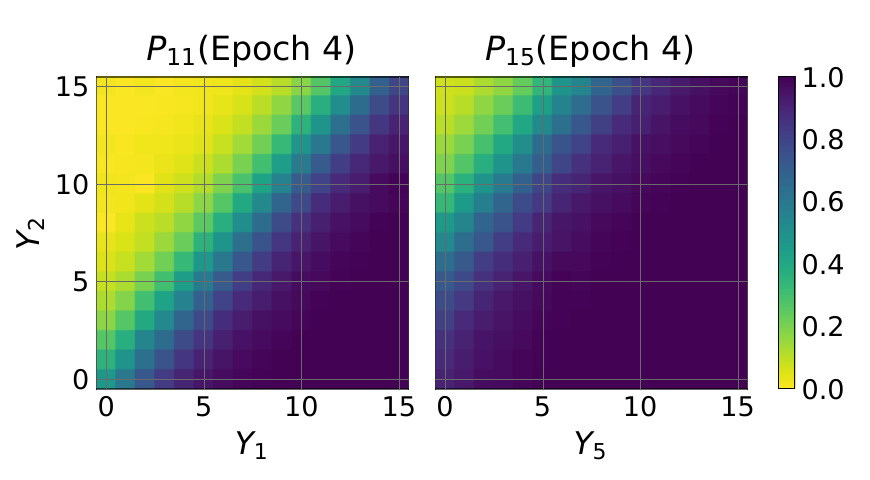}
    \caption{{Policy visualization of class 1 patients on varying $Y$. The vertical axis shows the value of $Y_2$; the horizontal axis shows the value of $Y_1$ in the left plot and of $Y_5$ in the right one. 
    }
    }
    \label{fig:policy_viz_y}
    \end{figure}

\subsection{SHAP Analysis for Nuanced Network Effects}
\label{subapp:shap}
To uncover more subtle and non-intuitive patterns, we further employ SHAP (SHapley Additive exPlanations) analysis to interpret the PPO policy. SHAP quantifies how each input feature of the policy network affects the predicted overflow probabilities for individual samples, with the magnitude of a SHAP value indicating a feature’s importance and its sign showing whether the feature increases (positive) or decreases (negative) the output.   The output of interest is the probability of \emph{not} overflowing a patient (e.g., $P_{ii}$ for class~$i$), while the input features include $(Q,Z,Y)$. Recall that $Q_i = 0 \vee (X_i - N_i)$ represents the number of waiting patients, and $Z_i = N_i \wedge X_i$ represents the number of in-service patients. We use $(Q,Z)$ instead of $X$ directly to disentangle the effects of waiting and in-service patients, as they play distinct operational roles in overflow decisions.

Figure~\ref{fig:shap_analysis} shows the top five most influential features for the decision not to overflow a class~2 patient ($P_{22}$) at 12~pm ($h=4$), and Figure~\ref{fig:shap_sharing} presents the corresponding results for class~4 ($P_{44}$). Source-ward features are excluded from both plots. In both cases, the most salient features prominently include those of the  relevant overflow wards (e.g., $Y_3$, $Z_5$ and $Y_5$ for $P_{22}$; $Y_3$, $Z_2$, $Y_1$ and $Y_2$ for $P_{44}$), confirming that the PPO policy captures realistic operational logic consistent with the direct visualization analysis in Section~\ref{subapp:direct visualization}. Beyond these intuitive dependencies, SHAP also highlights additional network-level effects: (i) \textbf{Anticipating future congestion:} Figure~\ref{fig:shap_anticipate} (for class~2) shows that the to-depart count from ward~1 ($Y_1$, positive correlation) and the occupancy of ward~1 ($Z_1$, negative correlation) are among the most influential features. This indicates a \emph{proactive preparation} effect. That is, when ward~1 is more congested (higher $Z_1$) and/or has fewer discharges (lower $Y_1$), ward~2 overflows its own patients more aggressively (lower $P_{22}$) to reserve capacity for future class~1 arrivals, since ward~2 is frequently used as a overflow ward for class~1. (ii) \textbf{Sharing overflow pools:} For class~4 ($P_{44}$), congestion ($Z_5$) in ward~5 strongly influences decisions, despite wards~4 and~5 not directly overflowing to each other. This reflects an indirect ``chain effect'' arising from the fact that both classes share wards~1, 2, and 3 as overflow choices. That is, when ward~5 becomes more congested (higher $Z_5$), it is more likely to send patients to these shared overflow wards. Anticipating this competition, the policy instructs ward~4 to overflow its own patients more aggressively (lower $P_{44}$).

\begin{figure}[htp]
\centering 
    \subfigure[$P_{22}$]{  
    \includegraphics[width=0.4\linewidth]{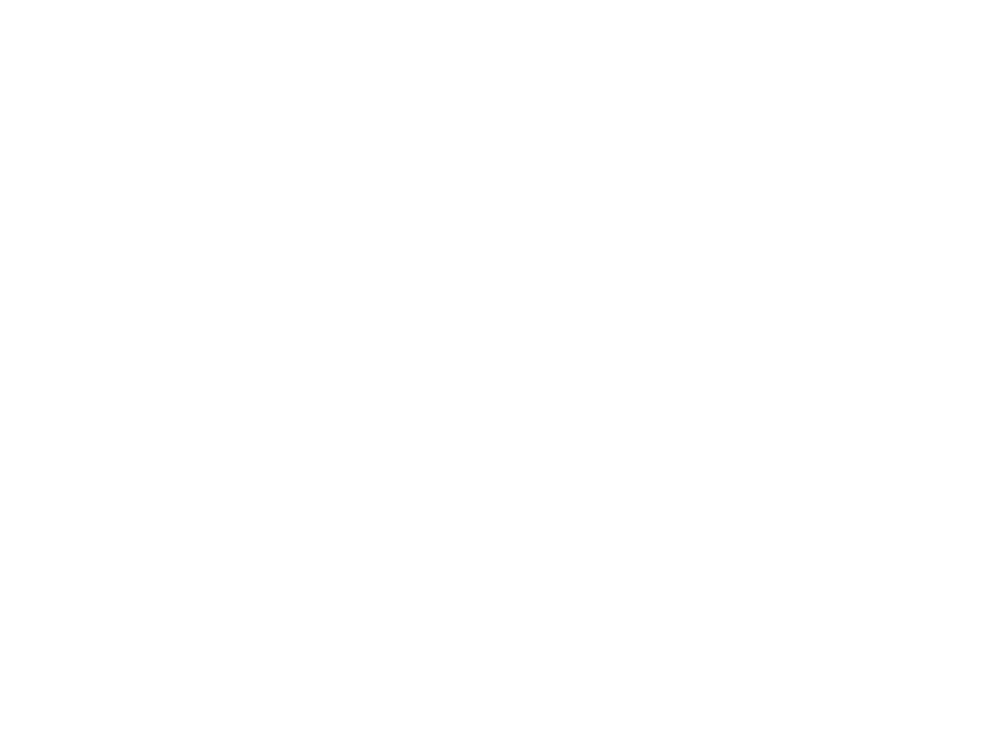}
     \label{fig:shap_anticipate}
    }
    \subfigure[$P_{44}$]{
    \includegraphics[width=0.4\linewidth]{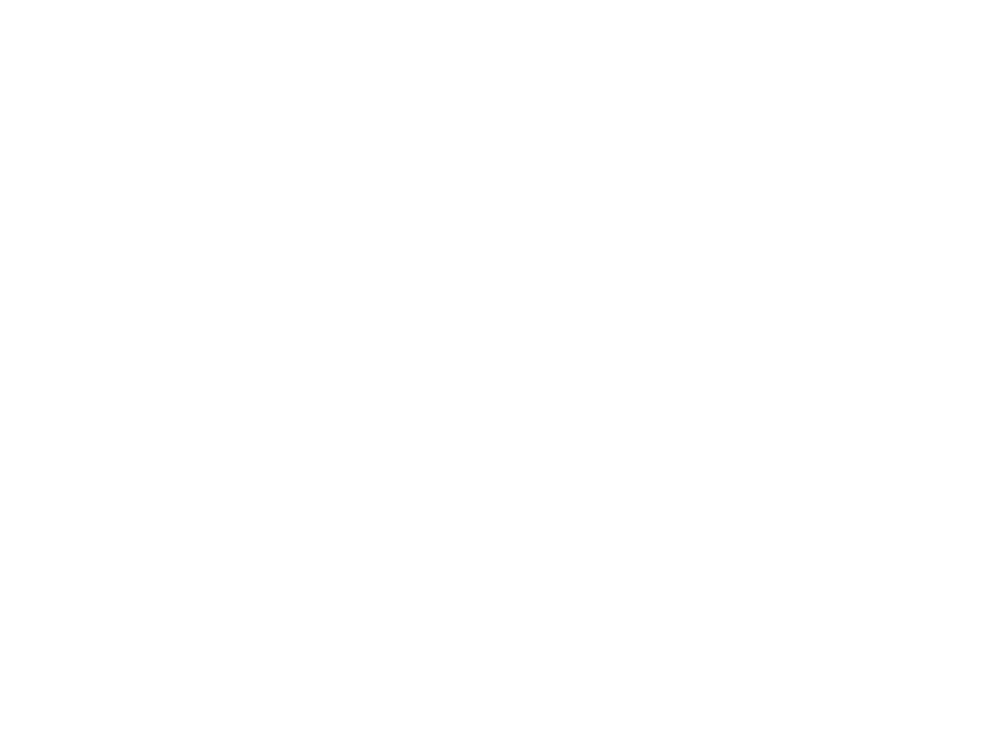} 
    \label{fig:shap_sharing}
    }
    \caption{SHAP analysis of probabilities of staying for epoch $h=4$. Left for class 2 patients staying (not overflow) in ward 2; right for class 4 patients staying (not overflow) in ward 4. Each figure shows the top five most important features, excluding the source-ward features.  
    } 
    \label{fig:shap_analysis}
    \end{figure}

\subsection{Comparison Between PPO Policy and ADP Policy}
\label{sub-app: visual compare with adp}
To directly compare what PPO and ADP policies prescribed, we visualize both under identical system states. For the deterministic ADP policy, we plot overflow proportions (Figure~\ref{fig:pol_vis_adp}) to align with the randomized PPO policy representation (Figure~\ref{fig:policy_viz}). The resulting patterns are broadly consistent across methods. Both display time-varying and state-dependent behaviors, such as more aggressive overflow when source wards are congested or overflow wards are less occupied.
\begin{figure}[htp]
    \centering
    \includegraphics[width=0.8\linewidth]{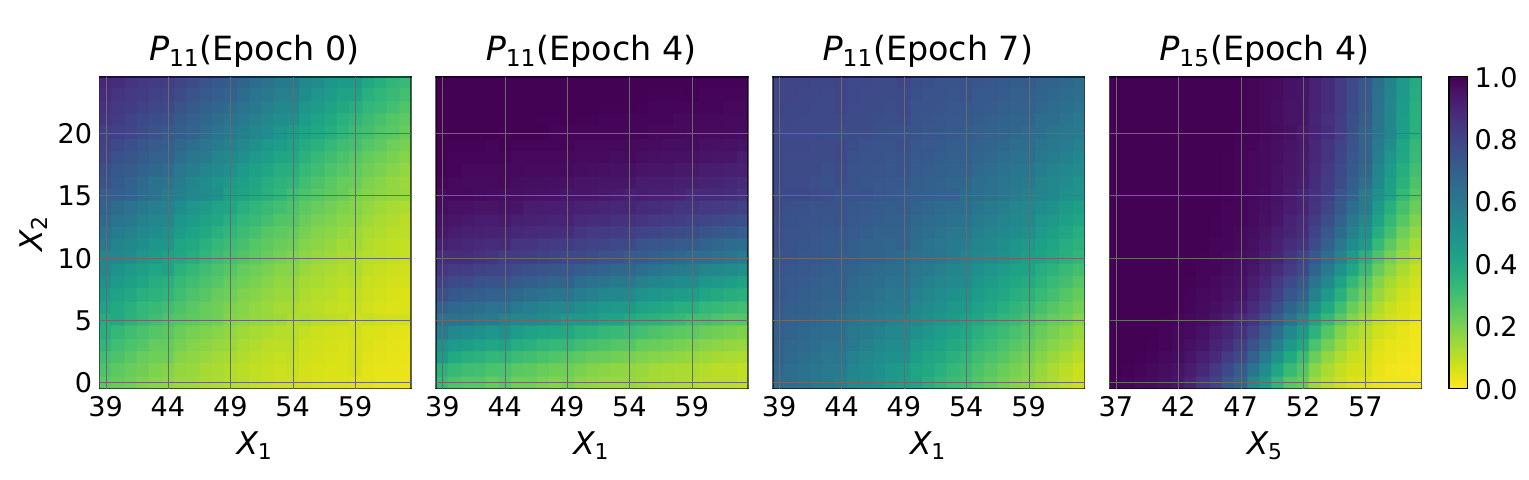}
    \caption{Policy visualization when varying $X$ of the ADP policy. The vertical axis shows the value of $X_2$; the horizontal axis shows the value of $X_1$ in the left three plots and of $X_5$ in the right-most one.}
    \label{fig:pol_vis_adp}
\end{figure}

\clearpage

\begin{center}\bfseries\LARGE Technical Companion for 
\\
``Inpatient Overflow Management
with Proximal Policy Optimization''
\end{center}

\setcounter{page}{1}
\setcounter{section}{0}

\medskip 

\section{Proof of Proposition~\ref{prop: periodic average cost upper bound}}
\label{app:proof}

Following~\cite{dai2022queueing}, the average cost gap between the two policies can be written as 
\begin{align*}
    &(\mu^{h}_{\theta})^T \tilde{\textbf{g}}_{\theta}^{h}-(\mu^{h}_\eta)^T \tilde{\textbf{g}}_\eta^{h}\\
    =&(\mu^{h}_{\theta})^T(\tilde{\textbf{g}}_{\theta}^{h}+(\tilde{\textbf{P}}_{\theta}^{h}-I)\mathbf{v}_\eta^{h}))-(\mu^{h}_\eta)^T \tilde{\textbf{g}}_\eta^{h}\\
    =&(\mu_\eta^{h} )^T(\tilde{\textbf{g}}_{\theta}^{h}-(\mu_\eta^{h})^T\tilde{\textbf{g}}_{\eta}^{h}\textbf{e}+(\tilde{\textbf{P}}_{\theta}^{h}-I)\mathbf{v}_\eta^{h})+(\mu_{\theta}^{h}-\mu_\eta^{h})^T(\tilde{\textbf{g}}_{\theta}^{h}+(\tilde{\textbf{P}}_{\theta}^{h}-I)\mathbf{v}_\eta^{h}) 
    \\
    =&N_1^h(\theta,\eta)+N_2^h(\theta,\eta),
\end{align*}
where the first equation holds since the stationary distribution $\mu_{\theta}^{h}$ satisfying $(\mu_{\theta}^{h})^T(\tilde{\textbf{P}}_{\theta}^{h}-I)=0$. The third equation holds since for any constant $\alpha$, $(\mu_{\theta}^{h}-\mu_\eta^{h})^T\alpha\textbf{e}=\alpha(\mu_{\theta}^{h})^T \textbf{e}-\alpha(\mu_\eta^{h})^T \textbf{e}=0$, and we choose $\alpha=(\mu_\eta^{h})^T\tilde{\textbf{g}}_{\eta}^{h}$. Next, we characterize the decay rate of $N_1^h$ and $N_2^h$.

For a given vector $\omega\in \mathcal{S}^{h}$, define two types of $\mathcal{V}$-norm as
\begin{equation} \|\omega\|_{\infty,\mathcal{V}}=\sum_{s\in\mathcal{S}^{h}} \frac{|\omega(s)|}{\mathcal{V}(s)},\quad \|\omega\|_{1,\mathcal{V}}=\sum_{s\in\mathcal{S}^h} |\omega(s)|\mathcal{V}(s).
    \label{equ: V-norm def}
\end{equation}
Note that slightly different from the definition given in~\cite{dai2022queueing}, we focus on each given $h$ and the summation is taken over the states in the subspace $\mathcal{S}^{h}$. Similarly, we slightly adapt the $\mathcal{V}$-norm for a given matrix $\Omega \in \mathcal{S}^{h}\times\mathcal{S}^{h}$ as 
\begin{equation}
    \|\Omega\|_{\mathcal{V}}=\sup_{ s\in\mathcal{S}^h} \frac{1}{\mathcal{V}(s)} \sum_{ s'\in\mathcal{S}^{h}}|\Omega(s,  s')|\mathcal{V}( s' ) . 
    \label{equ: V-norm matrix def}
\end{equation} 
We further denote 
$$\tilde N^h(\theta,\eta)=\tilde{\textbf{g}}_{\theta}^{h}-(\mu_\eta^{h})^T\tilde{\textbf{g}}_{\eta}^{h}\textbf{e}+(\tilde{\textbf{P}}_{\theta}^{h}-I)\mathbf{v}_\eta^{h}, 
$$ 
which is an $\mathcal{S}^h$-dimensional vector.
Following~\cite{dai2022queueing}, we can bound the absolute value of the two scalars $N_1^h$ and $N_2^h$ as 
\begin{equation}
    \begin{aligned}
        |N_1^h(\theta,\eta)|
    \le & (\mu_\eta^{h})^T\mathcal{V} \cdot \|\tilde N^h(\theta,\eta)\|_{\infty,\mathcal{V}}, \\
   |N_2^h(\theta,\eta)|
   \le &\|\mu_{\theta}^{h}-\mu_\eta^{h}\|_{1,\mathcal{V}} \cdot\| \tilde N^h(\theta,\eta)\|_{\infty,\mathcal{V}}. 
    \end{aligned}
    \label{equ: bounds for N1 N2}
\end{equation}

For the bound of $N_1^h$, only the term $\|\tilde N^h(\theta,\eta)\|_{\infty,\mathcal{V}}$ relates to the new parameter $\theta$, while for $N_2^h$, it contains an additional term $\|\mu_{\theta}^{h}-\mu_\eta^{h}\|_{1,\mathcal{V}}$ which also relates to $\theta$. 
According to Theorem~1 and Lemma~5 in \cite{dai2022queueing}, we have $\|\mu_{\theta}^{h}-\mu_\eta^{h}\|_{1,\mathcal{V}}=O(\|(\tilde{\textbf{P}}_{\theta}^{h} - \tilde{\textbf{P}}_\eta^{h}) Z_\eta\|_{\mathcal{V}})$ and $\|(\tilde{\textbf{P}}_{\theta}^{h} - \tilde{\textbf{P}}_\eta^{h}) Z_\eta\|_{\mathcal{V}} \rightarrow 0$ as $\theta\rightarrow \eta$, where the matrix $Z_\eta$ is 
         defined as
             \begin{equation}
       Z_\eta:=\sum_{n=0}^{\infty}(\tilde{\textbf{P}}_\eta^{h}-\Pi_\eta^{h})^n.   
           \label{equ: fondamental matrix}
       \end{equation}
To get the exact order of $N_1^h, N_2^h$ with respect to $\|r_{\theta,\eta}^h-1\|$, we need to further check the order of $\|\tilde N^h(\theta,\eta)\|_{\infty,\mathcal{V}}$ and $\|(\tilde{\textbf{P}}_{\theta}^{h} - \tilde{\textbf{P}}_\eta^{h}) Z_\eta\|_{\mathcal{V}} $.

For $\|\tilde N^h(\theta,\eta)\|_{\infty,\mathcal{V}}$, we can rewrite and bound it as follows:
\begin{equation}
    \begin{aligned}
 \|\tilde N^h(\theta,\eta)\|_{\infty,\mathcal{V}} 
    \le &  \|(\tilde{\textbf{P}}_{\theta}^{h} - \tilde{\textbf{P}}_\eta^{h}) \mathbf{v}_\eta^{h}||_{\infty,\mathcal{V}} + ||\tilde{\textbf{g}}_{\theta}^{h} - \tilde{\textbf{g}}_{\eta}^{h}\|_{\infty,\mathcal{V}}.
\end{aligned}
\label{Equ: UB tilde N}
\end{equation}
Our remaining task is to analyze the order of $\|(\tilde{\textbf{P}}_{\theta}^{h} - \tilde{\textbf{P}}_\eta^{h}) \mathbf{v}_\eta^{h}\|_{\infty,\mathcal{V}}$ and $\|\tilde{\textbf{g}}_{\theta}^{h} - \tilde{\textbf{g}}_{\eta}^{h}\|_{\infty,\mathcal{V}}$. Note that the second term on the right-hand side does not exist in \cite{dai2022queueing} since they have action-independent cost. 
Also note that they characterize the decay rates of $N^h(\theta,\eta)$ in terms of $D_{\theta,\eta}^h:=\|(\tilde{\textbf{P}}_{\theta}^{h} - \tilde{\textbf{P}}_\eta^{h}) Z_\eta\|_{\mathcal{V}}$. However, since we cannot directly bound $\|\tilde{\textbf{g}}_{\theta}^{h} - \tilde{\textbf{g}}_{\eta}^{h}\|_{\infty,\mathcal{V}}$ with respect to $D_{\theta,\eta}^h$, we need to take one step further to bound $N_1, N_2$ directly to $\|r_{\theta,\eta}^h-1\|_\infty$ instead of $D_{\theta,\eta}^h$.

First, we consider $\|(\tilde{\textbf{P}}_{\theta}^{h} - \tilde{\textbf{P}}_\eta^{h}) \mathbf{v}_\eta^{h}\|_{\infty,\mathcal{V}}$.  According to Lemma 3 of \cite{dai2022queueing}, the relative value function $\textbf{v}_\eta^{h}$ can be rewritten as 
$$\textbf{v}_\eta^{h}=Z_\eta \big( \tilde{\textbf{g}}_{\eta}^{h}-(\mu_\eta^h)^T\tilde{\textbf{g}}_\eta^h\textbf{e} \big). 
$$
Therefore, we can bound 
\begin{equation}
    \begin{aligned}
\|(\tilde{\textbf{P}}_{\theta}^{h} - \tilde{\textbf{P}}_\eta^{h}) \mathbf{v}_\eta^{h}\|_{\infty,\mathcal{V}}  
\le & \|\tilde{\textbf{P}}_{\theta}^{h}-\tilde{\textbf{P}}_\eta^{h}\|_{\mathcal{V}}
\|Z_\eta\|_{\mathcal{V}}\cdot 
\Big(  \|\tilde{\textbf{g}}_{\eta}^{h}\|_{\infty,\mathcal{V}}+(\mu_\eta^h)^T\tilde{\textbf{g}}_\eta^h
\Big) 
. 
\end{aligned}
\label{equ: Pv bound}
\end{equation} 
Then, to bound the term $\|\tilde{\textbf{P}}_{\theta}^{h}-\tilde{\textbf{P}}_\eta^{h}\|_\mathcal{V}$, we recall that the one-day transition matrices are specified as:
\begin{equation}
    \begin{aligned}
    &\tilde{\textbf{P}}_{\eta}^h=\textbf{P}_\eta^{h,h+1}\textbf{P}_\eta^{h+1,h+2}\cdots\textbf{P}_\eta^{h-1,h}, \\
&\tilde{\textbf{P}}_{\theta}^h=\textbf{P}_\theta^{h,h+1}\textbf{P}_\eta^{h+1,h+2}\cdots\textbf{P}_\eta^{h-1,h} . 
    \end{aligned}
\end{equation}
We denote the elements of these one-day transition matrices as $\{\tilde p_\eta^h(s'|s), s,s'\in\mathcal{S}^h\}$ and $\{\tilde{p}_{\theta}^h(s'|s), s,s'\in\mathcal{S}^h\}$, respectively. The probability $\tilde p_\theta^h(s'|s)$ follows 
\begin{align*}
    \tilde p_\theta^h(s'|s)=&\sum_{s^{h+1}\in\mathcal{S}^{h+1},\cdots,s^{h-1}\in\mathcal{S}^{h-1}}p_\theta^{h,h+1}(s^{h+1}|s)p_\eta^{h+1,h+2}(s^{h+2}|s^{h+1})\cdots p_\eta^{h-1,h}(s'|s^{h-1})
    \\
    =&\sum_{f\in\mathcal{A}(s)}\pi_\theta(f\mid s)\sum_{s^{h+1}\in\mathcal{S}^{h+1},\cdots,s^{h-1}\in\mathcal{S}^{h-1}}p^{h,h+1}(s^{h+1}|s,f)p_\eta^{h+1,h+2}(s^{h+2}|s^{h+1})\cdots p_\eta^{h-1,h}(s'|s^{h-1}), 
\end{align*}
where $\{p^{\ell, \ell+1}_\eta(s'|s), s\in\mathcal{S}^\ell, s'\in\mathcal{S}^{\ell+1}\}$ is the set of elements of one-epoch transition matrix $\textbf{P}^{\ell,\ell+1}_\eta$, and $p^{h,h+1}(s^{h+1}|s,f)$ is the one-epoch transition probability from state $s\in\mathcal{S}^h$ to $s^{h+1}\in\mathcal{S}^{h+1}$ given action $f$. Here we use the fact that, after the action $f$ is fixed, the transition only depends on the arrivals and departures during this epoch and no longer depends on the action or policy. We denote $$\tilde p_\eta^h(s'|s,f)=\sum_{s^{h+1}\in\mathcal{S}^{h+1},\cdots,s^{h-1}\in\mathcal{S}^{h-1}}p^{h,h+1}(s^{h+1}|s,f)p_\eta^{h+1,h+2} (s^{h+2}|s^{h+1})\cdots p_\eta^{h-1,h}(s'|s^{h-1}),$$ which was also used in Equation~\eqref{equ: def adv one-day h} and introduced there. This term is independent of the new policy parameter $\theta$. 
Using this term, we can write that 
$\tilde p_\theta^h(s'|s) =\sum_{f\in\mathcal{A}(s)}\pi_\theta(f\mid s)\tilde p_\eta^h(s'|s,f)$.  
Similarly, we have 
$\tilde p_\eta^h(s'|s)
=\sum_{f\in\mathcal{A}(s)}\pi_\eta(f\mid s)\tilde p_\eta^h(s'|s,f)$.

Therefore, the term $\|\tilde{\textbf{P}}_{\theta}^{h}-\tilde{\textbf{P}}_\eta^{h}\|_\mathcal{V}$ can be bounded as
\begin{align*}
    \|\tilde{\textbf{P}}_{\theta}^{h}-\tilde{\textbf{P}}_\eta^{h}\|_{\mathcal{V}}=&\sup_{s\in\mathcal{S}^h}\frac{1}{\mathcal{V}(s)}\sum_{s'\in\mathcal{S}^h}|\tilde{p}^h_\theta(s'|s)-\tilde p^h_\eta(s'|s)|\mathcal{V}(s')\\
    =&\sup_{s\in\mathcal{S}^h}\frac{1}{\mathcal{V}(s)}\sum_{s'\in\mathcal{S}^h}
    \bigg| \sum_{f\in\mathcal{A}(s)}\pi_\theta(f\mid s)\tilde p_\eta^h(s'|s,f)-\sum_{f\in\mathcal{A}(s)}\pi_\eta(f\mid s)\tilde p_\eta^h(s'|s,f)
    \bigg| \mathcal{V}(s')\\
    =& \sup_{s\in\mathcal{S}^h} \frac{1}{\mathcal{V}(s)}\sum_{s'\in\mathcal{S}^h}\bigg|\sum_{f\in\mathcal{A}(s)}(r_{\theta,\eta}(f\mid s)-1)\pi_\eta(f\mid s)\tilde p_\eta^h(s'|s,f)\bigg|\mathcal{V}(s
    ')\\
    \le& \|\textbf{r}^h_{\theta,\eta}-1\|_\infty \sup_{s\in\mathcal{S}^h}\frac{1}{ \mathcal{V}(s)} \sum_{s'\in\mathcal{S}^h} \tilde p_\eta(s'|s)\mathcal{V}(s')
    \\
    <& \|\textbf{r}^h_{\theta,\eta}-1\|_\infty\sup_{s\in\mathcal{S}^h}\frac{1}{\mathcal{V}(s)}(b\mathcal{V}(s)+d\textbf{1}_{C}(s))
    \\
    \le & \|\textbf{r}^h_{\theta,\eta}-1\|_\infty\sup_{s\in\mathcal{S}^h}(b+\frac{d}{\mathcal{V}(s)})
    \\
    \le &\|\textbf{r}^h_{\theta,\eta}-1\|_\infty(b+d),
\end{align*}
where the second inequality holds because of the drift condition in Assumption~\ref{ass: drift}, and the last inequality holds since we assume that $\mathcal{V}\ge 1$. 

By plugging this upper bound into Equation~\eqref{equ: Pv bound}, we have
\begin{align*}
    \|(\tilde{\textbf{P}}_{\theta}^{h} - \tilde{\textbf{P}}_\eta^{h}) \mathbf{v}_\eta^{h}\|_{\infty,\mathcal{V}}
    \le&   \|\textbf{r}^h_{\theta,\eta}-1\|_\infty (b+d)\|Z_\eta\|_{\mathcal{V}} \cdot( \|\tilde{\textbf{g}}_{\eta}^{h}\|_{\infty,\mathcal{V}}+(\mu_\eta^h)^T\tilde{\textbf{g}}_\eta^h),
\end{align*}
 where $\|Z_\eta\|_\mathcal{V}<\infty$ from Theorem 16.1.2 in \cite{meyn2012markov}, $\|\tilde{\textbf{g}}_{\eta}^{h}\|_{\infty,\mathcal{V}}\le 1$ since $\mathcal{V}\ge \tilde{\textbf{g}}_{\eta}^{h}$ according to Assumption~\ref{ass: drift}, and $(\mu_\eta^h)^T\tilde{\textbf{g}}_\eta^h<\infty$ according to Lemma 1 of \cite{dai2022queueing}. As a result, we have $\|(\tilde{\textbf{P}}_{\theta}^{h} - \tilde{\textbf{P}}_\eta^{h}) \mathbf{v}_\eta^{h}\|_{\infty,\mathcal{V}}=O(\|\textbf{r}^h_{\theta,\eta}-1\|_\infty)$.

Next, we try to bound $\|\tilde{\textbf{g}}_{\theta}^{h} - \tilde{\textbf{g}}_{\eta}^{h}\|_{\infty,\mathcal{V}}$. Recall that the expected one-day cost vectors are: 
\begin{equation}
    \begin{aligned}
        & \tilde{\textbf{g}}_{\theta}^{h}=\textbf{g}_\theta^h+\textbf{P}_\theta^{h,h+1}\textbf{g}_\eta^{h+1}
        +
        (\textbf{P}_\theta^{h,h+1}\textbf{P}_\eta^{h+1,h+2})\textbf{g}_\eta^{h+2}
        +\cdots +(\textbf{P}_\theta^{h,h+1}\textbf{P}_\eta^{h+1,h+2}\cdots\textbf{P}_\eta^{h-2,h-1})\textbf{g}_\eta^{h-1}\\
        & \tilde{\textbf{g}}_{\eta}^{h}=\textbf{g}_\eta^h+\textbf{P}_\eta^{h,h+1}\textbf{g}_\eta^{h+1}
        +
        (\textbf{P}_\eta^{h,h+1}\textbf{P}_\eta^{h+1,h+2})\textbf{g}_\eta^{h+2}
        +\cdots +(\textbf{P}_\eta^{h,h+1}\textbf{P}_\eta^{h+1,h+2}\cdots\textbf{P}_\eta^{h-2,h-1})\textbf{g}_\eta^{h-1}.
    \end{aligned}
\end{equation}
We denote the elements of the one-day cost vector $\tilde{\textbf{g}}_\eta^h$ and $\tilde{\textbf{g}}_{\theta}^h$ as $\{\tilde g_\eta^h(s), s\in\mathcal{S}^h\}$ and $\{\tilde{g}_{\theta}^h(s), s\in\mathcal{S}^h\}$, respectively. 
Moreover, we denote
\begin{align*}
    \tilde g^h(s,f)=&g^h(s,f)+\sum_{s^{h+1}\in\mathcal{S}^{h+1}}p^{h,h+1}(s^{h+1}|s,f) \Big( g_\eta^{h+1}(s^{h+1})+\sum_{s^{h+1}\in\mathcal{S}^{h+1},s^{h+2}\in\mathcal{S}^{h+2}}p_\eta^{h+1,h+2}(s^{h+2}|s^{h+1})g_\eta^{h+2}(s^{h+2})\\
&+\cdots+\sum_{s^{h+1}\in\mathcal{S}^{h+1},\cdots,s^{h-1}\in\mathcal{S}^{h-1}}p_\eta^{h+1,h+2}(s^{h+2}|s^{h+1})\cdots p_\eta^{h-2,h-1}(s^{h-1}|s^{h-2})g_\eta^{h-1}(s^{h-1}) 
\Big) 
\end{align*}
which was also used in Equation~\eqref{equ: def adv one-day h}. Following the same argument as for the transition probabilities, this term is independent of the new policy parameter $\theta$, 
which implies that
\begin{align*}
    \tilde g_\theta^h(s)
    =&\sum_{f\in\mathcal{A}(s)}\pi_\theta(f\mid s)\tilde g_\eta^h(s,f),
\quad 
    \tilde g_\eta^h(s)
    = \sum_{f\in\mathcal{A}(s)}\pi_\eta(f\mid s)\tilde g_\eta^h(s,f) . 
\end{align*}
Then, the term $\|\tilde{\textbf{g}}_{\theta}^{h} - \tilde{\textbf{g}}_{\eta}^{h}\|_{\infty,\mathcal{V}}$ can be bounded as \begin{equation}
    \begin{aligned}
    \|\tilde{\textbf{g}}_{\theta}^{h} - \tilde{\textbf{g}}_{\eta}^{h}\|_{\infty,\mathcal{V}}=&\sup_{s\in\mathcal{S}^h}\frac{|\tilde g_{\theta}^h(s)-\tilde g_\eta^h(s)|}{\mathcal{V}(s)}\\
    =&\sup_{s\in\mathcal{S}^h}\frac{|\sum_{f\in\mathcal{A}(s)}\pi_\theta(f\mid s)\tilde g_\eta^h(s,f)-\sum_{f\in\mathcal{A}(s)}\pi_\eta(f\mid s)\tilde g_\eta^h(s,f)|}{\mathcal{V}(s)}\\
    =&\sup_{s\in\mathcal{S}^h}\frac{|\sum_{f\in\mathcal{A}(s)}(r_{\theta,\eta}(f\mid s)-1)\pi_\eta(f\mid s)\tilde g_\eta^h(s,f)|}{\mathcal{V}(s)}\\
    \le&\sup_{s\in\mathcal{S}^h}\frac{\|\textbf{r}^h_{\theta,\eta}-1\|_\infty |\sum_{f\in\mathcal{A}(s)}\pi_\eta(f\mid s)\tilde g_\eta^h(s,f)|}{\mathcal{V}(s)}\\
    =&\|\textbf{r}^h_{\theta,\eta}-1\|_\infty\|\tilde{\textbf{g}}_{\eta}^{h}\|_{\infty,\mathcal{V}}\\
    \le& \|\textbf{r}^h_{\theta,\eta}-1\|_\infty,
\end{aligned}
\label{equ: bound cost gap}
\end{equation}
which implies $\|\tilde{\textbf{g}}_{\theta}^{h} - \tilde{\textbf{g}}_{\eta}^{h}\|_{\infty,\mathcal{V}}=O(\|\textbf{r}^h_{\theta,\eta}-1\|_\infty)$ as well.
 The second-to-last inequality holds since in our setting $\tilde g_\eta^h(s,f)\ge 0,\forall (s,f)$, and the last inequality holds because $\mathcal{V}\ge \tilde{\textbf{g}}_{\eta}^{h}$ implies that $\|\tilde{\textbf{g}}_{\eta}^{h}\|_{\infty,\mathcal{V}}\le 1$.

By far, we have shown that both terms in the bound for $\|\tilde N^h(\theta,\eta)\|_{\infty,\mathcal{V}}$ has the same order $O(\|\textbf{r}^h_{\theta,\eta} - 1\|_\infty)$, so we also have $\|\tilde N^h(\theta,\eta)\|=O(\|\textbf{r}^h_{\theta,\eta} - 1\|_\infty)$.
As a result, according to~\eqref{equ: bounds for N1 N2} and the fact that  $\|\mu_{\theta}^{h}-\mu_\eta^{h}\|_{1,\mathcal{V}}
=O(\|\tilde{\textbf{P}}_{\theta}^{h}-\tilde{\textbf{P}}_\eta^{h}\|_{\mathcal{V}}
\|Z_\eta\|_{\mathcal{V}})=O(\|\textbf{r}^h_{\theta,\eta} - 1\|_\infty)$, we have $N_1(\theta,\eta)=O(\|\textbf{r}^h_{\theta,\eta} - 1\|_\infty)$, and  $N_2(\theta,\eta)=O(\|\textbf{r}^h_{\theta,\eta} - 1\|_\infty^2)$, which completes the proof.
\unskip\nobreak\hfill $\square$
\endproof  

\section{Illustration of PPO in Two-pool Setting}
\label{App: ppo explanability}
In this section, we present an illustration of the mechanism behind PPO in our specific context -- the overflow assignment for inpatients. For illustration purpose, we focus on a simple two-pool midnight model with randomized atomic action. Furthermore, we focus on illustrating the updates for one given state $s$ and assume the policy at other states remain unchanged. That is, the objective in this showcase example is to minimize 
\begin{equation}
\hat N_1(\theta,s) := \mathop{\mathbb{E}}\limits_{f\sim\pi_\eta(\cdot|s)}    [r_{\theta,\eta}(f\mid s)  \hat A_\eta(s,f)]=\mathop{\mathbb{E}}\limits_{f\sim\pi_\theta(\cdot|s)} [ \hat A_\eta(s,f)]
\label{eq:policy-gradient-appendix}
\end{equation}
for the given state $s$. The clipping function can be easily added to this. By analyzing the gradient of $\hat N_1(\theta,s)$, we showcase how the overflow policy will change with different model parameters $(B,C,\mu,\lambda)$, 
providing some \emph{explainability} of the mechanism behind PPO.

We consider a simple two-pool midnight MDP with one decision epoch each day $(m=1)$. The state is simplified as $s=(x_1,x_2)\in\mathbb{R}^2$, since we do not need to track the to-depart counts when $m=1$. Correspondingly, the transition dynamics from current state to the state of the next day given overflow action $f=\{f_{i,j}\}$ can be specified as 
\begin{align*}
    x_j'=x_j+a_j-d_j+\sum_{i=1,i\ne j}^2 f_{i,j}-\sum_{\ell=1,\ell\ne j}^2 f_{j,\ell},\quad j=1,2.
\end{align*}
where $a_j,d_j$ denote the number of new arrivals and departures within a day. Here, $a_j$ is a realization of the random variable $A_j$ which follows $\text{Poisson}(\Lambda_j)$, and $d_j$ is a realization of the random variable $D_j$ which follows $\text{Bin}(q_j, \mu_j)$. 

\subsection{Policy Gradient}
\label{app-subsec: policy gradient assumptions}

To facilitate the gradient analysis, we make additional assumptions. 

\begin{assumption} (Symmetric Two-pool Midnight MDP)
\label{ass: model} The two pools have
  $(N_j, \lambda_j, \mu_j) = (N, \lambda, \mu)$ for $j=1,2$;  $C_1=C_2=C$, $B_{12}=B_{21}=B$.
\end{assumption}
Under Assumption~\ref{ass: model}, 
we can define two state subspaces according to the feasible actions:
    \begin{align*}
        \mathcal{S}_1=\{(x_1,x_2)\in\mathcal{S}: x_1> N, x_2< N\}; \quad 
        \mathcal{S}_2=\{(x_1,x_2)\in\mathcal{S}: x_1< N, x_2> N\} . 
    \end{align*}
Recall that the system-level action takes the form $f=\{f_{i,j}, i,j=1,2\}$, where $f_{i,j}$ represents the number of assignments from class $i$ to pool $j$. According to the definition of feasible action defined in Equation~\eqref{Equ: feasibility of action}, for $s\in\mathcal{S}_1$, the feasible action space 
is $\{\{q_1-f_{1,2},f_{1,2}, 0, 0\}: f_{1,2}=0,1,\dots,\min\{q_1,N-x_2\}\}$, where $q_1=(x_1-N)\vee 0$ denotes the queue length of class $1$. Similarly, for $s\in\mathcal{S}_2$, the feasible action space is $\{\{0, 0, f_{2,1},q_2-f_{2,1}\}: f_{2,1}=0,1,\dots,\min\{q_2,N-x_1\}\}$, where $q_2=(x_2-N)\vee 0$ denotes the queue length of class $2$. For any state $s\in\mathcal{S}\setminus(\mathcal{S}_1\cup\mathcal{S}_2)$, the only feasible action is no-overflow (action $\{0,0,0,0\}$). 
Without loss of generality, we focus on $\mathcal{S}_1$ in the following analysis, as the results can be easily extend to $\mathcal{S}_2$ due to symmetry in Assumption~\ref{ass: model}.

\begin{assumption} (Parametric Randomized Atomic Action)
\label{ass: policy}
    \begin{itemize}
\item [(i)] Batched setting: For a given pre-action state $s$, each atomic action $a^n$ depends on $s$, i.e., not affected by the previous atomic assignment. 

\item [(ii)] Parametric logistic model: The routing probability for the atomic action of a class 1 customer is parameterized as
$$\kappa_\theta(1|s,1)=\frac{1}{1+\exp(\theta_1 x_1+\theta_2 x_2+\theta_0)}
, \quad \kappa_\theta(2|s,1) = 1-\kappa_\theta(1|s,1).
$$ 
\end{itemize}
\end{assumption}
Under a randomized policy $\pi_\theta$ satisfying Assumption~\ref{ass: policy}, for a given pre-action state $s\in\mathcal{S}_1$ and an associated feasible action $f=(q_1-f_{1,2},f_{1,2},0,0)$, the overflow quantity $f_{1,2}$ follows a binomial distribution $Bin(q_1,\kappa_\theta(2|s,1))$ (note that we allow overflow assignments to a full server here). Therefore, the aim of PPO is to update the parameters $\theta=(\theta_0,\theta_1,\theta_2)'\in\mathbb{R}^3$ to minimize $\hat N_1(\theta,s)=\mathbb{E}_{f\sim\pi_\theta(\cdot|s)}[ \hat A_\eta(s,f)]$ in~\eqref{eq:policy-gradient-appendix} through the policy gradient approach. 

\begin{assumption} (Advantage function approximation)
\label{ass: adv}
    \begin{itemize}
        \item [(i)] Linear approximation: The value function $v_\eta$ is approximated with linear combinations of a set of linear and quadratic basis functions, i.e., $$\hat v_\eta(s)=\hat\beta_1x_1+\hat\beta_2x_2+\hat\beta_3 x_1^2+\hat\beta_4x_2^2,$$
        where $\hat\beta_k,k=1,2,3,4$ are the coefficient parameters. 
        \item [(ii)] Transition probability approximation: The waiting customers in buffers can be served with the same service time distribution as in server pools, and leave the system after service.
    \end{itemize}
\end{assumption}
In Section~\ref{Subsec: policy evaluation} of the main paper, we approximate the relative value function with the linear combinations of a set linear/quadratic basis as well as some queueing-based basis. Here for simplicity, in Assumption~\ref{ass: adv}(i), we focus on the simpler linear and quadratic basis. Assumption~\ref{ass: adv}(ii) is made to simplify the evaluation of the advantage function $\hat A_\eta$, which can be computed according to
\begin{equation}
    \hat A_\eta(s, f)=g(s,f)-\gamma_\eta+\mathbb{E}_{s'\sim p(\cdot|s,f)}[\hat v_\eta (s')]-\hat v_\eta(s).
    \label{equ: adv-general}
\end{equation}

\subsection{Policy gradient}
\label{app-subsec: policy gradient computation}
We state the policy gradient result in the following lemma, with its detailed proof in Section~\ref{app:proof-lemma2}.

\begin{lemma} 
Under Assumptions~\ref{ass: model}-\ref{ass: adv}, for any $s=(x_1,x_2)\in\mathcal{S}_1$ and $f=(q_1-f_{1,2},f_{1,2}, 0, 0)$,
    \begin{align*}
        \frac{\partial\hat N_1(\theta,s)}{\partial \theta_0}
    &= \nabla_0\hat N_1(\theta,s), 
    \\
    \frac{\partial \hat N_1(\theta)}{\partial \theta_k}
    &=\nabla_0\hat N_1(\theta,s)\cdot x_k, \quad k=1,2 . 
    \end{align*}
Here, 
    \begin{equation}
        \begin{aligned}
        \nabla_0\hat N_1(\theta,s)&=\sum_{f_{1,2}=0}^{q_1} \pi_\theta(f\mid s)\big( f_{1,2} - q_1\kappa_\theta(2|s,1) \big) \hat A_\eta(s,f)\\
        &=q_1\kappa_{\theta}(2|s,1)\big(1- \kappa_{\theta}(2|s,1)\big) \Big( 2\hat\beta_3(1-\mu)^2\big(2(q_1-1)\kappa_{\theta}(2|s,1)+x_2-x_1+1\big)+B-C\Big).
    \end{aligned}
    \label{equ: policy gradient with adv}
    \end{equation}
    \label{prop: policy gradient}
    \end{lemma} 

This closed form for policy gradient  allows us to examine the optimal action that minimizes $\hat N_1(\theta,s)$. Through this examination, we generate insights into how the policy gradient approach is guiding us to find a good action under different model and cost parameters $(\lambda, \mu, B, C)$.

We start by analyzing the monotonicity of $\nabla_0 \hat N_1(\theta,s)$ w.r.t. $\kappa_\theta(2|s,1)$, which depends on the sign of $\hat\beta_3$. In the rest of the analysis, we focus on the case where $\hat\beta_3>0$ since it leads to non-trivial policy updates. In this case, the new overflow probability obtained by minimizing $\hat N_1(\theta,s)$ should either equals to 0 or 1, or make $\nabla_0 \hat N_1(\theta^*,s)=0$ hold. The latter first-order condition gives us $\kappa_{\theta^*}(2|s,1)=\max\big(0,\min(1,\kappa^*(s)) \big)$, where 
$$
\kappa^*(s)=\frac{N-x_2-(B-C)/2\hat\beta_3(1-\mu)^2}{2(q_1-1)}+\frac{1}{2}.
$$ 
We discuss the property of $\kappa^*(s)$ separately when $B\leq C$ or $B>C$.

When $B\leq C$, $\kappa^*(s)\ge \frac{N-x_2+q_1-1}{2(q_1-1)}\ge \frac{N-x_2}{q_1}$. If $\kappa_\theta(2|s,1)=\frac{N-x_2}{q_1}$, the expected number of overflow assignments equals to the number of idle servers in pool 2. This essentially corresponds to the complete-overflow policy, which is expected to be optimal when the overflow cost is cheap. 

When $B>C$, $\kappa^*(s)$ decreases in $x_2$, which follows the intuition about a ``good'' policy, i.e., the more crowded pool 2 gets, the less overflow should be assigned from class 1 t pool 2. For how this action changes with $x_1$, we focus on examining the mean of overflow assignment, i.e., $q_1\kappa^*(s)$. We have 
$$
q_1\kappa^*(s) =\frac{N-x_2-(B-C)/2\hat\beta_3(1-\mu)^2+1}{2}+\frac{N-x_2-(B-C)/2\hat\beta_3(1-\mu)^2}{q_1-1}+\frac{1}{2}(q_1-1), 
$$ 
which increase with $q_1$ when $q_1\ge 2(N-x_2)-\frac{B-C}{\hat\beta_3(1-\mu^2)}$, but decrease with $q_1$ otherwise. Therefore, when $x_2$ is close to $N$, the critical point $2(N-x_2)-\frac{B-C}{\hat\beta_3(1-\mu^2)}\le 1$, so $q_1\kappa^*(s) $ increase with $q_1$, which also follows the intuition about a ``good'' policy since we need to overflow more to balance the load when there are more waiting patients. However, when $x_2$ is small, which means there are enough idle beds, the mean value of overflow assignments firstly decrease then increase with $x_1$. This policy is desired since when $x_1$ is very large, load balancing is the first-order issue, so we overflow more when $x_1$ is large; in contrast, when $x_1$ is relatively small, we need to trade-off between holding cost and undesirable overflow assignments, so when $x_1$ is larger, even with the same mean value of overflow assignments, there is a larger possibility that it will conduct a large number of overflow assignments and occupy too much class 2 servers in that case, causing a very large future cost according to ``snowball effect''. Therefore, mean value of overflow assignments should decrease.

In addition, a critical term in $q_1\kappa^*(s)$ is the term $\frac{(B-C)}{2\hat\beta_3(1-\mu)^2}$. Through some argument, we can show that when $B>C$, $\alpha$ in with $B-C$ and $\mu$. It implies that the willingness of overflow decrease with $B-C$ and $\mu$. These results also follow our intuition because when overflow cost is closer to holding cost (the gap $(B-C)$ is smaller), we prefer to overflow more to help balance the system; when the busy servers completes jobs faster (larger $\mu$), we prefer to let customers wait since they can be admitted into primary ward within a shorter time.

\subsection{Proof of Lemma~\ref{prop: policy gradient}}
\label{app:proof-lemma2}

For a given $s=(x_1,x_2)\in\mathcal{S}_1$, a feasible action $f$ takes the form of $(q_1-f_{1,2},f_{1,2}, 0, 0)$. 
Under the assumptions in Section~\ref{app-subsec: policy gradient assumptions}, we get from~\eqref{eq:policy-gradient-appendix} that $N_1(\theta,s) = \mathbb{E}_{f\sim\pi_\theta(\cdot|s)}[ \hat A_\eta(s,f)] = \sum_{f_{1,2}=0}^{q_1} \pi_\theta(f\mid s)\hat A_\eta(s,f)$. 
Therefore, taking derivative of $\hat N_1(\theta,s)$ w.r.t. $\theta_k,k=0,1,2$, respectively, we get 
  \begin{align}
      \frac{\partial }{\partial \theta_k}\hat N_1(\theta,s)
      &=\sum_{f_{1,2}=0}^{q_1} \frac{\partial \pi_\theta(f\mid s)}{\partial \theta_k}\hat A_\eta(s,f).  \label{Equ: gradient N 3}
  \end{align}
From Assumption~\ref{ass: policy}(i), the number of overflow quantity from class $1$ to pool $2$, $f_{1,2}$ follows $Bin(q_1,\kappa_\theta(2|s,1))$ under policy $\pi_\theta$, we can rewrite $\pi_\theta(f\mid s)$ as
  \begin{equation}
      \pi_\theta(f\mid s)=\left(\begin{array}{c}
           q_1  \\
           f_{1,2} 
      \end{array}\right) \kappa_\theta(2|s,1)^{f_{1,2}}(1-\kappa_\theta(2|s,1))^{q_1-f_{1,2}}.
      \label{Equ: action probability, binomial}
  \end{equation}
Then, by using some algebra, we have
\begin{align}
\frac{\partial }{\partial \theta_k}\pi_\theta(f\mid s) 
=& \pi_\theta(f\mid s)\left(\frac{f_{1,2}}{\kappa_\theta(2|s,1)}-\frac{q_1-f_{1,2}}{1-\kappa_\theta(2|s,1)}\right)\frac{\partial\kappa_\theta(2|s,1)}{\partial\theta_k}, 
\quad k=0,1,2. 
\label{Equ: action probability gradient 2}
\end{align}
Furthermore, recall that from Assumption~\ref{ass: policy}(ii), $\kappa_\theta(2|s,1) $ is parameterized as a logistic function. 
Therefore, we can further write out the following form for the gradients of $\kappa_\theta(2|s,1)$. For the gradient w.r.t. $\theta_0$, we have
\begin{equation}
    \begin{aligned}
    \frac{\partial \kappa_\theta(2|s,1)}{\partial \theta_0}=&-\frac{\text{exp}(-(\theta_1 x_1+\theta_2 x_2+\theta_0))\cdot(-1)}{(1+\text{exp}(-(\theta_1 x_1+\theta_2 x_2+\theta_0)))^2}\\
    =&
    \kappa_\theta(2|s,1)(1-\kappa_\theta(2|s,1)).
\end{aligned}
\label{Equ: atomic probability gradient 0}
\end{equation}
Similarly, for $\theta_1, \theta_2$, we get 
\begin{equation}
    \frac{\partial \kappa_\theta(2|s,1)}{\partial \theta_k}=\kappa_\theta(2|s,1)(1-\kappa_\theta(2|s,1)) x_k.
    \label{Equ: atomic probability gradient 1}
\end{equation}
Combining Equations~\eqref{Equ: action probability gradient 2} through~\eqref{Equ: atomic probability gradient 1} and plugging them back to~\eqref{Equ: gradient N 3}, we get the final results of policy gradient as follows.
\begin{align*}
\frac{\partial \hat N_1(\theta,s)}{\partial\theta_0}=&\sum_{f_{1,2}=0}^{q_1}\pi_\theta(f\mid s)\left(f_{1,2}-q_1\kappa_\theta(2|s,1)\right)\hat A_\eta(s,f), 
\\
\frac{\partial \hat N_1(\theta,s)}{\partial\theta_k}=&\sum_{f_{1,2}=0}^{q_1}\pi_\theta(f\mid s)\left(f_{1,2}-q_1\kappa_\theta(2|s,1)\right)x_k\hat A_\eta(s,f),\quad k=1,2. 
\end{align*}
For simplicity, we use $\nabla_0\hat N_1(\theta,s)$ to denote
\begin{equation}
    \sum_{f_{1,2}=0}^{q_1}\pi_\theta(f\mid s)\left(f_{1,2}-q_1\kappa_\theta(2|s,1)\right)x_k\hat A_\eta(s,f).
    \label{equ: tmp N0}
\end{equation}
As a result, the policy gradient can be rewritten as
  \begin{align*}
        \frac{\partial\hat N_1(\theta,s)}{\partial \theta_0}
    &= \nabla_0\hat N_1(\theta,s), 
    \\
    \frac{\partial \hat N_1(\theta)}{\partial \theta_k}
    &=\nabla_0\hat N_1(\theta,s)\cdot x_k, \quad k=1,2 . 
    \end{align*}
Next, to derive the closed form of the policy gradient, we need to derive the closed form of $\hat A_\eta$ and plug it into~\eqref{equ: tmp N0}. Recall that given a pre-action state $s\in\mathcal{S}_1$ and a feasible action $f=(f_{1,2}, q_1-f_{1,2}, 0, 0)$ with $0\le f_{1,2}\le q_1$, the advantage function $\hat A(s,f)$ can be computed via
\begin{equation}
    \hat A_\eta(s,f)=g(s,f)+\mathbb{E}_{s'\sim p(\cdot|s,f)}[\hat v_\eta(s')],
    \label{equ: est adv}
\end{equation}
where the current cost follows 
$$g(s,f)=C (q_1-f_{1,2})+B f_{1,2},$$ and according to Assumption~\ref{ass: adv}, the estimated value function follows 
 $$\hat v_\eta(s)=\hat\beta_1x_1+\hat\beta_2x_2+\hat\beta_3 x_1^2+\hat\beta_4x_2^2.$$
According to Assumption~\ref{ass: model}(i), the two-pool system is symmetric, so the parameters $\{\hat\beta_i,i=1,...,4\}$ for estimating $\hat v_\eta$ should also be symmetric, i.e., 
$$\hat v_\eta=\hat\beta_1(x_1+x_2)+\hat\beta_3(x_1^2+x_2^2).$$
To compute the closed form of cost-to-go $\mathbb{E}_{s'}[\hat v_\eta(s')]$, we need to specify the transition dynamics in our simplified two-pool setting. That is, given $(s,f)$, the next state $s'=(x_1',x_2')$ follows 
\begin{equation}
    x_1'=x_1-f_{1,2}+A_1-D_1,\quad
    x_2'=x_2+f_{1,2}+A_2-D_2,
\end{equation}
where from Assumption~\ref{ass: model}(i)(ii), the number of new arrivals $A_1,A_2$ both follow Poisson distribution with parameter $\lambda$, and the number of new departures $D_1,D_2$ follow distributions $Bin(x_1-f_{1,2},\mu)$ and $Bin(x_2+f_{1,2},\mu)$, respectively. Therefore, we have
\begin{equation}
    \begin{aligned}
  &\mathbb{E}_{s'\sim p(\cdot|s,f)}[\hat v_\eta (s')]
  \\
  = &  \mathbb{E}_{s'\sim p(\cdot|s,f)}[\hat \beta_1 (x_1'+x_2')+\hat\beta_3((x_1')^2+(x_2')^2)]
  \\
 =&\hat\beta_1\big(x_1+x_2+2\lambda-(x_1+x_2)\mu\big)
 +\hat\beta_3\mathbb{E}\big[ (x_1-f_{1,2}+A_1-D_1)^2+(x_2+f_{1,2}+A_2-D_2)^2 \big] .  
\end{aligned}
\label{Equ: cost-to-go function}
\end{equation}

Via some algebra to evaluate the expectation term in~\eqref{Equ: cost-to-go function}, we have 
\begin{align*}
\mathbb{E}_{s'\sim p(\cdot|s,f)}[\hat v_\eta (s')] 
=&\hat\beta^3(1-\mu)^2[(x_1-f_{1,2})^2+(x_2+f_{1,2})^2]+(\hat\beta_1+\hat\beta_3(2\lambda-\mu))(1-\mu)(x_1+x_2)+2\hat\beta_1\lambda+2\hat\beta_3(\lambda+\lambda^2)
\end{align*}
Plugging the formulas of $g(s,f)$ and $\mathbb{E}[\hat v_\eta (s')]$ back into~\eqref{equ: est adv}, we get 
\begin{align}
     \hat A_\eta(s, f)=&g(s,f)-\gamma+\mathbb{E}_{s'\sim p(\cdot|s,f)}[\hat v_\eta (s')]-\hat v_\eta(s)\notag
     \\
     =&(B-C)f_{1,2}+\hat\beta_3(1-\mu)^2[2f_{1,2}^2-2(x_1-x_2)f_{1,2}]+Const(s),  \label{equ: adv closed form}
\end{align}
where $Const(s)$ is a constant that depends on $s$ but is independent of $f$. Finally, by plugging~\eqref{equ: adv closed form} into~\eqref{equ: tmp N0}, we can rewrite the policy gradient $\nabla_0 \hat N_1(\theta,s)$ as 
\begin{align*}
    \nabla_0 \hat N_1(\theta,s)=&\sum_{ f_{1,2}=1}^{q_1}\pi_{\theta}( f\mid s)\left(f_{1,2}-q_1\kappa_{\theta}(2|s,1)\right)\hat A_\eta(s, f)\\
    =&\sum_{ f_{1,2}=1}^{q_1}\pi_{\theta}( f\mid s)\left(f_{1,2}-q_1\kappa_{\theta}(2|s,1)\right)\left((B-C)f_{1,2}+\hat\beta_3(1-\mu)^2[2f_{1,2}^2-2(x_1-x_2)f_{1,2}]+Const(s)\right)\\
   =&q_1\kappa_{\theta}(2|s,1)\big(1- \kappa_{\theta}(2|s,1)\big) \Big( 2\hat\beta_3(1-\mu)^2\big(2(q_1-1)\kappa_{\theta}(2|s,1)+x_2-x_1+1\big)+B-C\Big). 
\end{align*}
Here, we have used the binomial distribution property for $f_{1,2}$ and we are able to eliminate the $Const(s)$ since 
$$\sum_{ f_{1,2}=1}^{q_1}\pi_{\theta}( f\mid s)\left(f_{1,2}-q_1\kappa_{\theta}(2|s,1)\right)Const(s)=(\mathbb{E}[f_{1,2}]-q_1\kappa_\theta(2|s,1))\cdot Const(s)=0.
$$
\unskip\nobreak\hfill $\square$
\endproof

\end{APPENDICES}

\end{document}